\documentclass{tran-l}

\newtheorem{theorem}{Theorem}
\newtheorem{lemma}[theorem]{Lemma}

\theoremstyle{definition}
\newtheorem{definition}[theorem]{Definition}
\newtheorem{example}[theorem]{Example}
\newtheorem{proposition}[theorem]{Proposition}

\newtheorem{xca}[theorem]{Exercise}

\theoremstyle{remark}
\newtheorem{remark}{Remark}

\numberwithin{equation}{section}

\newcommand{\abs}[1]{\lvert#1\rvert}

\newcommand{\blankbox}[2]{%
  \parbox{\columnwidth}{\centering
    \setlength{\fboxsep}{0pt}%
    \fbox{\raisebox{0pt}[#2]{\hspace{#1}}}%
  }%
}

\def\1{\mbox{1\hspace{-0.25em}l}}

\def\gH{{{\mathbf H}}}
\def\cI{{{\mathcal L}}}
\def\cI{{{\mathcal I}}}

\newcommand \A[1]{{\bf (#1)}}

\def\leftB{[\![}
\def\rightB{]\!]}
\def\btheta{{\boldsymbol{\theta}}}
\def\bxi{{\boldsymbol{\xi}}}
\def\bXi{{\boldsymbol{\Xi}}}

\def\X{{\mathbf{X}}}

\def\x{{\mathbf{x}}}

\def\y{{\mathbf{y}}}
\def\z{{\mathbf{z}}}

\def\0{{\mathbf{0}}}

\def\gF{{\mathbf{F}}}

\def\K{{\mathbf{K}}}
\def\gR{{\mathbf{R}}}

\newcommand{\F}{\mathcal{F}}  
\renewcommand{\P}{\mathbb{P}}

\newcommand{\E}{\mathbb{E}}

\newcommand{\N}{\mathbb{N}}  
  
\newcommand{\T}{\mathbb{T}}  

\newcommand{\mW}{\mathcal{W}}  

\newcommand{\R}{\mathbb{R}}  
\newcommand{\I}{\mathbb{I}}

\def\bphi{{\boldsymbol{\phi}}}

\usepackage{color}

\begin{document}

\title[Regularization by noise for a chain of ODEs]{Regularization effects of a noise propagating through a chain of differential equations: an almost sharp result}

\author{Paul-Eric Chaudru de Raynal}
\address{UMR CNRS 6629, Laboratoire de Mathématiques Jean Leray, 2, rue de la Houssini\`ere
44322 Nantes cedex 3, France.}
\curraddr{}
\email{Paul-Eric.Deraynal@univ-nantes.fr}

\author{St\'ephane Menozzi}
\address{UMR CNRS 8071, Laboratoire de Mod\'elisation Math\'ematique d'Evry (LaMME), Universit\'e d'Evry Val d'Essonne, Paris-Saclay, 23 Boulevard de France 91037 Evry, France and Laboratory of Stochastic Analysis, HSE University,
Pokrovsky Blvd, 11, Moscow, Russian Federation}
\email{stephane.menozzi@univ-evry.fr}
\thanks{The article  was prepared within the framework of the HSE University Basic Research Program}


\date{October 9, 2017 and, in revised form, May 1, 2021.}



\begin{abstract}
We investigate the effects of the propagation of a non-degenerate Brownian noise through a chain of deterministic differential equations whose coefficients are \textit{rough} and satisfy a weak like H\"ormander structure (i.e. a non-degeneracy condition w.r.t. the components which transmit the noise).
In particular  we characterize, through suitable counter-examples, almost sharp regularity exponents that ensure that weak well posedness  holds for the associated SDE. As a by-product of our approach, we also derive some density estimates of Krylov type for the weak solutions of the considered SDEs.
\end{abstract}

\maketitle

\section{Introduction and Main Results}
In this work we are interested in studying the weak regularization effects of a Brownian noise propagating through a chain of $n$ $d$-dimensional oscillators. Namely, we establish weak uniqueness for Stochastic Differential Equations (SDEs in short) of the following type:

\begin{equation}
\label{SYST}
\begin{array}{l}
\displaystyle dX_t^1 = F_1(t,X_t^1,\dots,X_t^n) dt + \sigma(t,X_{t}^1,\dots,X_{t}^n) dW_t,
\\
\displaystyle dX_t^2 = F_2(t,X_t^1,\dots,X_t^n) dt,
\\
\displaystyle dX_t^3 = F_3(t,X_t^2,\dots,X_t^n) dt,
\\
\vdots
\\
\displaystyle dX_t^n = F_n(t,X_t^{n-1},X_t^n) dt \textcolor{black}{.}
\end{array}
\quad t \geq 0.
\end{equation}
In the above equation, $(W_t)_{t \geq 0}$ stands for a $d$-dimensional Brownian motion and the components $(X_t^i)_{i\in \leftB1,n\rightB}$ are $\R^d $-valued as well. We suppose that the $(F_i)_{i\in \leftB 2,n\rightB}$ satisfy a kind of weak H\"ormander condition, i.e. the matrices $\big(D_{x_{i-1}}F_i(t,\cdot)\big)_{i\in \leftB 2,n\rightB} $ have full rank. However, the coefficients $ (F_i)_{i\in \leftB 2,n\rightB}$ can be rather \textit{rough} in their other entries. Namely, H\"older continuous or even  in a suitable $L^q-L^p $ space for $F_1$, \textcolor{black}{where the parameter $q$ relates to the time integrability and $p$ to the spatial one}. 
We assume as well that the diffusion coefficient $\sigma $ is bounded \textcolor{black}{from above and below} and spatially H\"older continuous.

We emphasize that, under these conditions, the Stroock and Varadhan Theory \textcolor{black}{for weak uniqueness} does not apply. This especially comes from the \textcolor{black}{specific degenerate framework} 
 \textcolor{black}{considered} here: the noise in the $i^{\rm th}$ component only comes from the $(i-1)^{\rm th}$ component, $2 \leq i \leq n$, \textcolor{black}{through the non-degeneracy of the gradients $\big(D_{x_{i-1}}F_i(t,\cdot)\big)_{i\in \leftB 2,n\rightB} $ (components which transmit the noise)}. We nevertheless show that the system is well posed, in a weak sense, when the drift of the first component is H\"older continuous or bounded in supremum norm or in suitable $L^q-L^p$ norm and the drift functions of the other components are only H\"older continuous with respect to the variable\textcolor{black}{s} that do not transmit the noise. Denoting by $(\beta_i^j)_{2\leq i \leq j \leq n}$ the H\"older index of the drift of the $i^{\rm th}$ component w.r.t. the $j^{\rm th}$ variable we assume $\beta_{i}^j\in \big( [(2i-3)/(2j-1)],1\big]$. We also show that these thresholds are (almost) sharp thanks to \textcolor{black}{appropriate} counter examples.

Also, as a by-product of our analysis, we prove that the density of the unique weak solution of the system satisfies Krylov-like estimate\textcolor{black}{s}.

\subsection*{Weak and strong regularization by noise} 

Strong and weak well posedness of stochastic systems outside the classical Cauchy-Lipschitz framework have motivated a lot of works since the last past four decades\footnote{\textcolor{black}{In the presentation below, we will mainly focus on Brownian driven SDEs. We can refer to the recent work of Priola \cite{prio:17} 
for the more general L\'evy driven case in the non-degenerate framework.}}.

Concerning the strong well posedness, the first result in that direction is due to Zvonkin \cite{zvonkin_transformation_1974} who showed that one-dimensional non degenerate \textcolor{black}{Brownian driven} stochastic differential equations with bounded and measurable drift and H\"older continuous diffusion matrix are well posed for H\"older index strictly greater than $1/2$. Then, Veretennikov \cite{veretennikov_strong_1980} generalized the result to the multidimensional case for \textcolor{black}{a} Lipschitz diffusion matrix. These results have been recently revisited in the work of Krylov and R\"ockner \cite{kryl:rock:05}, where SDE\textcolor{black}{s} with additive \textcolor{black}{Brownian} noise and locally integrable drift are shown to be \textcolor{black}{strongly} well posed and Zhang \cite{zhang_well-posedness_2010} \textcolor{black}{who} extended \textcolor{black}{the} Krylov and R\"ockner result to SDEs with multiplicative noise and weakly Lipschitz diffusion matrix (\emph{i.e.} in Sobolev Sense). 
\textcolor{black}{Similar issues are handled as well} in \cite{fedr:flan:11}. Also, we can mention the recent work by Davie \cite{davie_uniqueness_2007} in which \emph{path-by-path} uniqueness is proved for non degenerate \textcolor{black}{Brownian} SDEs with bounded drift and the approach of Catellier and Gubinelli \cite{catellier_averaging_2012} (which also relies on \emph{path-by-path} uniqueness) where SDEs with additive fractional Brownian noise are investigated. \textcolor{black}{Finally, let us mention the work \cite{gradinaru_existence_2013} where the strong well-posedness of a particular one dimensional system with singular inhomogeneous drift is proved.} We refer the reader to the Saint Flour Lecture notes of Flandoli \cite{flandoli_random_2011} where a very interesting and general account on the topic is given.

On the other hand, \textcolor{black}{and still in the Brownian framework}, it has been shown that non degenerate stochastic system\textcolor{black}{s} are well posed in a weak sense as soon as the drift function is measurable and bounded and the diffusion matrix only a continuous (in space) function. This is the celebrated theory of the \emph{martingale problemma} put on complete mathematical framework by Stroock and Varadhan, \textcolor{black}{see}  \cite{stro:vara:79}. Weak well posedness results for non degenerate SDE with additive \textcolor{black}{noise} have also been explored recently: Flandoli, Issoglio and Russo showed in \cite{flandoli_multidimensional_2014} that multidimensional non degenerate SDE\textcolor{black}{s} with non-homogeneous distributional drift are well posed as soon as the regularity index is strictly greater than $-1/2$. At the same time, Delarue and Diel proved in \cite{delarue_rough_2015} that the result still holds when the regularity index is strictly greater than $-2/3$ in the one-dimensional case. This last work has then been generalized by Cannizzaro and Chouk \cite{cannizzaro_multidimensional_2015} to the multidimensional \textcolor{black}{setting}. \textcolor{black}{Note however that, in the two last mentioned works, the Authors assumed that the drift can be enhanced into a rough path structure.}

All the above strong and weak results \textcolor{black}{deeply} rely on the non-degeneracy assumption imposed to the noise and \textcolor{black}{illustrate what is usually called, following the terminology of Flandoli, a \emph{regularization by noise phenomenon}}. Here, the regularization has to be \textcolor{black}{understood} as follows: while \textcolor{black}{an} ordinary differential system could be ill-posed when the drift is less than Lipschitz (or weakly Lipschitz \cite{dipe:lion:89}), the analogous non degenerate stochastic system is well posed (in a strong or a weak sense). To obtain this kind of result, the noise plays a central role. A striking example to understand the phenomenon is the Peano Example : while the deterministic \textcolor{black}{scalar} ODE 
\begin{equation}\label{eq:peano}
\dot{Y}_t = {\rm sign}(Y_t)|Y_t|^\alpha dt,\quad Y_0=0,\quad  \alpha\in (0,1),
\end{equation}
has an infinite number of solutions that could still \textcolor{black}{be} trapped in the singularity for any given time, the corresponding \textit{\textcolor{black}{Brownian}} SDE is strongly well posed. In \cite{delarue_transition_2014}, Delarue and Flandoli put the phenomenom in light: in short time, the fluctuations of the noise dominate the system so that the solution leaves the singularity and in long time, the drift dominates and constrains the solution to fluctuate around one of the extrem\textcolor{black}{al} solution\textcolor{black}{s} of the Peano Example. Hence, there is a strong competition, in short time, between the irregularity of the drift and the fluctuations of the noise.\\

Here, our result 
\textcolor{black}{mostly emphasizes}
a \emph{regularization phenomenon coming from a degenerate noise} (\emph{i.e.} when $n\geq 2$ in \eqref{SYST}). In view of the above discussion, it is clear that the degeneracy may dramatically damage the \emph{regularization by noise} properties and, in order to preserve some regularization effect, 
\textcolor{black}{the noise still needs}
to have a way to propagate through the system. \textcolor{black}{Such kind of behavior will typically hold when the system satisfies a so-called H\"ormander condition for
hypoellipticity, see H\"ormander's seminal work on the topic \cite{horm:67}}. 

In our case, we suppose the drift of each component to be differentiable w.r.t. its first variable and the resulting gradient to be non-degenerate, but only H\"older continuous in the other variable. This last \textcolor{black}{non-degeneracy assumption} is the reason why this kind of condition is called \emph{weak H\"ormander condition}. \textcolor{black}{Namely, \textcolor{black}{up to some regularization of the diffusion coefficient,} the drift is needed to span the space through Lie Bracketing}. Also, in comparison with the \textit{\textcolor{black}{general}}  H\"ormander setting, \textcolor{black}{the specific drift structure we consider here is such that at each  \textcolor{black}{level} of the chain we only require one additional Lie bracket to generate the corresponding directions,  \textcolor{black}{up to some regularization of the diffusion coefficient again}}. This setting allows us to recover some regularization effect of the noise at each  \textcolor{black}{level} of the chain\textcolor{black}{. We also refer for similar issues to Section 2 in \cite{lan:poli:94}}.

Concerning the strong regularization effects of a degenerate noise \textcolor{black}{in a weak H\"ormander setting,}
one of the first result has been given by Chaudru de Raynal  in \cite{chau:17} and concerns strong well posedness of the above system \textcolor{black}{\eqref{SYST}} when $n=2$. It is shown in that case  that  the system is well posed as soon as the drift coefficients are H\"older continuous with H\"older exponent strictly greater than $2/3$ w.r.t. the degenerate variable and when the diffusion matrix is Lipschitz continuous in space. \textcolor{black}{This result was then extended by Wang and Zhang \cite{wang:zhan:16} under H\"older-Dini conditions with the same critical H\"older threshold $2/3$. We also mention, again for two oscillators and when the degenerate component only depends linearly of the non-degenerate variable and not on the degenerate component, the recent work by Fedrizzi, Flandoli, Priola and Vovelle \cite{fedr:flan:prio:vove:17} who address the case of  a weakly differentiable non-degenerate drift with order of weak differentiation strictly greater than $2/3$. The critical case corresponding to the exponent $2/3 $ has been discussed by Zhang \cite{zhan:16}.}

\textcolor{black}{From the weak regularization by noise viewpoint, in our current weak H\"ormander setting,  one of the first results is the work by Menozzi \cite{meno:10}. The key-point there is to exploit some smoothing effects of a suitable parametrix kernel, associated with a Gaussian linearization of \eqref{SYST}, which had already been used by Delarue and Menozzi in \cite{dela:meno:10} to derive heat-kernel bounds for the solution of \eqref{SYST}. In \cite{meno:10}, it is shown that the system \eqref{SYST} is (weakly) well posed for a spatially Lipschitz continuous drift satisfying the previously mentioned non-degeneracy condition, and a spatially H\"older continuous diffusion coefficient. The result was then extended in \cite{meno:17} for a spatially continuous diffusion coefficient}, \textcolor{black}{following the martingale problem approach establishing some suitable Calder\'on-Zygmund estimates for a degenerate Gaussian kernel and appropriate non-standard localization arguments.}  Also, in the case of two oscillators, Zhang showed in \cite{zhan:16} that when the degenerate component \textcolor{black}{only depends linearly on} the non-degenerate variable and not on the degenerate component, the system is weakly well posed as soon as the drift of the first component satisfies some local integrability conditions and when the diffusion coefficient is continuous. 
At the same time, Chaudru de Raynal showed in \cite{chau:16} that when $n=2$ the system is well posed in a weak sense as soon as the drift of each component are at least $1/3$ H\"older continuous in the degenerate variable and showed that this result is (almost) sharp for the drift of the second oscillator thanks to an appropriate counter example.

Hence, the minimal threshold obtained for the H\"older regularity of the drift is not an artefact: this can be seen as the price to pay to balance the degeneracy of the noise. Especially, in view of the previous discussion on the Delarue and Flandoli work, it \textcolor{black}{is related} to the fact that the fluctuations of the degenerate noise are not strong enough to push the solution away from the singularity if the drift is too irregular. As said above, this is illustrated in \cite{chau:16} where a counter example is built thanks to this observation. Namely, it is shown that any \textcolor{black}{stochastic} perturbation of the Peano Example \eqref{eq:peano} has to have (at least) fluctuations of order strictly 
\textcolor{black}{lower} than $1/(1-\alpha)$ in order to restore (weak) uniqueness.
\textcolor{black}{Hence, for two oscillators, assuming that the dynamics of the degenerate component is driven by \eqref{eq:peano} perturbed by the integral of the Brownian source plugged in the non-degenerate component in \eqref{SYST}, we have that
the typical variance of the noise is of order $t^{3/2}$ at time $t$. From the above condition, we indeed find $1/(1-\alpha)>3/2 \iff \alpha>1/3 $.}
\vspace*{-.2cm}
\subsection*{Organization of our paper}
The paper is organized as follows. Our assumptions and main results are stated at the end of the current Section. We present in Section \ref{SEC2} the main tools that allow to derive our results. Namely, a suitable  Gaussian linearization of the initial model \eqref{SYST} around a deterministic Cauchy-Peano flow of the initial system of ODEs (corresponding to  \eqref{SYST} taking $\sigma=0 $). In particular, since we consider \textit{rough} coefficients, we establish therein measurability properties and bi-Lipschitz like regularity for such flows.  The well posedness of the martingale problem for the operator associated with \eqref{SYST} is then obtained in Section \ref{THE_SEC_MP}. Section \ref{SEC_COUNTER_EX} is eventually dedicated to a class of counter examples which emphasize the \textit{almost} sharpness of our results.

\subsection*{Assumptions and main result} Our assumptions are the following:
 \begin{trivlist}
  \item[ \A{UE} ] 
  There exists $\kappa\ge 1$ s.t. for all $(t,\x)\in \R_+\times \R^{nd},\ z\in \R^d$,
$$ \kappa^{-1}|z|^2\le \langle \sigma\sigma^*(t,\x) z,z\rangle \le \kappa |z|^2,$$ 
where  $|\cdot|$ denotes the Euclidean norm, $\langle \cdot, \cdot \rangle $ 
is the inner product and  ${}^* $ stands for the transpose.
\item[\A{S}] The coefficients $\sigma(t,\cdot), \big( F_i(t,\textcolor{black}{\0})\big)_{i\in \leftB 2,n\rightB} $ are assumed to be bounded measurable in time.
Also, the diffusion coefficient $\sigma(t,\cdot) $ is  uniformly  $\eta $-H\"older continuous in space, for some $\eta>0 $ uniformly in time. The drift entries $\big(F_i(t,\cdot) \big)_{i\in\leftB 2,n\rightB} $ are s.t. for all $(z_i, \cdots, z_n)\in \R^{d(n-(i-1))} $, the mapping $z\in \R^d\mapsto F_i(t,z,z_i, \cdots, z_n) $ is in $C^{1+\eta}(\R^d,\R^d) $ uniformly in time\footnote{For the sake of clarity we chose the same regularity index for $\sigma$ and $(D_{x_{i-1}}F_i)_{i\in \leftB 2,n\rightB}$,  but the result remains true for any  $\eta_\sigma$ H\"older continuous $\sigma$ and  $\eta_{\gF,i}$ H\"older continuous  $D_{x_{i-1}}F_i$, provided $\eta_\sigma$ and $\eta_{\gF,i}$ belong to $(0,1]$.}  \textcolor{black}{and w.r.t. $(z_i,\ldots,z_n)$.}  \textcolor{black}{Moreover,} \textcolor{black}{$D_{x_{i-1}}F_i(t,\cdot) $ is bounded} . Eventually, the mappings $(z_i,\cdots,z_n)\in \R^{d(n-(i-1))} \mapsto F_i(t,z,z_i,\cdots,z_n) $ are, \textcolor{black}{for $j\in \leftB i,n\rightB $}, $\beta_i^j $-H\"older continuous in the variable $z_j$, with $\beta_i^j>0 $, uniformly in time and in $z\in \R^d$.
\item[\A{D}] The first entry of the drift $F_1$ is supposed 
to satisfy 
one of the following assumptions:
\begin{enumerate}
\item[(a)] \textcolor{black}{The measurable mapping  $t\in \R_+\mapsto F_1(t,\0) $ is bounded and}  $F_1(t,\cdot) $ is H\"older continuous in space\footnote{Actually one can assume that $F_1$ is $\beta_1^j$ H\"older continuous in the $j^{\rm th}$ variable for any $\beta_1^j$ in $(0,1]$.} uniformly in time.  
\item[(b)] \textcolor{black}{The measurable mapping  $(t,\x)\in \R_+\times \R^{nd}\mapsto F_1(t,\x)$ }is bounded.
\item[(c)] $F_1\in L^q( \R_+
,L^p( \R^{nd})), \frac{n^2d}{p}+\frac 2q<1, p\ge 2,q>2 $.
\end{enumerate}
\textcolor{black}{Observe that case (b) can be viewed as a particular case of (c), corresponding to $p=q=\infty $.  Since the techniques used to address those two cases are rather different (see Section \ref{EX_D_B_C}), we prefer to consider them separately.} 
\item[\A{H}] There exists a closed convex subset ${\mathcal E}_{i-1} \subset GL_{d}(\R)$
(set of invertible $d \times d$  matrices)
 s.t., for all
$t \geq 0$ and $(x_{i-1},\dots,x_n) \in \R^{(n-i+2)d}$, $D_{x_{i-1}} F_i(t,x_{i-1},\dots,x_n)
\in {\mathcal E}_{i-1}$.
For example, ${\mathcal E}_i$, $ i \in \leftB 1, n-1\rightB$, may be a closed ball
included in $GL_{ d}(\R)$, which is an open set.

\end{trivlist}

We say that assumption \A{A} is in force whenever \A{UE}, \A{S}, \A{H} and at least one of the three items in \A{D} hold.

\begin{theorem}[Weak Uniqueness and H\"older continuity indexes]\label{THM1}
Assume \textbf{(A)} and that the following conditions hold:
\begin{equation}
\label{COND_HOLDER}
\forall i \in \leftB 2, n\rightB,\ j\in \leftB i,n\rightB,     \beta_{i}^j\in \Bigg( \frac{2i-3}{2j-1},1\Bigg].
\end{equation}
 Then, the martingale problem associated with $(L_t)_{t\ge 0} $ where for all $\phi\in C_0^{2}(\R^{nd},\R) $, $ \x\in \R^{nd}$,
 \begin{equation}\label{DEF_L}
  L_t \varphi(\x)= \langle F(t,\x), D_\x \phi(\x) \rangle  +\frac 12 {\rm tr}\big(a(t,\x) D_{x_1}^2\phi(\x)\big)\textcolor{black}{,\quad  a := \sigma\sigma^*,}
  \end{equation}
 is well posed, i.e. there exists a unique probability measure $\P $ on $C(\textcolor{black}{\R_+},\R^{nd}) $ s.t. denoting
 by $(\X_t)_{t\ge 0}$ the associated canonical process, for every \textcolor{black}{$\varphi\in C_0^{1,2}(\textcolor{black}{\R_+}\times \R^{nd},\R) $} and conditionally to $\X_t=\x $ for $t\ge 0$, $\Big(\varphi(s,\X_s) -\varphi(t,\x)-\int_t^s (\partial_u+L_u)\varphi(u,\X_u) du\Big)_{s\ge t}$ is a $\P $-martingale. In particular weak uniqueness holds for the SDE \eqref{SYST}.

The transition probability $P(t,s,\x,\cdot) $, determined by $(L_s)_{s\ge 0} $, is s.t. for a given $T>0$, almost all $\textcolor{black}{s \in (t,T]}$ and all $\Gamma \in {\mathcal B}(\R^{nd}) $:
$P(t,s,\x,\Gamma)=\int_{\Gamma} p(\textcolor{black}{t,s},\x,\y)d\y$. 

Furthermore, we have the following Krylov-like estimate: for all fixed $T>0$ and every $f\in L^{q'}([0,T],L^{p'}(\R^{nd}))$ with $\frac{n^2d}{p'}+\frac {\textcolor{black}{2}}{q'}<2, p'>1,q'>1 $, \textcolor{black}{$(t,\x)\in [0,T]\times \R^{nd} $}: 
\begin{equation}
\label{KRYLOV_1}
|\E^{\P_{t,\x}} [\int_{t}^T f(s,\X_s)ds]|\le  C\|f\|_{L^{q'}(([0,T],L^{p'}(\R^{nd}))},
\end{equation}
where $\E^{\P_{t,\x}}$
denotes the expectation w.r.t. 
$$\P_{t,\x}[\cdot]:=\P[\cdot|\X_t=\x]\ \  {\rm \it and}\ \ C:=C(\A{A},p',q',T).$$
\end{theorem}

Hence, our Theorem allows to recover almost all the previously mentioned works on weak well posedness and provides an extension for the full chain. It also permits us to avoid any regularity assumption on the drift of the diffusion component so that we recover the Stroock and Varadhan result in the case $n=1$ up to an arbitrary small H\"older exponent on the continuity of the diffusion matrix. Concerning this last point, we feel that using the localization strategy proposed by Menozzi in \cite{meno:17} we may be able to get rid of this assumption and only assume the diffusion coefficient to be continuous in space. \textcolor{black}{Indeed, using our results (say Lemma \ref{lemme:bilipflow} below together with condition \eqref{COND_HOLDER}) should allow one to adapt the approach of \cite{meno:17} and extend Theorem \ref{THM1} to continuous diffusion matrix.}

We also underline that our result allows to deal with a large class of different drifts for the non degenerate component: the system can be globally with sub-linear growth (Assumption \A{A}-(a)), rough and bounded  (Assumption \A{A}-(b)) or only suitably integrable and rough (Assumption \A{A}-(c)).

Moreover, the following result shows that Theorem \ref{THM1} is \textit{almost} sharp. By \textit{almost}, we mean that the critical lower threshold\textcolor{black}{s} in \eqref{COND_HOLDER} and in \A{D}-(c) are not yet handled. Namely we have:
\begin{theorem}[Almost sharpness]\label{THMAS}
There exists $\gF$ satisfying \A{UE}, \A{S}, \A{H} and such that:
\begin{equation}
\label{COND_HOLDERAS}
\exists i \in \leftB 2, n\rightB,\ j\in \leftB i,n\rightB,     \beta_{i}^j < \frac{2i-3}{2j-1},
\end{equation}
or 
\begin{equation}\label{KrylovCondition}
F_1\in L^q([0,T],L^p( \R^{nd})),\quad  \frac{n^2d}{p}+\frac 2q>1,\quad p\ge 2,q>2 
\end{equation}
for which weak uniqueness fails for the SDE \eqref{SYST}.
\end{theorem}

We first emphasize that there are already some results in that direction: in \cite{beck_stochastic_2014} the Authors show that when $n=1$ and when the integrability condition \eqref{KrylovCondition} is not satisfied (\emph{i.e.} in the supercritical case) equation \eqref{SYST} does not have a weak solution. \textcolor{black}{Another counter example to that case can be found in \cite{gradinaru_existence_2013}.} Note that in comparison with the results in \cite{flandoli_multidimensional_2014}, \cite{delarue_rough_2015}, \cite{cannizzaro_multidimensional_2015}, the almost sharpness of the integrability condition \eqref{KrylovCondition} has to be understood for drifts assumed to be functions and not distributions. 

Secondly, it has been proved in \cite{chau:16}, that for all $i$ in $\leftB 2,n\rightB$ the H\"older exponents $\beta_i^i$, are also almost sharp, thanks to a class of counter examples based on stochastic perturbations of the Peano example \eqref{eq:peano}. 

Thirdly, we feel that the H\"older continuity assumption assumed on  $D_{x_{i-1}}F_i(t,\cdot)$ \textcolor{black}{is a technical artefact. Nevertheless, relaxing this assumption to consider \textcolor{black}{the} $\big(D_{x_{i-1}}F_i(t,\cdot)\big)_{i \in \leftB 2,n\rightB}$ are just continuous functions of $x_{i-1}$ is definitely more tricky. Indeed, in that case, our approach based on parametrix fails and the natural approach, \textcolor{black}{relying} on harmonic analysis techniques, seems very involved.}

\textcolor{black}{And last, but not least, let us notice that the two thresholds for the drift component (say \eqref{COND_HOLDER} and condition \A{A}-(c)) will appear many times throughout this work as a minimal value for making our proof work (see the proofs of Lemmas \ref{lemme:bilipflow}, \ref{convergence_dirac}, \ref{GROS_LEMME} and \ref{LEMME_U_F1}). This underlines the sharpness of the exponent for the strategy we used and explains why the critical case of these conditions is not investigate\textcolor{black}{d} here. It seems indeed clear for us that the critical case requires different tools as those presented here.}

\begin{remark}
 \textcolor{black}{Before entering into the proof,}  \textcolor{black}{we indicate that, from a more analytic viewpoint, the techniques we develop in the current work could also be used in order to derive well-posedness results in the mild sense (see \cite{stro:vara:79}, \cite{kolo:11}) for the corresponding degenerate parabolic PDE, which involves \textit{rough} coefficients in a weak H\"ormander setting,  when the source term belongs to appropriate Lebesgue spaces.}
\end{remark}

\section{Strategy and key tools}\label{SEC2}

Our strategy relies on \textcolor{black}{the} martingale problem approach. Hence, we face two \textcolor{black}{difficulties}: firstly, we have to show the existence of a solution to the martingale problem in our current setting, which becomes quite tricky under \A{D}-(c) while it is quite obvious under \A{D} - (a) and (b); secondly we have to show that the solution is unique which is the real core of this paper.\\

\textbf{About uniqueness.}  Usual approaches to uniqueness for the martingale problem associated  with a given operator are based on a perturbative method. 
Let us detail two of the main strategies developed in the literature. The historical one due to Stroock 	and Varadhan \cite{stro:vara:79} consists in exploiting some $L^p$ controls on the derivatives of a suitable Gaussian heat kernel (\emph{parametrix}).
It allows, in the non-degenerate diffusive case, to establish well posedness provided the diffusion coefficient is solely continuous. As a by-product of this approach, Krylov like estimates of type \eqref{KRYLOV_1} are obtained, emphasizing that the canonical process associated with the solution actually possesses a density which enjoys integrability properties up to a certain threshold. Extensions of \textcolor{black}{these} types of results to the chain \eqref{SYST} are available in \cite{meno:17}.

On the other hand, a more recent approach is due to Bass and Perkins \cite{bass:perk:09}. In the non-degenerate setting, under the stronger assumption of H\"older continuity of the diffusion coefficient, it only requires pointwise controls of an underlying \emph{parametrix} kernel. This approach has then been successfully extended \textcolor{black}{under the considered weak type H\"ormander setting} to a chain of type \eqref{SYST} in \cite{meno:10} in the diffusive case and in \cite{huan:meno:15} for more general stable driven degenerate SDEs of type \eqref{SYST} with dimension \textcolor{black}{restrictions}.
It is actually more direct than the approach of Stroock and Varadhan. However, its drawback is that it does not provide \textit{a priori} information on the density of the underlying canonical process. 

Let us underline that in both cases, the \emph{parametrix} plays a central role.
This approach consists in expanding the generator of a given stochastic process around a suitable \textit{proxy} \textcolor{black}{generator} which can be well handled. The point is then \textcolor{black}{to} control in a suitable way the associated approximation error. In our current degenerate diffusive setting, since the SDE is Brownian driven, the 
difficulty is to exhibit an appropriate Gaussian process that fulfills the previously indicated constraints.

When the drift $\gF$ is \textit{smooth} in addition to \A{A}, say globally Lipschitz continuous, it has been shown in \cite{dela:meno:10}, \cite{meno:10}, \cite{meno:17} that a good proxy is provided by the linearization around the deterministic flow associated with \eqref{SYST} (\emph{i.e.} when $\sigma=0$ therein) leading to consider a multi-scale Gaussian process as parametrix. It is therefore a natural candidate for the current work. Anyhow, under \A{A}, we do not have anymore a deterministic flow in the usual    
Cauchy-Lipschitz sense. A first difficulty is therefore to deal with non-smooth and non-unique Cauchy-Peano flow\textcolor{black}{s}. It actually turns out that any measurable flow solving \eqref{SYST} with $\sigma=0 $ is a good candidate to make our machinery work. The specific controls associated with those objects are presented in Section \ref{CAUCHY_PEANO_FLOWS}.

Also, in order to handle very rough drifts for the non degenerate component, from \A{D}-(c) $F_1\in L^{q}([0,T],L^{p}(\R^{nd}))$, 
we are led to apply the Girsanov transform to the equation with $F_1=0 $. To do so requires to have some \textit{a priori} knowledge of the corresponding density. This is why, to achieve our goal, the Stroock and Varadhan approach leading to estimate \eqref{KRYLOV_1}, seems to be the natural framework.

In comparison with the approach based on \textcolor{black}{the} Zvonkin Transform initiated in \cite{chau:16}, our  approach allows to obtain a clever analysis of the chain in the sense that we are here able to enlight the almost sharp regularity needed for each component of the drift and w.r.t. each variable. This last point is not possible \emph{via} the Zvonkin Transfom which is more global and does not permit this distinction. Accordingly to the works \cite{chau:17}, \cite{fedr:flan:prio:vove:17}, the \textcolor{black}{Zvonkin} approach seems more suited to derive strong uniqueness. In that case a global threshold appears for each variable at each  \textcolor{black}{level} of the chain.\\

\textbf{About existence.} Concerning the existence part, our proof consists in adapting to our degenerate setting the idea introduced by Portenko \cite{port:90} and used by Krylov and R\"ockner \cite{kryl:rock:05} as well to build local weak solutions in the non-degenerate case.\\

\textbf{Usual notations.} In what follows, we denote a quantity in $\R^{nd}$ by a bold letter: i.e. $\0$, stands for zero in $\R^{nd}$ and \textcolor{black}{we denote by $(X_t^1,\dots,X_t^n)_{t \geq 0}$ the components of $({\mathbf X}_t)_{t \geq 0}$}.
Introducing the embedding matrix $B$ from $\R^d$ into $\R^{nd} $, i.e. $B= (I_{d \times d} , 0, \dots, 0)^*$, where ``$*$'' stands for the transpose, we rewrite \textcolor{black}{accordingly} \eqref{SYST} in the shortened form 
\begin{equation*}
d{\mathbf X}_t = {\mathbf F}(t,{\mathbf X}_t) dt+ B \sigma(t,{\mathbf X}_{t}) dW_t,
\end{equation*}
where ${\mathbf F}=(F_1,\dots,F_n)$ is an $\R^{nd}$-valued function.\\

\textbf{The deterministic backward flow.}
In the following, we will first assume for the sake of simplicity that assumption \A{D}-(a) is in force. The extension to cases (b) and (c) will be discussed later on.  Introduce now, for fixed $ T>0,\ \y \in \R^{nd}$ and $t\in [0,T]$ the backward flow:

\begin{equation}
\label{DET_SYST}
\overset{.}{\btheta}_{t,T}(\y)=\gF(t,\btheta_{t,T}(\y)),\ \btheta_{T,T}(\y)=\y .
\end{equation} 
\begin{remark}
We mention carefully that from  the Cauchy-Peano theorem, a solution to  \eqref{DET_SYST} exists. Indeed, the coefficients are continuous and have at most linear growth. 
\end{remark}

\subsection{Linearized Multi-scale Gaussian Process and Associated Controls}\label{sec:freezing}
We now want to introduce the forward linearized flow around a solution of \eqref{DET_SYST}. Namely, we consider for $s\ge 0$ the deterministic ODE  
\begin{equation}
\label{LIN_AROUND_BK_FLOW} 
\overset{.}{\tilde{\bphi}}_s = {\mathbf F}(s,\btheta_{s,T}(\y))
+ D {\mathbf F}(s,\btheta_{s,T}(\y))[\tilde{\bphi}_s
- \btheta_{s,T}(\y)], 
\end{equation}
where for all $\z\in \R^{nd}$,\\ 
$$D\gF(s,\z)=\left (\begin{array}{ccccc}0 & \cdots & \cdots &\cdots  & 0\\
D_{z_1}F_2(s,\z) & 0 &\cdots &\cdots &0\\
0 & D_{z_2} F_3(s,\z)& 0& 0 &\vdots\\
\vdots &  0                      & \ddots & \vdots\\
0 &\cdots &     0      & D_{z_{n-1}}F_n(s,\z) & 0
\end{array}\right) 
 $$ 
denotes the subdiagonal of the Jacobian matrix ${D_{\mathbf \z} \gF(s,\cdot)} $ at point $\z$. Introduce now for a given $(T,\y)\in \textcolor{black}{\R_{+*}}\times \R^{nd}$, the resolvent $(\tilde \gR^{T,\y}(t,s))_{s,t\ge 0} $ associated with the partial gradients $(D {\mathbf F}(t,\btheta_{t,T}(\y)))_{t \ge0} $ which satisfies for $(s,t)\in (\R_{+})^2 $:
\begin{equation}
\begin{split}
\partial_{s}\tilde \gR^{T,\y}(s,t)&={ D \gF}(s,\btheta_{s,T}(\y))\tilde \gR^{T,\y}(s,t) , 
\ \tilde \gR^{T,\y}(t,t)=I_{nd\times nd}.
\end{split}
\label{DYN_RES}
\end{equation}
Note in particular that since the partial gradients are subdiagonal $ {\rm det}(\tilde \gR^{T,\y}(t,s))=1$.\\

We consider now the stochastic linearized dynamics $(\tilde \X_s^{T,\y})_{s\in [t,T]} $:
\begin{eqnarray}
&&\hspace*{-.5cm}d\tilde \X_s^{T,\y}=[\gF(s,\btheta_{s,T}(\y))+ D\gF(s,\btheta_{s,T}(\y))(\tilde \X_s^{T,\y}-\btheta_{s,T}(\y))]ds +B\sigma(s,\btheta_{s,T}(\y)) dW_s,\nonumber\\
&& \hspace*{.75cm}\forall s\in  [t,T],\
 \tilde \X_t^{T,\y}=\x. \label{FROZ}
 \end{eqnarray} 
 
 From equations \eqref{LIN_AROUND_BK_FLOW} and \eqref{DYN_RES} we explicitly integrate \eqref{FROZ} to obtain for all $v\in [t,T] $:
 \begin{equation}
 \label{INTEGRATED}
\begin{split}
\tilde  \X_v^{T,\y}&=\textcolor{black}{\tilde\gR}^{T,\y}(v,t)\x+\int_t^v  \textcolor{black}{\tilde \gR}^{T,\y}(v,s)\Big( \gF(s,\btheta_{s,T}(\y))-D\gF(s,\btheta_{s,T}(\y))\btheta_{s,T}(\y)\Big )ds\\
 & +\int_t^v \textcolor{black}{\tilde \gR}^{T,\y}(v,s)B\sigma(s,\btheta_{s,T}(\y)) dW_s.
 \end{split}
 \end{equation}
Denoting by $\tilde \btheta_{s,t}^{T,\y}(\x) $ the solution of \eqref{LIN_AROUND_BK_FLOW} with starting point $\tilde \btheta_{t,t}^{T,\y}(\x)=\x $ we rewrite:
 \begin{equation}
 \label{INTEGRATED_FLOW}
\tilde  \X_v^{T,\y}=\tilde \btheta_{v,t}^{T,\y}(\x) +\int_t^v \textcolor{black}{\tilde \gR}^{T,\y}(v,s)B\sigma(s,\btheta_{s,T}(\y)) dW_s, \ v\in [t,T].
\end{equation}

An important correspondence is now given by the following Proposition.
\begin{proposition}[Density of the linearized dynamics]\label{THE_PROP}
Under \A{A}, we have that, for all $v\in (t,T]$ the random variable $\tilde  \X_v^{T,\y} $ in \eqref{INTEGRATED_FLOW} admits a Gaussian density $\tilde p^{T,\y}(t,v,\x,\cdot) $ which writes:
\begin{eqnarray}
\forall \z \in \R^{nd}, \ \tilde p^{T,\y}(t,v,\x,\z) \label{CORRESP}
\\
:=\frac{1}{(2\pi)^{\frac{nd}2}\det(\tilde \K_{v,t}^{T,\y})^{\frac 12}}\exp\left( -\frac 12 \left\langle (\tilde \K_{v,t}^{T,\y})^{-1} (\tilde \btheta_{v,t}^{T,\y}(\x)-\z),\tilde \btheta_{v,t}^{T,\y}(\x)-\z\right\rangle\right),\notag
\end{eqnarray}
where
\begin{align*}
\tilde \K_{v,t}^{T,\y}&:=\int_t^v \textcolor{black}{\tilde  \gR}^{T,\y}(v,s) Ba(s,\btheta_{s,T}(\y))B^*\textcolor{black}{\tilde \gR}^{T,\y}(v,s)^* ds,\\
 a(s,\btheta_{s,T}(\y))&\textcolor{black}{:=}\sigma\sigma^*(s,\btheta_{s,T}(\y)) .
\end{align*}
Also, there exists $C:=C(\A{A},T)>0$ s.t. for all $k\in \leftB 0,2\rightB $, $i\in \leftB 1,n\rightB $, 
\begin{eqnarray}
&&|D_{x_i}^k\tilde p^{T,\y}(t,T,\x,\y)|\notag\\
&\le& \frac{C}{(T-t)^{k\big((i-1)+\frac{1}{2}\big) +\frac{n^2d}2}}\exp\left(-C^{-1}(T-t) \big|\T_{T-t}^{-1}\big(\x-\btheta_{t,T}(\y)\big)\big|^2\right)\notag\\
&=:&\frac{C}{(T-t)^{k\big((i-1)+\frac{1}{2}\big)}}\bar p_{C^{-1}}(t,T,\x,\y),\label{CTR_GRAD}
\end{eqnarray}
where for all $u>0$, we denote by $\T_u $ the important scale matrix:
\begin{equation}
\label{DEF_T_ALPHA}
\T_u=\left( \begin{array}{cccc}
uI_{d\times d}& 0_{d\times d}& \cdots&0_{d\times d}\\
0_{d\times d}   &u^2I_{d\times d}&0_{}& \vdots\\
\vdots & \ddots&\ddots & \vdots\\
0_{d\times d}& \cdots & 0_{}& u^{n}I_{d \times d}
\end{array}\right).
\end{equation}
\end{proposition}
\begin{proof}
Expression \eqref{CORRESP} readily follows from \eqref{INTEGRATED_FLOW}. We recall as well that, under \A{A}, the covariance matrix $\tilde \K_{v,t}^{T,\y} $ enjoys, uniformly in $\y\in \R^{nd} $ a \textit{good scaling property} in the sense of Definition 3.2 in \cite{dela:meno:10} (see also Proposition 3.4 of that reference). That is: for all fixed $T>0$, there exists $C_{\ref{GSP}}:=C_{\ref{GSP}}(\A{A},T)\ge 1$ s.t. for all $0\le v<t\le T $, for all $\y \in \R^{nd} $:
\begin{equation}
\label{GSP}
\forall \bxi \in \R^{nd},\   C_{\ref{GSP}}^{-1} (v-t)^{-1}|\T_{v-t} \bxi|^2\le \langle  \tilde \K_{v,t}^{T,\y}\bxi,\bxi\rangle \le C_{\ref{GSP}} (v-t)^{-1}|\T_{v-t} \bxi|^2.
\end{equation} 
\textcolor{black}{\begin{remark}[On the H\"ormander like non-degeneracy assumption and the good scaling property]\label{RM_GSP}
We carefully point out that the boundedness and non-degeneracy conditions expressed on the derivatives $(D_{x_{i-1}}F_i(t,\cdot))_{i\in \leftB 2,n\rightB} $\textcolor{black}{, as well as on the diffusion coefficient $\sigma(t,\cdot)$,} in assumptions \A{S}, \A{H} and \textcolor{black}{\A{UE}} are precisely explicitly used to derive this bound. The natural analogue bounds of \eqref{GSP} also hold for the inverse matrices, see again \cite{dela:meno:10}, Definition 3.2 and Lemma 3.6 therein as well as the proof of \eqref{CTR_LINEA_RETRO} below.
\end{remark}
}

Rewrite now from \eqref{INTEGRATED} and \eqref{INTEGRATED_FLOW}:
\begin{eqnarray}
&&\left\langle (\tilde \K_{T,t}^{T,\y})^{-1} (\tilde \btheta_{T,t}^{T,\y}(\x)-\y),\tilde \btheta_{T,t}^{T,\y}(\x)-\y\right\rangle\notag\\
&&=\left\langle (\tilde \gR^{T,\y}(T,t)^* (\tilde \K_{T,t}^{T,\y})^{-1} \tilde \gR^{T,\y}(T,t)(\x-\tilde \btheta_{t,T}^{T,\y}(\y)),\x-\tilde \btheta_{t,T}^{T,\y}(\y) \right\rangle,\label{PREAL_IDENT_FLOW}
\end{eqnarray}
where, \textcolor{black}{accordingly with the previous notations for the forward linearized flow $\tilde \btheta_{T,t}^{T,\y}(\x) $}, we denote:
\begin{equation*}
\tilde \btheta_{t,T}^{T,\y}(\y)=\textcolor{black}{\tilde \gR}^{T,\y}(t,T)\y-\int_t^T  \textcolor{black}{\tilde \gR}^{T,\y}(t,s)\Big( \gF(s,\btheta_{s,T}(\y))-D\gF(s,\btheta_{s,T}(\y))\btheta_{s,T}(\y)\Big )ds.
\end{equation*}
\textcolor{black}{for the corresponding linearized-backward flow starting from $\y$ at time $T$}.

Observe now \textcolor{black}{that $\tilde \btheta_{t,T}^{T,\y}(\y)=\btheta_{t,T}(\y) $}. 
Indeed, from \eqref{DYN_RES}:
\begin{eqnarray*}
 \partial_t \tilde \btheta_{t,T}^{T,\y}(\y)=D\gF(t,\btheta_{t,T}(\y)) \tilde \btheta_{t,T}^{T,\y}(\y)+\Big( \gF(t,\btheta_{t,T}(\y))-D\gF(t,\btheta_{t,T}(\y)) \btheta_{t,T}(\y)\Big),
\end{eqnarray*}
so that:
\begin{eqnarray*}
\partial_t \tilde \btheta_{t,T}^{T,\y}(\y)-\partial_t \btheta_{t,T}(\y)=D\gF(t,\btheta_{t,T}(\y)) \Big(\tilde \btheta_{t,T}^{T,\y}(\y)-\btheta_{t,T}(\y)\Big).
\end{eqnarray*}
Since $\tilde \btheta_{T,T}^{T,\y}(\y)=\btheta_{T,T}(\y)=\y $, we deduce from Gronwall's Lemma that $\tilde \btheta_{t,T}^{T,\y}(\y)=\btheta_{t,T}(\y) $ for all $t\in [0,T]$.

We carefully point out that, even though the solution to the ODE \eqref{DET_SYST} is not unique, once we have chosen a solution and consider the associated flow to construct  our linearized Gaussian model, we precisely get the identification $\tilde \btheta_{t,T}^{T,\y}(\y)=\btheta_{t,T}(\y) $ for all $t\in [0,T]$ with the \textbf{chosen} flow.

We thus get from the previous identification, equations \eqref{PREAL_IDENT_FLOW}, \eqref{GSP} and Remark \ref{RM_GSP} that there exists $C:=C(\A{A},T)>0$, s.t. for all $t\in [0,T) $,
\begin{eqnarray}
\label{CTR_LINEA_RETRO}
&&C^{-1}(T-t)|\T_{T-t}^{-1}(\x-\btheta_{t,T}(\y))|^2  \le \left\langle (\tilde \K_{T,t}^{T,\y})^{-1} (\tilde \btheta_{T,t}^{T,\y}(\x)-\y),\tilde \btheta_{T,t}^{T,\y}(\x)-\y\right\rangle\\
&& \hspace*{5cm} \le C(T-t)|\T_{T-t}^{-1}(\x-\btheta_{t,T}(\y))|^2.\notag
\end{eqnarray}
\textcolor{black}{Indeed, from \eqref{GSP} it is easily derived that the spectrum of $\widehat{ \tilde \K}_1^{T,t,T,\y}:= (T-t)\T_{T-t}^{-1} \tilde \K_{T,t}^{T,\y}\T_{T-t}^{-1}$ lies in \textcolor{black}{$[C_{\ref{GSP}}^{-1},C_{\ref{GSP}}] $}. So does the spectrum of $(\widehat{ \tilde \K}_1^{T,t,T,\y})^{-1}:=(T-t)^{-1}\T_{T-t} (\tilde \K_{T,t}^{T,\y})^{-1}\T_{T-t} $. From Lemma 3.6 in \cite{dela:meno:10} (see also Lemma 6.2 in \cite{meno:17} for notations closer to the current framework) we can also write $\tilde \gR^{T,\y}(T,t)=\T_{T-t}{\widehat {\tilde \gR}}^{T,t,T,\y}(1,0)\T_{T-t}^{-1} $ where ${\widehat{\tilde  \gR}}^{T,t,T,\y}(1,0) $ is a non-degenerate bounded matrix whose bounds do not depend on $t,T$, i.e. ${\widehat {\tilde \gR}}^{T,t,T,\y}(1,0)$ is a non-degenerate \textit{macro} matrix. Therefore, introducing $\tilde \gH_{T,t}^{T,\y}=\tilde \gR^{T,\y}(t,T) \tilde \K_{T,t}^{T,\y} \tilde \gR^{T,\y}(t,T)^* $, we can rewrite from \eqref{PREAL_IDENT_FLOW}  that: 
\begin{eqnarray}
&&
\left\langle (\tilde \K_{T,t}^{T,\y})^{-1} (\tilde \btheta_{T,t}^{T,\y}(\x)-\y),\tilde \btheta_{T,t}^{T,\y}(\x)-\y\right\rangle\notag\\
&&\hspace*{1.55cm}=\left \langle (\tilde \gH_{T,t}^{T,\y})^{-1} (\x- \btheta_{T,t}(\y)),\x- \btheta_{T,t}(\y)\right\rangle,\notag\\
&&(\tilde \gH_{T,t}^{T,\y})^{-1}= \tilde \gR^{T,\y}(T,t)^*(\tilde \K_{T,t}^{T,\y})^{-1} \tilde \gR^{T,\y}(T,t)\notag\\
&&\hspace*{1.55cm}=(T-t)\T_{T-t}^{-1}\big({\widehat {\tilde \gR}}^{T,t,T,\y}(1,0)\big)^*
(\widehat{ \tilde \K}_1^{T,t,T,\y})^{-1} {\widehat {\tilde \gR}}^{T,t,T,\y}(1,0) \T_{T-t}^{-1}.\label{GSP_H}
\end{eqnarray}
The previous non-degeneracy properties of $(\widehat{ \tilde \K}_1^{T,t,T,\y})^{-1}, {\widehat {\tilde \gR}}^{T,t,T,\y}(1,0)$ then readily give \eqref{CTR_LINEA_RETRO}. Put it differently, the matrix $\tilde \gH_{T,t}^{T,\y} $ also satisfies a \textit{good scaling property}.}

We now deduce from \textcolor{black}{\eqref{GSP_H}} and \eqref{GSP} that \eqref{CTR_GRAD} holds for $k=0 $.
\textcolor{black}{We now write}  from \eqref{CORRESP}:
\begin{align*}
&D_{x_i}\tilde p^{T,\y}(t,T,\x,\y)\\
=&\frac{1}{(2\pi)^{\frac{nd}2}\det(\tilde \K_{T,t}^{T,\y})^{\textcolor{black}{\frac 12}}}D_{x_i}\exp\left( -\frac 12 \left\langle (\tilde \gH_{T,t}^{T,\y})^{-1} (\x-\btheta_{t,T}(\y)),\x-\btheta_{t,T}(\y)\right\rangle\right)\\
=&\frac{\Big(-(\tilde \gH_{T,t}^{T,\y})^{-1} (\x-\btheta_{t,T}(\y))\Big)_i}{(2\pi)^{\frac{nd}2}\det(\tilde \K_{T,t}^{T,\y})^{\textcolor{black}{\frac 12}}}\exp\left( -\frac 12 \left\langle (\tilde \gH_{T,t}^{T,\y})^{-1} (\x-\btheta_{t,T}(\y)),\x-\btheta_{t,T}(\y)\right\rangle\right).
\end{align*}
\textcolor{black}{We thus derive from \eqref{CTR_LINEA_RETRO} and \eqref{GSP_H}}:
\begin{eqnarray*}
|D_{x_i}\tilde p^{T,\y}(t,T,\x,\y)|\le \frac{C}{(T-t)^{(i-1)+\frac 12+\frac{n^2d}2}} \exp\left( -C^{-1}(T-t) \left|\T_{T-t}^{-1}\big(\x-\btheta_{t,T}(\y)\big)\right|^2\right),
\end{eqnarray*}
which proves \eqref{CTR_GRAD} for $k=1$. The case $k=2 $ is derived similarly.
\end{proof}

\subsection{Regularity and measurability of the Cauchy-Peano flow} \label{CAUCHY_PEANO_FLOWS}
\textcolor{black}{Let us recall, as indicated before equation \eqref{DET_SYST}, that we work first under assumption \A{D}-(a). In this setting,}
we mention that the delicate part here \textcolor{black}{consists in dealing} with the nonlinear flow $\btheta_{t,s}(\y) $. Because of our low H\"older regularity, we face two problems: \textcolor{black}{one has to}
choose a measurable flow of \eqref{DET_SYST} (which is very important to make licit any integration of this flow along the terminal condition) \textcolor{black}{and}
this flow \textcolor{black}{must} have the appropriate regularity to deal with our parametrix kernel, say \textcolor{black}{e.g.} bi-Lipschitz as in \cite{dela:meno:10}, \cite{meno:10}, \cite{meno:17}. 

\textcolor{black}{The first issue is addressed by Lemma \ref{lemme:measu} below}. 
\textcolor{black}{The second problem} is quite involved and requires also a careful  analysis. Indeed our approach, based on parametrix kernel, makes an intensive use of the gradient estimate of the frozen transition density $\tilde{p}^{T,\y}$ given in \eqref{CORRESP}. This leads us to study the space integral of the Gaussian like function $\bar p_{C^{-1}}$ defined by \eqref{CTR_GRAD} w.r.t. the backward variable $\y$. The crucial point is that such an integral then involves the backward flow with argument the integration variable. In a smooth setting, such a problem is easily handled through a change of variable. When working with non-continuously differentiable coefficients, one may also use the bi-Lipschitz property of the flow to change its argument from the integration variable to the fixed initial one (see e.g. \cite{dela:meno:10} where $\gF$ is Lipschitz in space). In the current setting the flow is not smooth enough either to perform a change of variable nor to use the bi-Lipschitz estimate. Nevertheless, using a careful regularization procedure which precisely works under the condition \eqref{COND_HOLDER} on the H\"older continuity exponents, we show in Lemma \ref{lemme:bilipflow} below that the chosen flow satisfies an \emph{approximate} bi-Lipschitz estimate. This \emph{approximate} bi-Lipschitz estimate is sufficient to deal with our parametrix kernel.\\

\begin{lemma}\label{lemme:measu}
\textcolor{black}{For a given $T>0$, 
there exists a measurable mapping $\textcolor{black}{(s,t,\x) \in [0,T]^2\times}$ $\R^{nd}\mapsto \btheta_{t,s}(\x)$ s.t. $\btheta_{t,s}(\x)=\x+\int_s^t \gF(v,\btheta_{v,s}(\x))dv $}.
\end{lemma}
\begin{proof}
\textcolor{black}{The proof follows} from the result of \cite{zubelevich_several_2012} and usual covering arguments.
\end{proof}
\textcolor{black}{From now on, we choose by simplicity to work with a given measurable flow $\btheta_{t,s}(\x) $ which exists by the previous lemma.} 
\begin{lemma}\label{lemme:bilipflow}
There exist constants $(C_{\ref{lemme:bilipflow}},C_{\ref{lemme:bilipflow}}'):=(C_{\ref{lemme:bilipflow}}, C_{\ref{lemme:bilipflow}}')(\A{A},T)\ge 1$ s.t. for all \textcolor{black}{$0\le t<s\le T $} small enough:
\begin{align}
\label{EQ_EQUIV_FLOW}
C_{\ref{lemme:bilipflow}}^{-1} (s-t)  |\T_{s-t}^{-1}(\btheta_{s,t}(\x)-\y)|^2 - C_{\ref{lemme:bilipflow}}' &\le (s-t) |\T_{s-t}^{-1}(\x-\btheta_{t,s}(\y))|^2\notag\\
&\le C _{\ref{lemme:bilipflow}}(s-t)  |\T_{s-t}^{-1}(\btheta_{s,t}(\x)-\y)|^2 + C_{\ref{lemme:bilipflow}}'.
\end{align}
Also, for any measurable flow $\check \btheta_{s,t} $ \textcolor{black}{satisfying the integral equation in Lemma \ref{lemme:measu} and} possibly different from the chosen one $\btheta_{s,t} $, it also holds that:
\begin{align}
\label{EQ_EQUIV_FLOW_DIFF}
C_{\ref{lemme:bilipflow}}^{-1} (s-t)  |\T_{s-t}^{-1}(\check\btheta_{s,t}(\x)-\y)|^2 - C_{\ref{lemme:bilipflow}}' &\le (s-t) |\T_{s-t}^{-1}(\x-\btheta_{t,s}(\y))|^2\notag\\
&\le C _{\ref{lemme:bilipflow}}(s-t)  |\T_{s-t}^{-1}(\check \btheta_{s,t}(\x)-\y)|^2 +C _{\ref{lemme:bilipflow}}'.
\end{align}
\end{lemma}
Lemma \ref{lemme:bilipflow} is a key tool for our analysis. It roughly says that, even though the drift coefficient is not smooth, we can still expect a kind of equivalence of the rescaled forward and backward flows (which has been thoroughly used in the papers \cite{dela:meno:10}, \cite{meno:10}, \cite{meno:17} for Lipschitz drifts) up to an additional constant contribution. \textcolor{black}{This is precisely the result of equation \eqref{EQ_EQUIV_FLOW}. Also, since uniqueness here possibly fails for the flows, equation \eqref{EQ_EQUIV_FLOW_DIFF} gives that the bound still holds for two arbitrary measurable flows. This specific property will be used later on  in the proof of Lemma \ref{convergence_dirac} below in Appendix \ref{sec:ProofLemma}}.

It turns out that, \textcolor{black}{the new contribution in \eqref{EQ_EQUIV_FLOW}, \eqref{EQ_EQUIV_FLOW_DIFF}} does not perturb the analysis of the parametrix kernels associated with the density of $\tilde \X_T^{T,\y} $ starting from $\x $ at time $t\in [0,T)$ given in Proposition \ref{THE_PROP}. We refer to Section \ref{WP_SEC_HOLDER} for details.

\begin{proof}
\textcolor{black}{We focus on the proof of \eqref{EQ_EQUIV_FLOW_DIFF}. Indeed, \eqref{EQ_EQUIV_FLOW} is derived as a special case taking the same flows, i.e. $\check \btheta=\btheta $}.
Considering now two measurable flows $\btheta,\check \btheta $ provided by \textcolor{black}{the integral equation in} Lemma \ref{lemme:measu}, we write from the \textcolor{black}{integral} dynamics:
\begin{align}
\label{DIFF_DYN_FLOW}
&(s-t)^{\frac 12} \T_{s-t}^{-1}(\x-\btheta_{t,s}(\y))\notag\\
&= (s-t)^{\frac 12} \T_{s-t}^{-1}\Big[(\textcolor{black}{\check\btheta_{s,t}(\x)}-\y)-\int_{t}^s \Big(\gF(u,\textcolor{black}{\check\btheta_{u,t}(\x)})-\gF(u,\btheta_{u,s}(\y))\Big)du\Big].\notag\\
&=(s-t)^{\frac 12} \T_{s-t}^{-1}\big(\textcolor{black}{\check \btheta_{s,t}(\x)}-\y\big)- {\cI}_{s,t}(\x,\y),\\
{\cI}_{s,t}(\x,\y)&\textcolor{black}{:=}(s-t)^{\frac 12} \T_{s-t}^{-1}\int_{t}^s \Big(\gF(u,\textcolor{black}{\check \btheta_{u,t}(\x)})-\gF(u,\btheta_{u,s}(\y))\Big)du .\notag
\end{align}
We aim at establishing that
\begin{equation}
\label{CTR_I1}
|\cI_{s,t}(\x,\y)|\le C\Big\{1\textcolor{black}{+(s-t)^{-1}\int_t^s (s-t)^{\frac 12}|\T_{s-t}^{-1}(\check \btheta_{u,t}(\x)-\btheta_{u,s}(\y))|   du}\Big\},
\end{equation}
which together with \eqref{DIFF_DYN_FLOW} \textcolor{black}{and the Gronwall lemma} gives the r.h.s. of \textcolor{black}{\eqref{EQ_EQUIV_FLOW_DIFF}}. The l.h.s. could be derived similarly to the analysis we now perform.

Since the function $\gF $ is not Lispchitz, we will thoroughly use, as crucial auxiliary tool, some appropriate mollified flows. 
We first denote by $\delta \in \R^n\otimes \R^n$, a matrix whose  entry $\delta_{ij} $ is strictly positive for indexes $i\in \leftB 2,n\rightB $ and $j\ge i$. We then define for all $
\ v\in [0,T],\ \z\in \R^{nd} $, $i\in \leftB 2,n\rightB $, 
\begin{align}
\label{DEF_CONV_ADHOC}
F_i^\delta(v,\z^{i-1,n})&:=F_i(v,\cdot)\star \rho_{\delta_{i,.}}(\z)\\
&=\int_{\R^{d(n-(i-1))}} F_i(v,z_{i-1},z_i-w_i,\cdots ,z_n-w_n)\rho_{\delta_{i,.}}(w)dw.\notag
\end{align} 
Here, for all  $w=(w_i,\cdots,w_n)\in \R^{d(n-(i-1))}$, 
$$\rho_{\delta_{i,.}}(w):=\frac{1}{\prod_{j=i}^n\delta_{ij}^{d}}\rho_i\left(\frac{w_i}{\delta_{ii}},\frac{w_{i+1}}{\delta_{i(i+1)}},\cdots,\textcolor{black}{\frac{w_{n}}{\delta_{in}}} \right),$$ where $ \rho_i:\R^{d(n-(i-1))}\rightarrow \R_{+}$ is a standard mollifier, i.e. $\rho_i $ \textcolor{black}{is smooth}, has compact support and $\int_{\R^{d(n-(i-1))}} \rho_i(\z)d\z=1 $. We denote 
$$\gF^\delta(v,\z)=(F_1(v,\z),F_2^\delta(v,\z), \cdots,F_n^\delta(v,\z)).$$ 

\textcolor{black}{Note carefully that, under the considered assumptions \A{A}, the  $(F^\delta_i)_{i\in \leftB 2,n\rightB} $ are thus Lipschitz continuous functions, with explosive Lipschitz constant w.r.t. the mollifying procedure in the variables $(x_i,\cdots,x_n) $. More precisely, standard arguments from approximation theory give that, under the current assumptions, the Lipschitz constant of $F^\delta_i$ w.r.t. its $j^{\rm{th}}$ variable blows up at rate $\delta_{ij}^{-1+\beta_i^j}$}.

\paragraph{\textbf{\textcolor{black}{Controls associated with the mollification procedure.}}}
The first key point is that the regularized drift $\mathbf{F}^\delta$ only appears \textcolor{black}{in a time integral for} 
our analysis \textcolor{black}{(see equations \eqref{DIFF_DYN_FLOW} and \eqref{DECOUP_I1})}. \textcolor{black}{So, the parameters $\delta_{ij}$ only have} to satisfy that there exists $C:=C(\A{A},T)>0$ such that for all $\z$ in $\R^{nd}$, for all $u$ in $[t,s]$:
\begin{eqnarray}\label{approxrescaleddrfit}
\Big|(s-t)^{\frac 12}\T_{s-t}^{-1}\Big(\gF(u,\z)-\gF^\delta(u,\z)\Big)\Big|\leq C(s-t)^{-1}.
\end{eqnarray}
From the definition of our regularization procedure in \eqref{DEF_CONV_ADHOC} this means that $\delta_{ij}$ must be such that
\begin{eqnarray}\label{eq:controleT1}
\textcolor{black}{\sum_{i=2}^n (s-t)^{\frac 12-i}\sum_{j=i}^n \delta_{ij}^{\beta_i^j} 
\leq C(s-t)^{-1}.}
\end{eqnarray}
Hence, one can choose $\textcolor{black}{\delta_{ij}= (s-t)^{(i-\frac 32)\frac{1}{\beta_i^j}}}$ which yields $\textcolor{black}{(s-t)^{\frac 12-i}\delta_{ij}^{\beta_i^j}=(s-t)^{-1}}$. This choice of $\delta_{ij}$ will be \textcolor{black}{kept for} the rest of the proof.

The second key point relies on the fact that, \textcolor{black}{for} this choice of the regularization parameter, the rescaled drift $\mathbf F^\delta$ satisfies an \emph{approximate} Lipschitz property whose Lipschitz constant, once the drift is integrated, does not yield any additional singularity. Namely, there exists a $C:=C(\A{A},T)$ such that for all \textcolor{black}{$u$} in $[s,t]$, for all $\z, \z'$ in $\R^{nd}$
\begin{eqnarray}\label{approximateLippourF}
&&\Big|(s-t)^{\frac 12}\T_{s-t}^{-1}\Big(\gF^\delta(\textcolor{black}{u},\z)-\gF^\delta(\textcolor{black}{u},\z')\Big)\Big| \\
&&\le C\Bigg( (s-t)^{-\frac 12}+ (s-t)^{-1}|(s-t)^{\frac 12}\T_{s-t}^{-1}(\z-\z')| \Bigg).\notag
\end{eqnarray}
Indeed, \textcolor{black}{as already underlined at the end of the previous paragraph,} the $(F^\delta_i)_{i\in \leftB 2,n\rightB} $ are Lipschitz continuous functions (with potentially explosive Lipschitz constant in the variables $(x_i,\cdots,x_n) $ for $F_i^\delta $ because of the regularization procedure) and $F_1 $ is $\beta_1^j>0$ H\"older continuous in the $j^{\rm th} $ variable for arbitrary $(\beta_1^j)_{j\in \leftB 1,n\rightB} $ in $\textcolor{black}{(}0,1] $. The Young inequality then yields that there exists $C_{\beta_1^j}>0$ s.t. for all $\x\in \R^{nd},\ |\x|^{\beta_1^j}\le C_{\beta_1^j}(1+|\x|)  $. Hence,
\begin{eqnarray*}
&&\Big|(s-t)^{\frac 12}\T_{s-t}^{-1}\Big(\gF^\delta(v,\z)-\gF^\delta(v,\z')\Big)\Big| \\
&\le & C\bigg((s-t)^{-\frac 12}(1+|(\z-\z')
|) \\
&&+\sum_{i=2}^n (s-t)^{\frac 12-i}\Big(|(\z-\z')_{i-1}| +\sum_{j=i}^n\frac{|(\z-\z')_j|}{\delta_{ij}^{1-\beta_i^j}}  \Big) \bigg)\\
&\le &C\Bigg( \textcolor{black}{(s-t)^{-\frac12}}+  |(s-t)^{\frac 12}\T_{s-t}^{-1}(\z-\z')|\\
&& \qquad \times \bigg(1+(s-t)^{-1}+\sum_{i=2}^n\sum_{j=i}^n\frac{(s-t)^{j-i}}{\delta_{ij}^{1-\beta_i^j}}\bigg)\Bigg).
\end{eqnarray*}
Hence \eqref{approximateLippourF} follows from the fact that, from our previous choice of $\delta_{ij}$, one gets
\begin{equation}\label{EXPOSANTS_QUI_MATCHENT_AUX_PETITS_OIGNONS_2}
\frac{(s-t)^{j-i}}{\delta_{ij}^{1-\beta_i^j}}=(s-t)^{(j-i)-(i-\frac 32)\frac{1}{\beta_i^j}(1-\beta_i^j)}\le C (s-t)^{-1},
\end{equation}
since from \textcolor{black}{the assumption} \eqref{COND_HOLDER} \textcolor{black}{on the indexes of H\"older continuity}:
 $$ \textcolor{black}{\beta_i^j>\frac{2i-3}{2j-1}\iff} (j-i)-(i-\frac 32)(1-\beta_i^j)/\beta_i^j>-1.$$

\paragraph{\textcolor{black}{\textbf{\textcolor{black}{Derivation of the final bound.}}}}
We are now in position to bound the term $\cI_{s,t}(\x,\y)$ defined in \eqref{DIFF_DYN_FLOW}. Under \A{A} we have:
\begin{eqnarray}
|\cI_{s,t}(\x,\y)|&\le& \int_{t}^{s} du  |(s-t)^{\frac 12}\T_{s-t}^{-1} \bigl(\gF(u,\textcolor{black}{\check \btheta_{u,t}(\x)})-\gF(u,\btheta_{u,s}(\y)) \bigr)|\notag\\
&\le & \int_t^s du \Big|(s-t)^{\frac 12}\T_{s-t}^{-1}\Big(\gF(u,\textcolor{black}{\check\btheta_{u,t}(\x)})-\gF^\delta(u,\textcolor{black}{\check \btheta_{u,t}(\x)})\Big)\Big|\notag\\
&&+\int_t^s du \Big|(s-t)^{\frac 12}\T_{s-t}^{-1}\Big(\gF^\delta(u,\textcolor{black}{\check \btheta_{u,t}(\x)})-\gF^\delta(u,\textcolor{black}{\btheta_{u,s}}\textcolor{black}{(\y)})\Big)\Big|\notag\\
&&+ \int_t^s du \Big|(s-t)^{\frac 12}\T_{s-t}^{-1}\Big(\gF^\delta(u,\textcolor{black}{\btheta_{u,s}(\y)})-\gF(u,\textcolor{black}{\btheta_{u,s}(\y)})\Big)\Big|\notag\\
&=:&\cI_{s,t}^{1}(\x,\y)+\cI_{s,t}^{2}(\x,\y)+\cI_{s,t}^{3}(\x,\y). \label{DECOUP_I1}
\end{eqnarray}
\textcolor{black}{Using} \eqref{approxrescaleddrfit},  we obtain that there exists $C:=C(\A{A},T)$ s.t. for all $0\le t< s\le T,\ \x,\y\in \R^{nd} $:
\begin{equation}
\label{I11}
\textcolor{black}{|\cI_{s,t}^{1}(\x,\y)|+|\cI_{s,t}^{3}(\x,\y)|}\le C.
\end{equation}
Finally, one can use \eqref{approximateLippourF}  to derive that 
for all $0\le t< s\le T,\ \x,\y\in \R^{nd} $:
\begin{eqnarray*}
|\cI_{s,t}^{2}(\x,\y)|&\le& C\Bigg[\textcolor{black}{(s-t)^{\frac 12}+\int_t^s (s-t)^{-1}\big((s-t)^{\frac 12}|\T_{s-t}^{-1}(\check \btheta_{u,t}(\x)-\btheta_{u,s}(\y))|\big)   du}\Bigg]\\
&\le & C\Big\{\textcolor{black}{1+(s-t)^{-1}\int_t^s (s-t)^{\frac 12}|\T_{s-t}^{-1}(\check \btheta_{u,t}(\x)-\btheta_{u,s}(\y))|   du}\Big\},
\end{eqnarray*}
\textcolor{black}{up to a modification of $C:=C(\A{A},T)$  for the last inequality}.
From this last equation together with \eqref{I11}, we therefore derive \eqref{CTR_I1}. The proof is complete.

\end{proof}

\subsection{Frozen Green kernels and associated PDEs}\label{sectionGreen}
In this paragraph we introduce useful tool\textcolor{black}{s} for the analysis of the martingale problem. Namely, we consider suitable green kernels associated with the previously defined frozen process and establish the \textit{Cauchy problemma} which it solves.

For all $ f\in C^{1,2}_0([0,T)\times \R^{nd},\R)$,  and $\varepsilon\ge 0 $ meant to be small, we define the Green function:
\begin{equation}\label{GREEN_KERNEL}
\forall (t,\x) \in [0,T)\times\R^{nd},  \tilde G^\varepsilon f(t,\x) = \int_{\textcolor{black}{(t+\varepsilon)}\wedge T}^T ds \int_{\R^{nd}} d\y \tilde{p}^{s, \y}(t,s, \x,\y) f(s,\y).
\end{equation}
\textcolor{black}{We point out that the measurability of the flow in $(s,\y)$ established in Lemma \ref{lemme:measu} precisely gives that  $\tilde{p}^{s, \y}(t,s, \x,\y)$ is a measurable function of these parameters and ensures that the above Green function is properly defined}.
Denote by  $(\tilde L_t^{s,\y})_{t\in [0,s]}$ the generator of $(\tilde \X_t^{s,\y})_{t\in [0,s]} $, i.e. for all $\varphi \in C_0^\infty(\R^{nd},\R),\  \x\in \R^{nd}$,
\begin{eqnarray}
 \tilde L_t^{s,\y} \varphi(\x)&:=&
\langle \gF(t,\btheta_{t,s}(\y))+ D \gF(t,\btheta_{t,s}(\y)) (\x-\btheta_{t,s}(\y)), D_\x \varphi(\x)\rangle\nonumber \\
&&+ \frac 12 {\rm tr}(\sigma\sigma^*(t,\btheta_{t,s}(\y) ) D_{x_1}^2 \varphi(\x) ).\label{DEF_TILDE_L}
\end{eqnarray}

One now easily checks that:
\begin{equation}\label{relation differentielle}
\forall (t,\x,\z) \in [0,s) \times (\R^{nd})^2,  \Big(\partial_t + \tilde{L}^{s, \y}_t \Big) \tilde{p}^{s,\y} (t,s,\x,\z) = 0.
\end{equation}
However, we carefully mention that some care is needed to establish the following lemma, whose proof is postponed to \textcolor{black}{Appendix} \ref{sec:ProofLemma}, which is crucial to derive that $\tilde G f\textcolor{black}{ := \tilde G^{0}f}$ actually solves an appropriate Cauchy \textit{like} problem.
\begin{lemma}[Dirac convergence of the frozen density]\label{convergence_dirac}
For all bounded continuous function $f:\R^{nd}\rightarrow \R, \x\in \R^{nd}$, setting for all $(\varepsilon,t)$ in $(\R_+\backslash\{0\}) \times [0,T]$, $f_{\varepsilon,t}(\x):=\int_{\R^{nd}} f(\y) \tilde{p}^{t+\varepsilon,\y}(t,t+\varepsilon,\x,\y)d\y $, we have:
\begin{equation}
\sup_{t \in [0,T]}\left| 
 f_{\varepsilon,t}(\x)-f(\x) \right| \underset{\varepsilon \downarrow 0} {\longrightarrow} 0.
\end{equation}
\end{lemma}
We emphasize that the above lemma is not a direct consequence of the convergence of the law of the frozen process towards the Dirac mass (see e.g. \eqref{relation differentielle}). Indeed, the integration parameter is also the freezing parameter which makes things more subtle.

We also have the following \textcolor{black}{related} result \textcolor{black}{again proved in Appendix \ref{sec:ProofLemma}}.
\begin{lemma}[$L^{q'}-L^{p'} $ convergence for the mollification with the frozen density]\label{convergence_LPLQ}
For {\color{black} all $f\in L^{q'}\big([0,T],L^{p'}(\R^{nd})\big)$, $q',p' <\infty$,}
setting\textcolor{black}{, for any $\varepsilon > 0$,} $f_\varepsilon(t,\x):=\int_{\R^{nd}} f(t+\varepsilon,\y)\textcolor{black}{\mathbb{I}_{t \in [0,T-\varepsilon]}} \tilde{p}^{t+\varepsilon,\y}(t,t+\varepsilon,\x,\y)d\y $, we have for all $p',q'>1 $:
\begin{equation}
 \|f_\varepsilon-f \|_{L^{q'}([0,T],L^{p'}(\R^{nd}))} \underset{\varepsilon \downarrow 0} {\longrightarrow} 0.
\end{equation}
\end{lemma}

 Introducing for all $f\in C^{1,2}_0([0,T)\times \R^{nd},\R),\ \varepsilon\ge 0 $ \textcolor{black}{and $(t,\x)\in [0,T)\times \R^{nd} $} the quantity:
\begin{equation}
\label{DEF_M_EPSILON}
\tilde M_{t
}^\varepsilon f(t,\x) =\int_{\textcolor{black}{(t+\varepsilon)\wedge T}}^Tds \int_{\R^{nd}} d\y \tilde{L}_t^{s,\y} \tilde{p}^{s,\y}(t,s,\x,\y)f(s,\y),
\end{equation}
we derive from \eqref{relation differentielle} and \textcolor{black}{Proposition \ref{THE_PROP}} 
that the following equality holds for all $\varepsilon>0 $ small enough:
\begin{equation}\label{relation differentielle 2}
\partial_t \tilde G^\varepsilon f(t,\x)+ \tilde M_{t
}^\varepsilon f(t,\x) = -f_\varepsilon(t,\x),\ \textcolor{black}{\forall (t,\x)\in [0,T)\times \R^{nd}},
\end{equation}
\textcolor{black}{denoting 
by} $f_\varepsilon(t,\x)=\int_{\R^{nd}} f(t+\varepsilon,\y)\textcolor{black}{\mathbb{I}_{t \in [0,T-\varepsilon]}} \tilde{p}^{t+\varepsilon,\y}(t,t+\varepsilon,\x,\y)d\y $, i.e. the time argument of $f$ is also shifted \textcolor{black}{and truncated above $T-\varepsilon$}. Observe here that the localization w.r.t. $\varepsilon$ is precisely needed to exploit directly \eqref{relation differentielle} and derive \eqref{relation differentielle 2}, for a given fixed $\varepsilon>0 $, by usual domination arguments. We mention that when $\varepsilon=0 $, it follows from the definition \eqref{GREEN_KERNEL} of our Green kernel that the smoothness on $f$ is not a sufficient condition to derive the smoothness of $\tilde G f=\tilde G^0f$. This comes from the dependence of the covariance matrix in $\tilde p^{s,\y}$ w.r.t. the integration variable (see \eqref{CORRESP}).

\begin{proposition}{Pointwise control of the Green Kernel.}\label{prop:pointwise_green_kernel}
There exists $C(T):=C(\A{A},T) \underset{T\rightarrow 0}{\longrightarrow} 0$ s.t. for all $(t,\x)\in  [0,T]\times \R^{nd} $ and all $f\in L^{q'}\big([0,T],L^{p'}(\R^{nd})\big)$ s.t. $\frac{n^2d}{p'}+\frac 2{q'}<2, p'>1,q'>1 $, $\varepsilon\in [0,T-t] $:
\begin{equation}
\label{CTR_GREEN_PONCT_1}
|\tilde G^\varepsilon f(t,\x)|\le C(T) \|f\|_{L^{q'}([0,T],L^{p'}(\R^{nd}))}.
\end{equation}
\end{proposition}
\begin{proof}
From \eqref{CTR_GRAD} with $k=0$ and Lemma \ref{lemme:bilipflow} we have that for all $\varepsilon \in [0,T-t] $,
\begin{eqnarray}
&&|\tilde G^\varepsilon f(t,\x)| \label{expourmesflow} \\
&\leq &\textcolor{black}{C}\int_t^T ds \int_{\R^{nd}} d\y |f(s,\y)| {\bar p}_{C^{-1}}(t,s,\x,\y)\notag \\
&\leq&\textcolor{black}{C} \int_t^Tds\int_{\R^{nd}} d\y |f(s,\y)| \frac{\exp\left(-C^{-1}(s-t)|\T_{s-t}^{-1}(\x-\btheta_{t,s}(\y))|^2 \right)}{ (s-t)^{\frac{n^2d}{2}}  }\notag\\
&\leq &\textcolor{black}{C} \int_t^Tds\int_{\R^{nd}} d\y |f(s,\y)| \frac{\exp\left(-C^{-1}\Big[ (s-t) |\T_{s-t}^{-1}\btheta_{s,t}(\x)-\y)|^2 + 1\Big]\right)}{ (s-t)^{\frac{n^2d}{2}}}, \notag
\end{eqnarray}
up to a modification of $C$.
So,  the result follows from the H\"older inequality and the condition on the exponents $p'$ and $q'$. Denoting by $\tilde p'$, $\tilde q'$ the conjugate of $p'$ and $q'$ respectively we indeed  have
\begin{eqnarray*}
&&\int_t^T ds\left|\int_{\R^{nd}}  d\y     \left|\frac{\exp\left(-C^{-1} (s-t) ( |\T_{s-t}^{-1}(\btheta_{s,t}(\x)-\y)|^2)\right)}{ (s-t)^{\frac{n^2d}{2}}}  \right|^{\tilde p'}\right|^{\frac{\tilde q'}{\tilde p'}}\\
&\leq & C \int_t^T ds \frac{1}{(s-t)^{\frac{n^2d}{2}(\tilde p'-1)(\frac{\tilde q'}{\tilde p'})}}  <+\infty \Leftrightarrow \frac{n^2d (\tilde p'-1)}{2} \frac{\tilde q'}{\tilde p'} <1 \Leftrightarrow \frac{n^2d}{2 p'} + \frac{1}{q'} <1. 
\end{eqnarray*}
\end{proof}

\section{Well posedness of the corresponding martingale problem }\label{THE_SEC_MP}

We have now \textcolor{black}{given} the main tools \textcolor{black}{needed} to prove our main results: the well posedness of the martingale problem associated with $(L_t)_{t \geq 0}$ defined \textcolor{black}{in} \eqref{DEF_L} and the corresponding Krylov-type estimates. This section is organized as follows: we first investigate the well posedness under Assumption \A{D}-(a). In that case, the existence part is not a challenge, since it  \textcolor{black}{readily} follows  from previous\textcolor{black}{ly} known results \textcolor{black}{based on compactness arguments that exploit the sublinear structure of the drift $\gF $}, while the uniqueness part is quite more delicate. \textcolor{black}{As a by-product of our approach to uniqueness we derive the Krylov like estimate}.

The scheme used for proving uniqueness under \textcolor{black}{\A{D}-(a)} will be a major tool to extend our result for uniqueness under \A{D}-(b) (the existence part under that assumption being a trivial application of \textcolor{black}{the} Girsanov Theorem) and then for existence and uniqueness under \A{D}-(c). Indeed, under this last assumption, even the existence part requires to derive first some Krylov type estimates.  We will precisely exploit those established under \A{A}-(a) considering first $F_1=0$ and then cope with the true $L^q-L^p $ drift through a Girsanov argument. \textcolor{black}{
The approach is in some sense similar to the one of Krylov and R\"ockner \cite{kryl:rock:05} or Fedrizzi \textit{et al.} \cite{fedr:flan:prio:vove:17} for the Girsanov part. The main difference is that in the quoted work the required Krylov like estimate readily followed from the explicit density of the unperturbed process at hand. The Brownian motion in \cite{kryl:rock:05}, the joint density of the Brownian motion and its integral in \cite{fedr:flan:prio:vove:17}. We here precisely show that, first under \A{A}-(a), the solution to the martingale problem has a density which satisfies a similar Krylov type estimate.
We actually prove that \textit{any} solution to the martingale problem satisfies such an estimate (see equation \eqref{KRYLOV_POUR_TOUS} below).}\\

\textcolor{black}{It is precisely to deal with $L^{q}-L^{p}$ drifts (under Assumption \A{D}-(c)) that we have chosen an approach inspired by the Stroock and Varadhan \textcolor{black}{original arguments} which explicitly provides the required Krylov like estimates.
Before going into the proof, let us briefly explain the main differences between our analysis and the strategy of \cite{stro:vara:79}.
In particular,
our approach differs from the original one because of the specific structure of our problem. }

In the original non degenerate setting \textcolor{black}{with bounded drifts} considered by Stroock and Varadhan, the Girsanov Theorem allows them to deal with the diffusive part of the equation only. Their main idea to obtain the desired control on their perturbed kernel goes through regularization argument\textcolor{black}{s}. The key point allowing them to get the estimation at the limit are: the strong convergence of the \textcolor{black}{driftless} Euler scheme (to keep track on pointwise estimate) and a localization argument. 

In our current setting things are a bit different: we are not allowed anymore to get rid of the drift, because of our degenerate structure. 


\textcolor{black}{Thus, our strategy is the following. In the H\"older framework of case \A{D}-(a), we manage to prove directly the existence of the density of the canonical process for any solution to the martingale problem through the associated Krylov-type estimate, for $f\in L^{q'}([0,T],L^{p'}(\R^{nd}))$ and $p',q'$ \textit{large enough}. This is a consequence of the  pointwise controls established in Lemma \ref{GROS_LEMME} below. This is enough to derive the well posedness of the martingale problem. In a second time, we complete the proof of the Krylov estimate \eqref{KRYLOV_1} on the indicated range for $p',q'$ through a regularization argument. Namely, we regularize the drift coefficient $\gF$ through convolution. For the regularized drift, it follows from \cite{dela:meno:10}, \cite{meno:10} that the corresponding process has a density. It then follows from the previous analysis that the process with mollified coefficients  satisfies uniformly w.r.t. the mollification parameter the Krylov estimate. The final statement then follows letting the mollification parameter tend to zero from the well posedness of the martingale problem. Cases (b) and (c) are handled from case (a) through an additional Girsanov type argument.}

\subsection{Well posedness with full H\"older drift, Assumption \A{D}-(a)}
\label{WP_SEC_HOLDER}

\subsubsection{Existence under Assumption \textcolor{black}{\A{D}{\rm -(a)}}}
The first step is to establish that there exists a solution to the martingale problem defined in Theorem \ref{THM1}.
From the definition of $(L_t)_{t\ge 0} $ in \eqref{DEF_L} it is easily seen that, under \A{A}, existence is obtained adapting to our current framework Theorem 6.1.7 in \cite{stro:vara:79}. The strategy is clear. An Euler like scheme can be considered. Fix first $T>0$ and consider the grid $\Lambda_m([0,T]):=\{(t_i:=ih)_{i\in \leftB 0,\textcolor{black}{m}\rightB}\}, h=T/m $, $\textcolor{black}{m\in \N }$. Introduce the corresponding ``discretization" scheme:
 \begin{align}
\label{EULER_MOD_NODRIFT}
 \X_{\textcolor{black}{s}}^m= \X_{t_i}^m +\int_{t_i}^{s}\Big( \gF(t_i, \X_{t_i}^m)
&+D \gF(t_i, \X_{t_i}^m)( \X_u^{m}-\X_{t_i}^m)\Big)du
\\
+&B\sigma(t_i, \X_{t_i}^m)(W_{s}-W_{t_i}),\ \textcolor{black}{s\in [t_i,t_{i+1}]}.\notag
\end{align}
This scheme defines a sequence of measures $(\P^m)_{m\ge 1} $ on $C([0,T],\R^{nd})$ which is tight and for which the continuity assumption on the coefficients \textcolor{black}{and the sub-linearity of the drift} allow to identify that any limit $\P$ solves the martingale associated with $( L_t)_{t\ge 0} $ on $[0,T]$. We refer to Section 6.1 of \cite{stro:vara:79} for details. To derive the existence of a solution to the martingale problem on the whole positive line, we can rely on a usual chaining \textcolor{black}{in time} argument, see e.g. Chapter 6 in \cite{stro:vara:79}.\\

\subsubsection{Uniqueness \textcolor{black}{under Assumption \A{D}{\rm{-(a)}}}}\label{susec:uniqueness_holdrift}{\color{black}
To establish uniqueness and the Krylov type control of Theorem \ref{THM1},  the key ingredient is to prove that an operator involving  $L$ and a suitable associated perturbation (based on the frozen process/generators of Section \ref{SEC2}) satisfy appropriate estimates. Namely, for $\varepsilon \ge 0 $ introduce:
\begin{eqnarray}\label{eq:def_de_R}
R^\varepsilon f(t,\x)&:=&(L_t \tilde  G^\varepsilon f-\tilde{  M}_t^\varepsilon f)(t,\x)\\
&=&\int_{(t+\varepsilon) \wedge T }^T ds \int_{\R^{nd}}d\y( L_t -\tilde { L}_t^{s,\y})\tilde {p}^{is s,\y}(t,s,\x,\y) f(s,\y),\notag
\end{eqnarray}
with $\tilde G^\varepsilon f, \tilde M^\varepsilon f $ and $\tilde p^{s,\y}$ defined in \eqref{GREEN_KERNEL}, \eqref{DEF_M_EPSILON} and \eqref{CORRESP} respectively. From now on, we also set $Rf(t,\x)=R^0 f(t,\x) $.
 \textcolor{black}{The aforementioned estimates are summarized in} 
the following general Lemma whose proof is postponed to the end of the current section.
\begin{lemma}[\textcolor{black}{Pointwise and} $L^{q'}-L^{p'} $ Control for $R^\varepsilon$]\label{GROS_LEMME}
\phantom{SAUT DE LIGNE}\hspace*{.5cm}\\
\textcolor{black}{There exists  $q_0':= q_0'(\A{A}),\ p_0':= p_0'(\A{A})$ s.t. for $q'\ge q_0', p'\ge p_0' $, it holds that for all $(t,\x)\in [0,T]\times \R^{nd}$ and $R^\varepsilon$ as in \eqref{eq:def_de_R}:
\begin{equation}
\label{CTR_PONCTUEL_RE}
|R^\varepsilon f(t,\x)|\le C \|f\|_{L^{q'}([0,T]\times \R^{nd},L^{p'}(\R^{nd}))},
\end{equation}
with $C:=C(\A{A},T,p_0',q_0') $}.\\


\textcolor{black}{Also,} for  $q',p'>1 $ s.t. $\frac{n^2d}{p'}+\frac{2}{q'}<2$, we have that for all $f\in L^{q'}([0,T]\times \R^{nd},L^{p'}(\R^{nd}))$, uniformly in $\varepsilon\in [0,\varepsilon_0] $, $\varepsilon_0 $ small enough:
\begin{equation}
\label{BD_R}
\|R^\varepsilon f\|_{L^{q'}([0,T]\times \R^{nd},L^{p'}(\R^{nd}))}\le C \|f\|_{L^{q'}([0,T]\times \R^{nd},L^{p'}(\R^{nd}))},
\end{equation}
with $C:=C(\A{A},T,p',q')\underset{T\rightarrow 0}{\longrightarrow} 0 $. In particular, equation \textcolor{black}{\eqref{BD_R}} implies that the operator $I-R^\varepsilon $ is invertible, with bounded inverse in $L^{q'}-L^{p'}$, provided $T$ is small enough.
\end{lemma}
\noindent Having this result at hand, we are now in position to derive uniqueness of the martingale problem.\\

The first step of our approach consists in showing that any solution to the martingale problem satisfies the Krylov like density estimate of Theorem \ref{THM1} for $p'$ and $q'$ therein \emph{large enough} but finite. Indeed, the first condition on $p'$, $q'$ precisely allows us to take benefit from the pointwise control  \eqref{CTR_PONCTUEL_RE} in Lemma \ref{GROS_LEMME}. Note that in this setting, thanks to an approximation argument, the Krylov type estimate \eqref{KRYLOV_1} reduces to
\begin{equation}\label{KRYLOV_POUR_TOUS}
\forall f\in C_0^\infty([0,T)\times \R^{nd},\R),\quad  \Bigg|\E^\P\Bigg[\int_{t}^T f(s,\X_s^{t,\x}) ds\Bigg]\Bigg|\le C(T)\|f\|_{L^{q'}([0,T],L^{p'}(\R^{nd}))}.
\end{equation}
\textcolor{black}{To establish the above control, we apply}
\textcolor{black}{the} It\^o formula on the Green kernel $\tilde G^\varepsilon f$, for $f$ in $ \textcolor{black}{C_0^{\infty}}([0,T)\times \R^{nd},\R)$, for $\varepsilon $ small enough (depending on the support of $f $) and the process $\mathbf X_s^{t, \mathbf x}$:

\begin{eqnarray*}
\tilde G^\varepsilon f(t,\x)+ \E\Bigg[\int_t^{T} (\partial_s+L_s)\tilde G^\varepsilon f(s,\X_s^{t,\x})   ds\Bigg]\\
= \tilde G^\varepsilon f(t,\x)+ \int_t^{T} \int_{\R^{nd}}(\partial_s+L_s)\tilde G^\varepsilon f(s,\y) P_{\X_s^{t,\x}}(d\y)   ds=0,
\end{eqnarray*}
where $P_{\X_s^{t,\x}}$ denotes the law of $\X_s^{t,\x} $. We exploit \eqref{relation differentielle 2} to write:
\begin{align}
\label{NUMBER_PREAL_REG_DENS}
\tilde G^\varepsilon f (t,\x)-& \int_t^{T} \int_{\R^{nd}}f_\varepsilon(s,\y) P_{\X_s^{t,\x}}(d\y)   ds\\
+&\int_{t}^T\int_{\R^{nd}}(L_s \tilde G^\varepsilon  f-\tilde M_s^\varepsilon f)(s,\y) P_{\X_s^{t,\x}}(d\y)   ds=0,\nonumber
\end{align}
\textcolor{black}{where $f_\varepsilon$ is defined in Lemma \ref{convergence_LPLQ}.}
Recall that thanks to Proposition \ref{prop:pointwise_green_kernel} we have that there exists $C(T):=C(\A{A},T) \underset{T\rightarrow 0}{\longrightarrow} 0$ s.t. for all $(t,\x)\in  [0,T]\times \R^{nd} $ and all $f\in L^{q'}([0,T],L^{p'}(\R^{nd})$:

\begin{equation}
\label{CTR_GREEN_PONCT}
|\tilde G^\varepsilon f(t,\x)|\le C(T) \|f\|_{L^{q'}([0,T],L^{p'}(\R^{nd}))}.
\end{equation}

\textcolor{black}{As we assume that $p',q' $ are \emph{large enough}, the pointwise control \eqref{CTR_PONCTUEL_RE} of Lemma \ref{GROS_LEMME} holds. From \eqref{NUMBER_PREAL_REG_DENS}, \eqref{CTR_GREEN_PONCT} and \eqref{CTR_PONCTUEL_RE} we hence readily get:
\begin{equation}
|\int_t^T \E[f_\varepsilon(\X_s^{t,\x})] ds |\le C\|f\|_{L^{q'}([0,T],L^{p'}(\R^{nd}))}.
\end{equation}
\textcolor{black}{Letting $\varepsilon $ go to zero,  we  thus derive from Lemma \ref{convergence_dirac} that, for any solution $\P$ of the martingale problem, provided $p',q'$ are large enough, inequality \eqref{KRYLOV_POUR_TOUS} holds.}\\
}

We are now in position to prove uniqueness to the martingale problem. \textcolor{black}{Let $\P$ be a solution and denote, thanks to the Krylov type estimate established above, its density by $p(t,s,\x,\y) $. We have from \eqref{NUMBER_PREAL_REG_DENS} that: for $f$ in $ \textcolor{black}{C_0^{\infty}}([0,T)\times \R^{nd},\R)$}
\begin{eqnarray}
&&\textcolor{black}{-\tilde G^\varepsilon f(t,\x)}\notag\\
&\textcolor{black}{=}&\textcolor{black}{- \E^{\P}\Bigg[\int_t^{T} f_\varepsilon(s,\X_s^{t,\x}) ds\Bigg]+\E^\P\Bigg[\int_{t}^T R^\varepsilon f(s,\X_s^{t,\x})  ds\Bigg]}\notag\\
&\textcolor{black}{=}&\textcolor{black}{- \int_t^{T} \int_{\R^{nd}}f_\varepsilon (s,\y) p(t,s,\x,\y) d\y  ds+\int_{t}^T\int_{\R^{nd}}R^\varepsilon f(s,\y) p(t,s,\x,\y) d\y  ds}\notag\\
&\textcolor{black}{=}&\textcolor{black}{-\int_t^{T} \int_{\R^{nd}}(I^\varepsilon-R^\varepsilon)f(s,\y) p(t,s,\x,\y) d\y  ds,} \label{REP_PRESQUE_FINI}
\end{eqnarray}
denoting by $I^\varepsilon f(s,\y):=f_\varepsilon (s,\y)$.
\textcolor{black}{It thus easily follows from Proposition \ref{prop:pointwise_green_kernel}, Lemmas \ref{convergence_LPLQ} 
and \ref{GROS_LEMME} eq. \eqref{BD_R}} that both sides are continuous with respect to $L^{q'}-L^{p'} $ norm, uniformly in $\varepsilon $. In particular, \eqref{REP_PRESQUE_FINI} holds for all $f\in L^{q'}-L^{p'} $. From Lemmas  \ref{convergence_LPLQ} and \ref{GROS_LEMME} we can apply  \eqref{REP_PRESQUE_FINI} to $(I^\varepsilon - R^\varepsilon)^{-1}f$ and deduce, still from these Lemmas, that
\begin{equation}
\label{PREAL_KRYLOV_AVANT_INVERSION}
\textcolor{black}{\E^\P\Bigg[\int_t^{T} f(s,\X_s^{t,\x}) ds\Bigg]=\tilde G\circ(I-R)^{-1}f(t,\x),}
\end{equation}
\textcolor{black}{by letting $\varepsilon $ go to zero. This gives uniqueness for $T$ small enough}. \textcolor{black}{Global well-posedness is derived from a chaining in time argument. }\\

To conclude the proof of Theorem \ref{THM1}, it now only remains to derive the Krylov like density estimate \eqref{KRYLOV_1} in full generality \emph{i.e.} for any $p'$, $q'$ satisfying $\frac{n^2d}{p'}+\frac{2}{q'}<2$ (whence not assuming that they are \emph{large enough} and finite).

Assume first that both $p'$ and $q'$ are not \emph{large enough}. \textcolor{black}{We now consider, for a parameter $\delta>0 $, the SDE \eqref{SYST} with drift $\gF^\delta(t,\x):=\big(\gF(t,\cdot)\star \Phi_\delta\big) (\x)$, where $\Phi_\delta(\cdot):=\delta^{-nd}\Phi(\cdot/\delta) $ with $\Phi\in C_0^\infty(\R^{nd},\R^+),\ \int_{\R^{nd}} \Phi(\z)d\z=1 $. It is known that, denoting by $\P^\delta $ the associated solution to the martingale problem, the canonical process enjoys two-sided multi-scale Gaussian bounds similar to those of \eqref{CTR_GRAD} (with $k=0$), see indeed again \cite{dela:meno:10}, \cite{meno:10}, \textcolor{black}{with constants possibly depending on $\delta $}.}

\textcolor{black}{Denoting by $p^\delta(t,s,\x,\y)$ its density, similarly to \eqref{NUMBER_PREAL_REG_DENS} we get for any $f$ in $ \textcolor{black}{C_0^{\infty}}([0,T)\times \R^{nd},\R)$,
\begin{align}
\label{NUMBER_REG_DENS}
\tilde G^{\delta,\varepsilon} f(t,\x)&- \int_t^{T} \int_{\R^{nd}}f_\varepsilon (s,\y) p^\delta(t,s,\x,\y) d\y  ds\\
&+\int_{t}^T\int_{\R^{nd}}(L_s^\delta \tilde G^\varepsilon f-\tilde M_s^{\delta,\varepsilon} f)(s,\y) p^\delta(t,s,\x,\y) d\y  ds=0,\notag
\end{align}
where in the above equation, $\tilde G^{\delta,\varepsilon} f, L_s^\delta, \tilde M_s^{\delta,\varepsilon}$ denote the frozen Green kernel, generator and frozen generator with mollified drift. Importantly, the pointwise bound \eqref{CTR_GREEN_PONCT} on the Green kernel and the controls of Lemma \ref{GROS_LEMME} are uniform with respect to this additional mollifying parameter.}
 
\textcolor{black}{Then, from \eqref{NUMBER_REG_DENS} and Lemma \ref{GROS_LEMME}, we deduce that 
\begin{eqnarray*}
C(T) \|f\|_{L^{q'}([0,T],L^{p'}(\R^{nd}))}\Big(1+\|p^\delta\|_{L^{\tilde q'}([0,T],L^{\tilde p'}(\R^{nd}))}\Big)\\
\ge \bigg|\int_{t}^{T}\int_{\R^{nd}}f_\varepsilon(s,\y) p^\delta(t,s,\x,\y) d\y ds\bigg|,
\end{eqnarray*}
where $\tilde q',\tilde p' $ are the conjugate exponents of $q',p'$ respectively}.

\textcolor{black}{Using now Lemma \ref{convergence_LPLQ} we get, for any $f  \in C_0^\infty([0,T)\times \R^{nd},\R)$
\begin{align*}
&\bigg|\int_{t}^{T}\int_{\R^{nd}}f(s,\y) p^\delta(t,s,\x,\y) d\y ds\bigg|\\
\le & C(T) \|f\|_{L^{q'}([0,T],L^{p'}(\R^{nd}))}\Big(1+\|p^\delta\|_{L^{\tilde q'}([0,T],L^{\tilde p'}(\R^{nd}))}\Big)\\
&+\|f-f_\varepsilon\|_{L^{q'}([0,T],L^{p'}(\R^{nd}))}\|p^\delta\|_{L^{\tilde q'}([0,T],L^{\tilde p'}(\R^{nd}))}\\
\le & C(T) \|f\|_{L^{q'}([0,T],L^{p'}(\R^{nd}))}\Big(1+\|p^\delta\|_{L^{\tilde q'}([0,T],L^{\tilde p'}(\R^{nd}))}\Big)\\
&+\frac 12\|f\|_{L^{q'}([0,T],L^{p'}(\R^{nd}))}\|p^\delta\|_{L^{\tilde q'}([0,T],L^{\tilde p'}(\R^{nd}))}
\end{align*}
as soon as $\varepsilon $, possibly depending on $f$, is small enough. Up to an additional density argument,}
\textcolor{black}{this  yields in particular from the Riesz representation theorem that, for $T$ small enough,  $\|p^\delta\|_{L^{\tilde q'}([0,T],L^{\tilde p'}(\R^{nd}))}\le C(T) $. Thus, for any $f$ in $L^{q'}([0,T],L^{p'}(\R^{nd}))$:
\begin{align*}
&\Bigg|\int_t^{T} \int_{\R^{nd}}f (s,\y) p^\delta(t,s,\x,\y)  d\y ds\Bigg|=\Bigg|\int_t^{T}\E^{\P^\delta} [f (s,\X_s^{t,\x})]  ds\Bigg|\\
 \le& C(T)\|f\|_{L^{q'}([0,T],L^{p'}(\R^{nd}))}.
\end{align*}
Estimate \eqref{KRYLOV_POUR_TOUS} is then obtained 
letting 
 $\delta $ go to zero from the well posedness of the martingale problem and Theorem 11.1.4 in \cite{stro:vara:79}}.

Next, the case when $q'=p'=\infty$ is direct and the remaining cases  are thus $q'<+\infty,\, p'=\infty $ and $q'=\infty,\, p'<+\infty $. The first one is again direct. The second follows from the fact that for any $q'\in (2+\infty]$, $\|f\|_{L^{q'}([0,T],L^{p'}(\R^{nd}))}\le \|f\|_{L^\infty([0,T],L^{p'}(\R^{nd}))} T^{1/q'}$. \qed \\
}

\vspace*{.2cm}
\paragraph{\textcolor{black}{\textbf{Proof of Lemma \ref{GROS_LEMME}}}}
We focus on the proof for $ Rf=R^0f$. The parameter $\varepsilon $ does not play here any role for the estimates. We have, by definition

\begin{eqnarray}
Rf(t,\x) &=& \int_t^Tds \int_{\R^{nd}} d\y ( L_t - \tilde{L}^{s,\y}_t )\tilde{p}^{s,\y}(t,s,\x,\y)f(s,\y)\nonumber\\
&=:&\int_t^T ds \int_{\R^{nd}} d\y H(t,s,\x,\y)\textcolor{black}{f}(s,\y).\label{DEF_I2}
\end{eqnarray}
Here, the operator $H$ is the so-called \textcolor{black}{backward }\textit{parametrix kernel}, \textcolor{black}{see \cite{mcke:sing:67}}. It already appeared, in a similar form but under stronger smoothness assumptions in \cite{dela:meno:10,meno:10}.

The bound \eqref{CTR_GRAD} of Proposition \ref{THE_PROP} now yields that there exists $C:=C(\A{A})$ such that:
\begin{align*}
&|Rf(t,\x) |\\
\leq&\int_t^T ds\int_{\R^{nd}} d\y |f(s,\y)| \Bigg\{\Big|F_1(t,\x)-F_1(t,\btheta_{t,s}(\y))\Big| |D_{x_1}\tilde p^{s,\y}(t,s,\x,\y)|\\
&+\sum_{i=2}^n \bigg\{\Big|F_i(t,\x)-\Big(F_i(t,\btheta_{t,s}(\y))- D_{x_{i-1}}F_i(t,\btheta_{t,s}(\y))(\x-\btheta_{t,s}(\y))_{i-1}   \Big)\Big|\\
&\times  |D_{x_i}\tilde p^{s,\y}(t,s,\x,\y)|\bigg\}+|a(t,\x)-a(t,\btheta_{t,s}(\y))| |D_{x_1}^2\tilde p^{s,\y}(t,s,\x,\y)|\Bigg\}\\
\le&C \int_t^T ds\int_{\R^{nd}} d\y |f(s,\y)| \Bigg\{\frac{\Big(1+|\x-\btheta_{t,s}(\y)|\Big)}{(s-t)^{\frac 12}}+\sum_{i=2}^n \bigg\{\Big[\Big|F_i(t,\x)-F_i(t,x_{i-1},\btheta_{t,s}(\y)^{i:n})\Big|  \\
&+\Big|F_i(t,x_{i-1},\btheta_{t,s}(\y)^{i:n})-\Big(F_i(t,\btheta_{t,s}(\y))- D_{x_{i-1}}F_i(t,\btheta_{t,s}(\y))(\x-\btheta_{t,s}(\y))_{i-1}   \Big)\Big|\Big]\\
&\times \frac{1}{(s-t)^{(i-1)+\frac 12}}\bigg\}+\frac{|\x-\btheta_{t,s}(\y)|^\eta}{s-t}\Bigg\} 
\frac{\exp\left(-C^{-1}(s-t)|\T_{s-t}^{-1}(\x-\btheta_{t,s}(\y))|^2 \right)}{ (s-t)^{\frac{n^2d}{2}}  },
\end{align*}
where we have denoted for $\z\in \R^{nd}, \z^{i:n}=(z_i,\cdots, z_n)\in \R^{(n-(i-1))d} $. From \A{A} and \eqref{COND_HOLDER} we thus derive, up to a modification of $C$:
\textcolor{black}{
\begin{align*}
| Rf(t,\x)|
&\le C \int_t^T ds\int_{\R^{nd}} d\y |f(s,\y)| \Bigg\{\Big((s-t)^{-\frac 12}+(s-t)^{\frac 12}|\T_{s-t}^{-1}(\x-\btheta_{t,s}(\y))|\Big)\\
&+\sum_{i=2}^n \bigg\{\Big[\sum_{j=i}^n\Big(\frac{|(\x-\btheta_{t,s}(\y))_j|}{(s-t)^{j-\frac 12}}\Big)^{\beta_i^j} (s-t)^{\beta_i^j(j-\frac 12)} \\
 &+\Big(\frac{|(\x-\btheta_{t,s}(\y))_{i-1}|}{(s-t)^{(i-1)-\frac 12}} \Big)^{1+\eta}   (s-t)^{(i-\frac 32)(1+\eta)}\Big|\Big]\frac{1}{(s-t)^{i-\frac 12}}\bigg\}\\
&+\frac{\big((s-t)^{\frac 12}|\T_{s-t}^{-1}(\x-\btheta_{t,s}(\y))|\big)^\eta}{(s-t)^{1-\frac \eta 2}}\Bigg\}\frac{\exp\left(-C^{-1}(s-t)|\T_{s-t}^{-1}(\x-\btheta_{t,s}(\y))|^2 \right)}{ (s-t)^{\frac{n^2d}{2}}  }.
\end{align*}
}
\textcolor{black}{
In the above equation we put each contribution coming from the difference of the generators at its \textit{intrinsic scale} w.r.t. exponential bounds. In other words, the terms: 
\begin{eqnarray*}
 \Big(\frac{|(\x-\btheta_{t,s}(\y))_j|}{(s-t)^{j-\frac 12}}\Big)^{\beta_i^j}&\le& \big((s-t)^{\frac 12}|\T_{s-t}^{-1}(\x-\btheta_{t,s}(\y))| \big)^{\beta_i^j}, \\
 \Big(\frac{|(\x-\btheta_{t,s}(\y))_{i-1}|}{(s-t)^{(i-1)-\frac 12}} \Big)^{1+\eta}&\le& \big((s-t)^{\frac 12}|\T_{s-t}^{-1}(\x-\btheta_{t,s}(\y))| \big)^{(1+\eta)},
 \end{eqnarray*} can be absorbed by the exponential.}
Therefore,
\begin{eqnarray}
| Rf(t,\x)|&\leq & C\Bigg\{\sum_{i=2}^n\sum_{j=i}^n \int_t^T \frac{ds}{(s-t)^{i-\frac 12 - \textcolor{black}{\beta_i^j(j-\frac12)}}} \int_{\R^{nd}}  d\y  
{\bar
 p}_{C^{-1}}(t,s,\x,\y)  |f(s,\y)| \nonumber\\
&&+ \int_t^T \frac{ds}{(s-t)^{1-\frac \eta 2}} \int_{\R^{nd}}d\y    {\bar p}_{C^{-1}}(t,s,\x,\y)  |f(s,\y)|\notag\\
&&+\sum_{i=2}^n\int_t^T \frac{ds}{(s-t)^{\textcolor{black}{1-\eta(i-\frac 3 2)}}} \int_{\R^{nd}}  d\y  {\textcolor{black}{\bar p} }_{C^{-1}}(t,s,\x,\y)  |f(s,\y)|\Bigg\}. \label{BD_Rprime}
 \end{eqnarray}
Note carefully that the condition \eqref{COND_HOLDER} precisely gives that for all $i\in \leftB 2,n\rightB, j\in \leftB i,n\rightB $, $1-\{(i-1)+\frac 12-\beta_i^j(j-1+\frac 12)> 0 $ so that all the above time singularity are integrable. Note also that, thanks to Lemma \ref{lemme:bilipflow} \textcolor{black}{(almost equivalence on the flows)}:
\begin{equation}
\label{CTR_SP_INT}
\int_{\R^{nd}} d\y   {\bar p}_{C^{-1}}(t,s,\x,\y) = \int_{\R^{nd}} d\y\frac{\exp\left(-C^{-1}(s-t)|\T_{s-t}^{-1}(\x-\btheta_{t,s}(\y))|^2 \right)}{ (s-t)^{\frac{n^2d}{2}}  }\le C_{\ref{CTR_SP_INT}}.
\end{equation}

The result\textcolor{black}{s in \eqref{CTR_PONCTUEL_RE} and \eqref{BD_R} now follow} from \eqref{BD_Rprime} and the following key Lemma.
\begin{lemma}[$L^{q'}-L^{p'} $ Controls for the singularized Green kernel]\label{CTR_N_GAMMA}
Introduce for $f\in L^{q'}\big([0,T],L^{p'}(\R^{nd})\big)$,
\textcolor{black}{with $p',q'>1$} and some $\textcolor{black}{\gamma\in [0,1)}$ the quantify:
\begin{equation}
\label{CTR_RESTE_GENE}
 N_\gamma f(t,\x) :=   \int_t^T \frac{ds}{(s-t)^{\gamma}} \int_{\R^{nd}} d\y  {\bar p}_{C^{-1}}(t,s,\x,\y)  |f(s,\y)|.
\end{equation}
\textcolor{black}{There exists  $q_0':= q_0'(\gamma)\ge 1,\ p_0':= p_0'(\gamma)\ge 1 $ and $C:=C(\A{A},T,p_0',q_0') $ s.t. for $q'\ge q_0', p'\ge p_0' $, it holds that for all $(t,\x)\in [0,T]\times \R^{nd}$:
\begin{equation*}
N_\gamma  f(t,\x)\le C \|f\|_{L^{q'}([0,T]\times \R^{nd},L^{p'}(\R^{nd}))}.
\end{equation*}
}

\textcolor{black}{Also,} there exists $C(T):=C(T,\A{A},\gamma) \underset{T\rightarrow 0}{\longrightarrow} 0$ s.t. for all $f\in L^{q'}\big([0,T],L^{p'}(\R^{nd})\big) $,
$$ \|N_\gamma f \|_{L^{q'}([0,T],L^{p'}(\R^{nd}))}\le C(T) \|f\|_{L^{q'}([0,T],L^{p'}(\R^{nd}))}.$$
\end{lemma}
\qed

\begin{proof}[Proof of Lemma \ref{CTR_N_GAMMA}]
\textcolor{black}{Let us first start with the pointwise estimate.} \textcolor{black}{For all $(t,\x)\in [0,T)\times \R^{nd} $, denoting by $\tilde p',\tilde q' $ the conjugate exponents of $p',q'$, write:
\begin{eqnarray*}
 N_\gamma f(t,\x) &\le&    \Big(\int_t^T \frac{ds}{(s-t)^{\tilde q'\gamma}} \big(\int_{\R^{nd}} d\y  {\bar p}_{C^{-1}}(t,s,\x,\y)^{\tilde p'}\big)^{\frac{\tilde q'}{\tilde p'}}  \Big)^{\frac 1{\tilde q'}}\|f\|_{L^{q'}([0,T],L^{p'}(\R^{nd}))}\\
 &\le &C\int_t^T \frac{ds}{(s-t)^{\tilde q'(\gamma+\frac{n^2d}2(1-\frac 1{\tilde p'}))}} \|f\|_{L^{q'}([0,T],L^{p'}(\R^{nd}))}.
\end{eqnarray*}
It is thus clear that the \textit{pointwise} statement of the lemma is fulfilled if $p',q'$ are \textit{large} enough, in order to guarantee that, for $(\tilde p')^{-1}+(p')^{-1}=1,\ (\tilde q')^{-1}+(q')^{-1}=1 $,  $\tilde p',\tilde q' $ are sufficiently close to 1 to have $\tilde q'(\gamma+\frac{n^2d}2(1-\frac{1}{\tilde p'}))<1 $.}

\textcolor{black}{Let us now turn to the estimate in $L^{q'}([0,T],L^{p'}(\R^{nd})) $ norm}. Setting 
$$K_{t,s,\gamma}\star f(s,\x):=\frac{1}{(s-t)^\gamma}\int_{\R^{nd}} d\y \bar p_{C^{-1}}(t,s,\x,\y) f(s,\y) d\y,$$ we derive from \eqref{CTR_RESTE_GENE} and the triangle inequality:
 \begin{equation}
\label{EQ_N_GAMMA}
 \Big \| N_\gamma f(t,\cdot)\Big\|_{L^{p'}(\R^{nd})}=\Big \|\int_t^T ds K_{t,s,\gamma}\star f(s,\cdot)\Big\|_{L^{p'}(\R^{nd})}\le \int_t^T ds \|K_{t,s,\gamma}\star f(s,\cdot)\|_{L^{p'}(\R^{nd})}.
 \end{equation}
From usual $L^1-L^{p'} $ convolution arguments we also get:
\begin{eqnarray*}
\|K_{t,s,\gamma}\star f(s,\cdot)\|_{L^{p'}(\R^{nd})}\le \frac{C}{(s-t)^\gamma} \|f(s,\cdot)\|_{L^{p'} (\R^{nd})}.
\end{eqnarray*} 
Plugging this estimate into \eqref{EQ_N_GAMMA} and using the H\"older inequality with exponents $q',\tilde {q}'>1 $ s.t.  $(q')^{-1}+(\tilde {q}')^{-1}=1$ we obtain:
\begin{eqnarray*}
\int_0^T dt \| N_\gamma f(t,\cdot)\|_{L^{p'} (\R^{nd})}^{q'}&\le& C\int_0^T dt \Big(\int_t^T \frac{ds}{(s-t)^\gamma} \|f(s,\cdot)\|_{L^{p'} (\R^{nd})}\Big)^{q'}\\
&\le &C_T \int_0^T dt \int_t^T ds \frac{C}{(s-t)^{\gamma}}\|f(s,\cdot)\|_{L^{p'}(\R^{nd})}^{q'},\\
\end{eqnarray*}
where $C_T:=C(\A{A},p',q',T)$ denotes a constant which is small as $T$ is. From the Fubini Theorem we eventually derive
\begin{eqnarray*}
\int_0^T dt \| N_\gamma f(t,\cdot)\|_{L^{p'} (\R^{nd})}^{q'}&\le&C_T
\int_0^T ds \|f(s,\cdot)\|_{L^{p'}(\R^{nd})}^{q'} \int_0^s dt\frac{C}{(s-t)^{\gamma}}\\
&\le& C_T \|f\|_{L^{q'}([0,T],L^{p'}(\R^{nd}))}^{q'},
\end{eqnarray*}
up to a modification of $C_T$ in the last inequality.
\end{proof}

Hence, using the above control  in \textcolor{black}{\eqref{BD_Rprime}}  yields:
$$\| Rf\|_{L^{q'}([0,T],L^{p'}( \R^{nd}))}\le C_T\|f\|_{L^{q'}([0,T],L^{p'}( \R^{nd}))}.$$
 \textcolor{black}{This concludes the proof of Lemma \ref{GROS_LEMME}. \qed}

\subsection{Existence and uniqueness under \A{D}-(b) and (c)} \label{EX_D_B_C}
Since \textcolor{black}{in that setting} no continuity is assumed on $F_1$ we will derive the well posedness through Girsanov arguments. It therefore clearly appears that the dynamics with 0 drift on the first non-degenerate component, and its associated linearization, plays a key role. 
We first introduce some notations used throughout this section. 

We first define for all $\x\in  \R^{nd},\ \bar \gF(t,\x):=({\mathbf 0}, F_2(t,\x),\cdots,F_n(t,\x))$, where $\bar \gF $ satisfies \A{S} and \A{H}.
\textcolor{black}{Recall from} Theorem \ref{THM1} \textcolor{black}{that}, under \A{A}-(a),  weak uniqueness holds for the SDE:
\begin{equation}\label{DYN_NO_DRIFT}
d\bar  \X_t=\textcolor{black}{\bar \gF}(t,\bar \X_t)dt+B\sigma(t,\bar \X_t)dW_t.
\end{equation}

For fixed $(T,\y)\in \R^{+} \times \R^{nd}$, we consider the following deterministic system to define our Gaussian proxy:
\begin{equation}
\label{SYST_DIFF_MOD_B}
\overset{.}{\bar \btheta}_{t,T}(\y)=\bar \gF(t,\bar \btheta_{t,T}(\y)),\ \bar \btheta_{T,T}(\y)=\y,\ 
\end{equation} 
and 
\begin{equation*}
\frac{d}{dt} \tilde{\bar \bphi}_t = \bar {\mathbf F}(t,\bar \btheta_{t,T}(\y))
+ D {\mathbf F}(t,\bar \btheta_{t,T}(\y))[\tilde{\bar \bphi}_t
- \bar \btheta_{t,T}(\y)], \quad t \geq 0.
\end{equation*}
Again, \textcolor{black}{in \eqref{SYST_DIFF_MOD_B}, we consider \textcolor{black}{a} Cauchy-Peano flow furnished by  Lemma \ref{lemme:measu} and which also satisfies the equivalence of rescaled norms of Lemma \ref{lemme:bilipflow}}.
The dynamics of the linearized Gaussian process associated with \eqref{DYN_NO_DRIFT} writes:
\begin{equation}
\label{LIN_SPEC1}
 d\tilde{\bar { \X}}_t^{T,\y}= \Big(\bar \gF(t,\bar \btheta_{t,T}(\y))+D\bar \gF(t,\bar \btheta_{t,T}(\y))(\tilde {\bar \X}_t^{T,\y}-\bar \btheta_{t,T}(\y))\Big) dt+B\sigma(t,\bar \btheta_{t,T}(\y))dW_t,
\end{equation} 
and we denote the associated generator by $(\tilde {\bar L}_t^{T,\y})_{t\in [0,T]} $ and by $\tilde {\bar p}^{T,\y}(t,s,\x,\cdot) $ the corresponding density at times $s>t$ when the process starts in $\x$ at time $t$. 

We point out that, with respect to the previously used notations, we choose to keep track of the driftless dynamics for $F_1$ adding bars on the associated objects: dynamics, generators, density.

For our strategy, recall that we aim at proving uniqueness for the initial SDE \eqref{SYST} through the well posedness of the martingale problem associated with $\textcolor{black}{(L_t)_{t \ge 0}} $. Once existence is known, the point is that we use a different Gaussian proxy than previously, namely the one considered in \eqref{LIN_SPEC1} associated with the driftless dynamics on the first component.

\subsubsection{Existence and Uniqueness under Assumption \A{D}-\textcolor{black}{{\rm (b)}}}
\label{EX_UN_D_B}
Under (b) (bounded measurable drift $F_1 $ on the non-degenerate component), existence is a direct consequence of the Girsanov theorem. We thus now focus on uniqueness. 

Repeating the previous approach (see subsection \ref{susec:uniqueness_holdrift}), using the family of random variables $(\tilde {\bar \X}_s^{s,\y})_{s\in[t,T]} $ defined in \eqref{LIN_SPEC1} with $\tilde {\bar \X}_t^{s,\y}=\x$ as Gaussian proxys, we have to bound analogously to the estimate of Lemma \ref{GROS_LEMME}: 
\begin{eqnarray}
\label{R_PROXY_MODIF}
\hspace{.9cm} Rf(t,\x)&:=&\int_t^{T} ds \int_{\R^{nd}} d\y (L_t-\tilde {\bar L}_t^{s,\y})\tilde {\bar p}^{s,\y}(t,s,\x,\y) f(s,\y)\\
 &=&\int_t^{T} ds \int_{\R^{nd}} d\y \Big\{(L_t-\textcolor{black}{\bar L_t})+(\bar L_t-\tilde {\bar L}_t^{s,\y})\Big\}\tilde {\bar p}^{s,\y}(t,s,\x,\y) f(s,\y)\notag\\
 &=:& \bar {\mathcal R}f(t,\x)+ \bar Rf(t,\x),\notag
\end{eqnarray}
where 
$$ \bar Rf(t,\x):=\int_t^{T} ds \int_{\R^{nd}} d\y (\bar L_t-\tilde {\bar L}_t^{s,\y})\tilde {\bar p}^{s,\y}(t,s,\x,\y) f(s,\y),$$
and 
\begin{align}
\label{BAR_R_CURSIF}
 \bar {\mathcal R}f(t,\x):=&(R-\bar R)f(t,\x)=\int_t^{T} ds \int_{\R^{nd}} d\y (L_t-\bar  L_t)\tilde {\bar p}^{s,\y}(t,s,\x,\y) f(s,\y)\\
=& \int_t^T ds\int_{\R^{nd}} d\y \langle F_1(t,\x),D_{x_1}\tilde {\bar p}^{s,\y}(t,s,\x,\y)\rangle f(s,\y).\notag
\end{align}
From Lemma \ref{GROS_LEMME}, we have already shown that $ \bar Rf(t,\x)$ is controlled in $L^{q'}-L^{p'} $ norm. It thus suffice to investigate the behavior of the $L^{q'}-L^{p'}  $ norm of $\bar {\mathcal R}f(t,\x)$ defined by \eqref{BAR_R_CURSIF}. Namely, our goal is to prove that
\begin{equation}
\label{CTR_TO_BE_ACHIEVED_GIRS}
\|\bar {\mathcal R}f\|_{L^{q'}([0,T],L^{p'}(\R^{nd}))}\le C(T) \|f\|_{L^{q'}([0,T],L^{p'}(\R^{nd}))},\ {\rm with}\ C(T)\underset{T\rightarrow 0}{\longrightarrow} 0.
\end{equation}

Since from \eqref{CTR_GRAD} and \eqref{BAR_R_CURSIF} we have for all $(t,\x)\in [0,T]\times \R^{nd}$:
\begin{align*}
|\bar {\mathcal R}f(t,\x)|\le  C \int_t^T \frac{ds}{(s-t)^{\frac 12}}\int_{\R^{nd}} d\y\frac{\exp\left(-C^{-1}(s-t)|\T_{s-t}^{-1}(\x-\bar \btheta_{t,s}(\y))|^2 \right)}{ (s-t)^{\frac{n^2d}{2}}  }|f(s,\y)|.
\end{align*}
The estimate \eqref{CTR_TO_BE_ACHIEVED_GIRS} then readily follows from Lemma \ref{CTR_N_GAMMA}.


\subsubsection{Existence and Uniqueness under Assumption \A{D}-{\textcolor{black}{{\rm(c)}}}}
We choose in this paragraph to address first the uniqueness, which is a rather direct extension of our previous approach, whereas the existence is a bit involved and requires to exploit the Krylov like inequality \eqref{KRYLOV_1} that has been established for the process $(\bar \X_t)_{t\ge 0} $ \textcolor{black}{with 0 drift in the non-degenerate component} introduced in \eqref{DYN_NO_DRIFT}.

\paragraph{\textbf{\textcolor{black}{Uniqueness under Assumption \A{D}-(c)}}}

With the notations of the previous paragraph it remains to control, in $L^{q'}-L^{p'} $ norm, the contribution ${\bar {\mathcal R}}f $ introduced in \eqref{BAR_R_CURSIF}. The term $\bar Rf$ in \eqref{R_PROXY_MODIF} is again controlled as under assumption \textcolor{black}{\A{D}-(a)}.
Similarly to the previous paragraph (see also \eqref{CTR_GRAD} and \eqref{BAR_R_CURSIF}), we have:
\begin{eqnarray}
\label{CTR_GRAD_LP}
|D_{x_1}\tilde {\bar p}^{\textcolor{black}{s},\y}(t,\textcolor{black}{s},\x,\y)|
&\le& \frac{C}{(\textcolor{black}{s}-t)^{\frac{n^2d}{2}+\frac 12}}
\times\exp\left(-\frac{C^{-1}}{2}(\textcolor{black}{s}-t)|\T_{\textcolor{black}{s}-t}^{-1}(\x-\bar \btheta_{t,\textcolor{black}{s}}(\y))|^2 \right)\notag\\
&\le & \frac C{(s-t)^{\frac 12}} \bar p_{C^{-1}}(t,s,\x,\y).
\end{eqnarray}
Uniqueness then follows from the following lemma which can be viewed as a refinement of Lemma \ref{CTR_N_GAMMA} and explicitly exploits the condition on $p$ and $q$ stated in \textcolor{black}{\A{D}-(c)}.
\begin{lemma}[Refined $L^{q}-L^p$ control of singularized Green kernels]\hspace*{2pt}\phantom{GRRRRR}\\
\label{LEMME_U_F1}
Introduce for all  $f,F_1\in L^q\big([0,T],L^p(\R^{nd})\big) $ and $(t,\x)\in[0,T]\times \R^{nd}$,
\begin{equation}\label{DEF_UF1_FIRST}
u_{F_1}(t,\x):=F_1(t,\x)\int_t^T \frac{ds}{(s-t)^{\frac 12}} \int_{\R^{nd}} \textcolor{black}{\bar p_{C^{-1}}} (t,s,\x,\y) f(s,\y) d\y.
\end{equation}
Then,  for all $p\ge 2, q>2 $ s.t. $\frac {n^2d}p+\frac{2}q<1 $,  there exists $C(T):=C(T,\A{A},p,q) \underset{T\rightarrow 0}{\longrightarrow} 0$ s.t. for all $f,F_1\in L^q\big([0,T],L^p(\R^{nd})\big) $
$$ \|u_{F_1}\|_{L^q([0,T],L^p(\R^{nd}))}\le C(T) \|F_1\|_{L^q([0,T],L^p(\R^{nd}))}\|f\|_{L^q([0,T],L^p(\R^{nd}))}.$$
\end{lemma}
\qed

\begin{proof}[Proof of Lemma \ref{LEMME_U_F1}]
With the notations of Lemma \ref{CTR_N_GAMMA} rewrite:
\begin{align}\label{DEF_UF1}
u_{F_1}(t,\x)&=F_1(t,\x)\int_t^T \frac{ds}{(s-t)^{\frac 12}} \int_{\R^{nd}} \bar p_{C^{-1}}(t,s,\x,\y) f(s,\y) d\y\notag\\
&=: F_1(t,\x)\int_t^T ds K_{t,s,\frac 12}\star f(s,\x).
\end{align}
The triangle inequality yields:
\begin{equation}
\label{NEW_THING}
\Big \|\int_t^T ds  F_1(t,\cdot)K_{t,s,\frac 12}\star f(s,\cdot)\Big\|_{L^p(\R^{nd})}\le \int_t^T ds \|F_1(t,\cdot)K_{t,s,\frac 12}\star f(s,\cdot)\|_{L^p(\R^{nd})}.
\end{equation}
The idea is here to reproduce the computations of Lemma \ref{CTR_N_GAMMA} integrating directly the singularized heat-kernel, i.e.  $K_{t,s,\frac 12}$, in the  $\y$ variable when performing the  H\"older inequality in order to make the product of the norms  $\|f(s,\cdot)\|_{L^p(\R^{nd})}\|F_1(t,\cdot)\|_{L^p(\R^{nd})} $ appear. Precisely:
\begin{align}
&\|F_1(t,\cdot)K_{t,s,\frac 12}\star f(s,\cdot)\|_{L^p(\R^{nd})}\notag\\
=& \Big(\int_{\R^{nd}}d\x |F_1(t,\x)|^p\textcolor{black}{\Big|}\int_{\R^{nd}} d\y K_{t,s,\frac 12}(\x,\y)f(s,\y) \textcolor{black}{\Big|^p} \Big)^{\frac 1p}\nonumber\\
\le& \Big(\int_{\R^{nd}}d\x |F_1(t,\x)|^p\Big\{\int_{\R^{nd}} d\y |f(s,\y)|^p \Big\}\Big\{\int_{\R^{nd}} d\y |K_{t,s,\frac 12}(\x,\y)|^{\tilde p} \Big\}^{\frac {p}{\tilde p}} \Big)^{\frac 1p}\notag,
\end{align}
where again $p^{-1}+(\tilde p)^{-1}=1 $. Observe now that usual Gaussian calculations give that \textcolor{black}{there exists $C_p>0$ such that, for all $\x$ in $\R^{nd}$}:
$$ \Big\{\int_{\R^{nd}} d\y |K_{t,s,\frac 12}(\x,\y)|^{\tilde p} \Big\}^{\frac {1}{\tilde p}}\le \frac{C_p}{(s-t)^{\frac{n^2 d}2\frac 1p+\frac 12}}.$$
Thus,
\begin{align}
\|F_1(t,\cdot)K_{t,s,\frac 12}\star f(s,\cdot)\|_{L^p(\R^{nd})} &\le & \|F_1(t,\cdot)\|_{L^p (\R^{nd})}\|f(s,\cdot)\|_{L^p (\R^{nd})}\frac{C_p}{(s-t)^{\frac{n^2 d}2\frac 1p+\frac 12}},\label{PREAL_HOLDER_UNIQUENESS_C}
\end{align}
which yields from \eqref{NEW_THING} that:
\begin{align*}
\Big \|\int_t^T &ds  F_1(t,\cdot)K_{t,s,\frac 12}\star f(s,\cdot)\Big\|_{L^p(\R^{nd})}\\
&\le C_p \int_t^T  \frac{ds}{(s-t)^{\frac{n^2 d}2\frac 1p+\frac 12}}\|F_1(t,\cdot)\|_{L^p (\R^{nd})}\|f(s,\cdot)\|_{L^p (\R^{nd})} .
\end{align*}
From the definition in \eqref{DEF_UF1} we eventually derive that:
\begin{align}
&\int_{0}^T dt \|u_{F_1}(t,\cdot)\|_{L^p(\R^{nd})}^q \notag\\
&\le  \textcolor{black}{C_p^q}\int_{0}^T \|F_1(t,\cdot)\|_{L^p (\R^{nd})}^{\textcolor{black}{q}} \Bigg(\int_{t}^T \frac{ds}{(s-t)^{(\frac {n^2d}2\frac 1p+\frac 12)\tilde q}} \Bigg)^{\frac q{\tilde q}}\int_{t}^T ds\|f(s,\cdot)\|_{L^p(\R^{nd})}^q\notag\\
&\le \textcolor{black}{C_p^q} \|F_1\|_{L^q([0,T],L^p(\R^{nd}))}^q\|f\|_{L^q([0,T],L^p(\R^{nd}))}^q \Bigg(\int_{t}^T \frac{ds}{(s-t)^{(\frac {n^2d}2\frac 1p+\frac 12)\tilde q}} \Bigg)^{\frac q{\tilde q}},\label{END_PROOF_UF1}
\end{align}
with $q^{-1}+(\tilde q)^{-1}=1  $. Let us now show that the remaining time integral in the above equation gives a small contribution in times. To do so, it suffices to show that   $(\frac {n^2d}2\frac 1p+\frac 12)\tilde q<1 $. Since $(\tilde q)^{-1}=1-\frac 1q $, we have that:
$$ \Big(\frac {n^2d}2\frac 1p+\frac 12\Big)\tilde q<1 \iff \frac {n^2d}2\frac 1p+\frac 12<1-\frac 1q \iff \frac {n^2d}p+\frac 2q<1,$$
which is precisely the condition appearing in \textcolor{black}{\A{D}-(c) and assumed in the current Lemma.}
\end{proof}

\paragraph{\textcolor{black}{\textbf{Existence under Assumption \A{D}-(c)}}}
We here consider a function $F_1 \in L^q([0,T],$ 
$L^p(\R^{nd}))$, where $p,q $ are as in \A{D}-(c). To prove the existence, the strategy is to \textcolor{black}{exploit} the idea introduced by Portenko \cite{port:90} and used by Krylov and R\"ockner \cite{kryl:rock:05} as well to build local weak solutions (before they also establish that they are actually strong solutions) in the non-degenerate case. We also refer for perturbed degenerate Ornstein-Uhlenbeck dynamics to \cite{fedr:flan:prio:vove:17} and \cite{zhan:16}. We adapt a bit this approach.\\

Recall that for the process $\textcolor{black}{(\bar \X_t)_{t\ge 0}}$ introduced in \eqref{DYN_NO_DRIFT} we have from Theorem \ref{THM1}, equation \eqref{KRYLOV_1}, the following 
density estimate.

 Denoting by $\bar P(t,s,\x,\cdot) $ the transition probability determined by $(\bar L_s)_{s\ge 0} $, it is s.t. for a given $T>0$, almost all $\textcolor{black}{s\in (t,T]}$ and all $\Gamma \in {\mathcal B}(\R^{nd}) $:
$\bar P(t,s,\textcolor{black}{\x},\Gamma)=\int_{\Gamma}\bar  p(\textcolor{black}{t,s},\x,\y)d\y$. More specifically,  for any $f\in  L^{q'}([0,T],L^{p'}( \R^{nd})), \frac{n^2d}{p'}+\frac {\textcolor{black}{ 2}}{q'}<2, p'>1,q'>1,$ \textcolor{black}{and $(t,\x)\in [0,T]\times \R^{nd} $}:
\begin{equation}
\label{EST_DENS_SOL_PB_MART_NO_DRIFT}
\Big |\bar \E^{\textcolor{black}{\bar \P_{t,\x}}}\Big[\int_{\textcolor{black}{t}}^Tf(s,\bar \X_s)ds \Big]\Big|\le \textcolor{black}{C_{\ref{EST_DENS_SOL_PB_MART_NO_DRIFT}}}
\|f\|_{L^{q'}([0,T],L^{p'}( \R^{nd}))},\ \textcolor{black}{C_{\ref{EST_DENS_SOL_PB_MART_NO_DRIFT}}}:=\textcolor{black}{C_{\ref{EST_DENS_SOL_PB_MART_NO_DRIFT}}}(\A{A},p',q',\textcolor{black}{T}).
 \end{equation}

We now state an exponential integrability result for the unique weak solution of \eqref{DYN_NO_DRIFT}. Such types of estimates were first proved by Khas'minskii in \cite{khas:59} for the Brownian motion. \textcolor{black}{We can also refer to Lemma 2.1 in Chapter 1 of the monograph by Sznitman \cite{szni:98}. Since the proof only relies on the Markov property,  it readily extends to the current inhomogeneous and non-Brownian framework}.

\begin{lemma}[Khas'minskii's type exponential integrability]\label{LEMME_KHAS}
Let $(\bar \X_t)_{t\ge 0}$ be the (unique-weak) solution to \eqref{DYN_NO_DRIFT}. Then, \textcolor{black}{for any fixed $T>0$ and  a positive Borel function $f:[0,T]\times \R^{nd}\rightarrow \R_+ $}  s.t.
$$\alpha_{\textcolor{black}{T}}:=\sup_{
\x\in \R^{nd }
}\bar \E^{\textcolor{black}{\bar \P_{\textcolor{black}{0},\x}}}\left[\int_{\textcolor{black}{0}}^T f(\textcolor{black}{s},\bar \X_s)ds\right] <1,$$
one also has:
$$\sup_{
\x\in \R^{nd}}\bar \E^{\textcolor{black}{\bar \P_{\textcolor{black}{0},\x}}}\left[\exp\left(\int_{\textcolor{black}{0}}^T f(\textcolor{black}{s},\bar \X_s)ds\right)\right]<\frac{1}{1-\alpha_{\textcolor{black}{T}}}.$$ 
\end{lemma}
As a corollary to the previous Lemma we obtain the following proposition \textcolor{black}{which will allow to apply the Girsanov Theorem to derive the existence of a solution to the martingale problem under the assumption \A{D}-(c) through a change of probability}. 
\begin{proposition}[Exponential integrability]\label{PROP_NOVIKOV}
Let $(\bar \X_t)_{t\ge 0}$ be the (unique-weak) solution to \eqref{DYN_NO_DRIFT}. Let $F_1 \in \textcolor{black}{L^{q}([0,T],L^{p}( \R^{nd}))}$ with $\frac{n^2d}{p}+\frac 2{q}<1, p\ge 2,q>2$. Then, \textcolor{black}{for any $\textcolor{black}{\lambda}>0 $}, there exists $K_{F_1,\textcolor{black}{\textcolor{black}{\lambda}}}:=K_{F_1,\textcolor{black}{\textcolor{black}{\lambda}}}(\A{A},T) $ s.t.
\begin{equation}
\label{PREAL_NOVIKOV}
\sup_{\x\in \R^{nd}}\bar \E^{\textcolor{black}{\bar \P_{0,\x}}}\left[\exp\left(\textcolor{black}{\textcolor{black}{\lambda}}\int_0^T |(\sigma^{-1} F_1)(\textcolor{black}{s},\bar \X_s)|^2 ds\right)\right]\le K_{F_1,\textcolor{black}{\textcolor{black}{\lambda}}}.
\end{equation} 
The constant $K_{F_1,\textcolor{black}{}\textcolor{black}{\lambda}}$ depends continuously on $\|F_1\|_{L^{q}([0,T],L^{p}( \R^{nd}))} $ \textcolor{black}{and $\textcolor{black}{\lambda} $}.
\end{proposition}
\textcolor{black}{We point out that in the above Proposition, the case $p=2$ can only be considered in the scalar non-degenerate case $d=n=1$}.
\begin{proof}
Observe that for $a>1$ s.t. $a(\frac{n^2d}{p}+\frac{2}{q})<1$, setting $p'=\frac{p}{2a},\ q'=\frac{q}{2a} $, so that indeed $\frac{n^2d}{p'}+\frac 2{q'}<2 $, one has:
\begin{align*}
\|\ |\sigma^{-1}&F_1|^{2a}\ \|_{L^{q'}([0,T],L^{ p'}(\R^{nd}))}=\Big(\int_{0}^T dt \Big\{\int_{\R^{nd}}d\x\textcolor{black}{(}|(\sigma^{-1}F_1)(t,\x)|^{2a})^{p'}\Big\}^{\frac{q'}{p'}}\Big)^{\frac{1}{q'}}\\
&\le  \textcolor{black}{\kappa^a}\Big(\int_{0}^T dt \Big\{\int_{\R^{nd}}d\x|F_1(t,\x)|^{p}\Big\}^{\frac{q}{p}}\Big)^{\frac{2a}{q}}\le \textcolor{black}{\kappa^a} \| F_1 \|_{L^{q}([0,T],L^{ p}(\R^{nd}))}^{2a}.
\end{align*}
From equation \eqref{EST_DENS_SOL_PB_MART_NO_DRIFT} we thus derive that for all $\x\in \R^{nd} $:
\begin{align}\label{CTR_PREAL_GIR_PROP}
\Big|\bar \E^{\bar \P_{0,\x}}&[\int_0^T |\sigma^{-1}F_1(s,\bar \X_s)|^{2a}ds]\Big|\\
&\le \textcolor{black}{C_{\ref{EST_DENS_SOL_PB_MART_NO_DRIFT}}} \|\ |\sigma^{-1}F_1|^{2a}\ \|_{L^{q'}([0,T],L^{ p'}(\R^{nd}))}\le C_{\ref{EST_DENS_SOL_PB_MART_NO_DRIFT}} \textcolor{black}{\kappa^a}\| F_1 \|_{L^{q}([0,T],L^{ p}(\R^{nd}))}^{2a} .\notag
\end{align}
For $a$ as above, write now for $\varepsilon\in (0,1) $ and from the Young inequality:
\begin{align*}
\bar \E^{\bar \P_{0,\x}}&\left[\exp\left(\textcolor{black}{\textcolor{black}{\lambda}}\int_0^T |(\sigma^{-1} F_1)(\textcolor{black}{s},\bar \X_s)|^2 ds\right)\right]\\
&\le \bar \E^{\bar \P_{0,\x}}\left[\exp\left(\textcolor{black}{\textcolor{black}{\lambda}}\int_0^T |(\sigma^{-1} F_1)(s,\bar \X_s)|^{2a}\frac{\varepsilon}a ds\right)\right] C(T,a,\varepsilon,\textcolor{black}{\textcolor{black}{\lambda}}). 
\end{align*}
The statement now directly follows from the above equation, \eqref{CTR_PREAL_GIR_PROP} and Lemma \ref{LEMME_KHAS} taking $\varepsilon:= (2\lambda \textcolor{black}{C_{\ref{EST_DENS_SOL_PB_MART_NO_DRIFT}}}\kappa^a \|\ F_1 \|_{L^{q}([0,T],L^{ p}(\R^{nd}))}^{2a} )^{-1}$.
\end{proof}

Proposition \ref{PROP_NOVIKOV} ensures that the Novikov condition is fulfilled in order to prove existence for the martingale problem associated with $(L_t)_{t\ge 0}$ for  $F_1 $ satisfying \A{D}-(c) starting from $\bar \P_{0,\x} $ and the associated dynamics \eqref{DYN_NO_DRIFT} of the canonical process. Set, 
$$\tilde W_t:= W_t-\int_0^t \textcolor{black}{(\sigma^{-1}F_1)}(s,\bar \X_s)ds,\ t\in [0,T].$$
From Proposition \ref{PROP_NOVIKOV}, we  derive that:
$$M_t:=\exp\left( \int_0^t \textcolor{black}{(\sigma^{-1}F_1)}(s,\bar \X_s)dW_s-\frac{1}{2}\int_{0}^t|\textcolor{black}{(\sigma^{-1}F_1)}(s,\bar \X_s)|^2 ds\right),\ t\in [0,T]$$
is \textcolor{black}{a} $\bar \P $-$ \F_t $ martingale (\textcolor{black}{here}  $(\F_s)_{s\in [0,T]} $ stands for the natural filtration associated with the canonical process $(\bar \X_s)_{s\in [0,T]} $ under $\bar \P $). It follows from the Girsanov theorem that $(\tilde W_t)_{t\in [0,T]} $ is a Wiener process on $ (\Omega,\F_T,(\F_s)_{s\in [0,T]},\P)$ where $\textcolor{black}{\left.[d\P/d\bar \P]\right|_{\F_T}}$ 
$:=M_T $. The dynamics of $(\bar \X_s)_{s\in [0,T]} $, writes under $\P $:
$$d\bar \X_t=\Big(BF_1(t,\bar \X_t)+\bar \gF(t,\bar \X_t)\Big)dt+B\sigma(t,\bar \X_t)d\tilde W_t=\gF(t,\bar \X_t)dt+B\sigma(t,\bar \X_t)d\tilde W_t,$$
that is $(\bar \X_t)_{t\in [0,T]} $ solves \eqref{SYST} under $\P $.

\subsection{Krylov bounds under \A{D}-(b) and (c)} 
\textcolor{black}{We aim here at proving estimate \eqref{KRYLOV_1} under assumptions \A{D}-(b) and (c). We mainly focus on case (c) which is the more involved. Case (b) can indeed be derived from the methodology developed under \A{D}-(a) and the controls obtained in Section \ref{EX_UN_D_B} or alternatively following the approach presented below for \A{D}-(c) which remains valid for a bounded drift on the non-degenerate component.}

\textcolor{black}{
Thanks to Proposition \ref{PROP_NOVIKOV}, the proof of the Krylov bound \eqref{KRYLOV_1} is similar to the proof of Lemma 3.3 in \cite{kryl:rock:05}. We provide a complete proof in our context for the sake of completeness. }

Write now for $f\in L^{q'}([0,T],L^{p'}(\R^{nd})), \frac{n^2d}{p'}+\frac{2}{q'}<2, p'>1,q'>1 $, $(t,\x)\in [0,T]\times \R^{nd} $ with the notations of the previous paragraph:
\begin{eqnarray*}
|\E^{\P_{t,\x}}[\int_t^T f(s,\X_s)ds]|&=&\Big| \bar \E^{\bar \P_{t,\x}}\Big[\frac{M_T}{M_t}\int_t^T f(s,\X_s)ds\Big]\Big|\\
&\le& \bar \E^{\bar \P_{t,\x}}\Big[\Big( \frac{M_T}{M_t}\Big)^\alpha\Big ]^{\frac 1\alpha}T^{\frac 1\alpha} \Big(\bar \E^{\bar \P_{t,\x}}\Big[\int_t^T |f(s,\bar \X_s)|^\beta ds\Big]\Big)^{\frac 1\beta},
\end{eqnarray*}
where $\alpha^{-1}+\beta^{-1}=1, \alpha>1,\beta>1 $. We have from Proposition \ref{PROP_NOVIKOV} that the exponential martingale has moments of all orders so that \textcolor{black}{there exists $C_T>0$, with $C_T \to 0$ when $T \to 0$ such that}:
\begin{align*}
|\E^{\P_{t,\x}}[\int_t^T f(s,\X_s)ds]|
\le &\textcolor{black}{C_T} \Big(\bar \E^{\bar \P_{t,\x}}\Big[\int_t^T |f(s,\bar \X_s)|^\beta ds\Big]\Big)^{\frac 1\beta}\\
 &\le \textcolor{black}{C_T}\| |f|^\beta \|_{L^{q''}([0,T],L^{p''}(\R^{nd}))}^{\frac 1\beta}\\ 
&\le \textcolor{black}{C_T}\|f\|_{L^{q'' \beta}([0,T],L^{p''\beta}(\R^{nd}))},
\end{align*}
if $p'',q'' \in [1,+\infty] $ are s.t. $n^2d/p''+2/q''<\textcolor{black}{2} $. Taking $\beta>1 $ sufficiently close to $1$ s.t. the previous condition holds for $ p''=p'/\beta, q''=q'/\beta$ eventually yields the result.

\section{Counter example}\label{SEC_COUNTER_EX}
This section is devoted to the almost sharpness of \textcolor{black}{the thresholds appearing in} Theorem \ref{THM1}. \textcolor{black}{This is the purpose of Theorem \ref{THMAS} which we now prove. We only focus here on the statement concerning the almost sharpness of the H\"older exponents $\beta_i^j$, $i\leq j$ in $\leftB 2,n\rightB^2$ in \eqref{COND_HOLDER}. \textcolor{black}{Indeed}, as emphasized in the introduction, \textcolor{black}{the sharpness of $L^q-L^p$ integrability conditions on the first component $F_1 $ of the drift follows from \textcolor{black}{Example 69 of} \cite{beck_stochastic_2014}} (see also Proposition 3.3 of \cite{gradinaru_existence_2013})}. We eventually recall that, the almost sharpness of the coefficients $\beta_i^i$, $i$ in $\leftB 2,n\rightB$, has been already proved in \cite{chau:16}.  \\


Let us first introduce the main idea of our counter example. As we already discussed, the \textcolor{black}{Peano} system \eqref{eq:peano} is ill posed as soon as $\alpha$ is in $(0,1)$ and $Y$ starts from 0 and well posed (in a strong sense) as soon as it is  suitably perturbed. In \cite{delarue_transition_2014}, the \textcolor{black}{a}uthors show that, in order to regularize, there must exist a transition time strictly less than one such that, \textcolor{black}{before} this time, the noise dominates \textcolor{black}{in the dynamics of} the system and \textcolor{black}{therefore} allows the solution to leave the singularity. This competition can be written explicitly and give\textcolor{black}{s} the following (heuristic) rule: the fluctuations \textcolor{black}{of order} $\gamma$ of the noise added in the system has to be strictly lower than $1/(1-\alpha)$. \textcolor{black}{We formalize these facts with our Proposition \ref{FAILURE_OF_WEAK_UNIQUENESS} below}.

\textcolor{black}{This proposition will be the key tool} to handle each H\"older threshold which depend\textcolor{black}{s} on the component and the variable. Hence, in order to \textcolor{black}{exhibit} the (almost) optimal threshold for the drift of the $i^{\rm th}$ component with respect to the $j^{\rm {th}}$ variable we need to build an \emph{ad hoc} Peano-like example. Focusing on the $i^{\rm th}$ component and the $j^{\rm {th}}$ variable, \textcolor{black}{for $i\in \leftB 2,n\rightB$ and $j\in \leftB i,n\rightB$}, we consider:
\begin{equation*}
\left\lbrace\begin{array}{l}
\displaystyle\dot{x}_t^1 = \dot{x}_t^2 = \ldots = \dot{x}_t^{i-1} = 0\\
\\
\displaystyle\dot{x}_t^i = {\rm sign}(x_t^j)|x_t^j|^{\beta_i^j} \\
\\
\displaystyle \dot{x}_t^{i+1} = x_t^i, \ldots   ,\dot{x}_t^{j} = x_t^{j-1} \\
\\
\displaystyle \dot{x}_t^{j+1} = \textcolor{black}{\dot x_t^{j+2}}= \ldots   =\dot{x}_t^{n} = 0\\
\end{array}\right.
\quad t \geq 0,
\end{equation*}
and $x_0^l=0$ for all $l$ in $\leftB 1,j\rightB$. \textcolor{black}{Each entry $(x_t^k)_{t\ge 0, k\in \leftB 1,n\rightB} $ of the above dynamics is scalar}. It is well seen that the \textcolor{black}{global well posedness} of this system relies on \textcolor{black}{the well posedness of the $i^{\rm th}$ equation} whose \textcolor{black}{extremal} solutions write $\pm c^{ij}_{\beta_i^j} t^{((j-i)\beta_i^j+1)/(\textcolor{black}{1-\beta_i^j})}$, \textcolor{black}{for some positive $c^{ij}_{\beta_i^j}$}. \textcolor{black}{For the considered initial point $\x_0=(x_0^1,\cdots, x_0^n)=(0,\cdots,0) $, the above dynamics can be rewritten in short as:
\begin{equation*}
d\x_t=[{\mathbf A}\x_t+F_{i}^j(\x_t)]dt,\ {\mathbf A}=\left( \begin{array}{ccccc} 
0                             & \cdots & \cdots & \cdots &0\\
1& 0        & \cdots& 0 & \vdots\\
0                              & 1&0        &\vdots &\vdots\\
\vdots                              & \ddots &\ddots        &0 & \vdots\\
 0                               &   \cdots          & 0&    1               &0
\end{array}\right),
\end{equation*}
where
$F_{i}^j(\x) =e_i  {\rm sign}(x_j)|x_j|^{\beta_i^j}$ and  $e_i$ is the $i^{\rm{th}}$ vector of the canonical basis of $\R^n$.
In that case, the \textcolor{black}{corresponding} stochastic perturbation \textcolor{black}{has} the form
\begin{equation}\label{SYSTPEANOPER}
d\X_t=[{\mathbf A}\X_t+F_{i}^j(\X_t)]dt +B dW_t,
\end{equation}
}
 This in particular means that the perturbation of the $i^{{\rm th}}$ component is done by the $(i-1)^{{\rm th}}$ iterated integrals (in time) of Brownian motion. Thus, focusing on the $i^{{\rm th}}$  \textcolor{black}{level} of the chain and the $j^{{\rm th}}$ component this means that we are interested in the following type of SDE

\begin{align}
Z_t^{\textcolor{black}{i,j}} &= x + \mW_t^{\textcolor{black}{i}}\label{eq:peanomulti_ij}\\
&+\int_0^t {\rm sign}\Bigg(\int_0^{s}  \ldots \int_0^{\textcolor{black}{s_{i+2}}} Z_{s_{i+1}}^{\textcolor{black}{i,j}} ds_{i+1}\ldots d s_{j}\Bigg) \Bigg|\int_0^{s}  \ldots \int_0^{s_{\textcolor{black}{i+2}}} Z_{s_{i+1}}^{\textcolor{black}{i,j}} ds_{i+1}\ldots d s_{j}\Bigg|^{\textcolor{black}{\beta_i^j}} d s\notag
\end{align}
\textcolor{black}{where $\mW^{i}$ will be chosen as the $(i-1)^{{\rm th}}$ iterated integral in time of Brownian motion.}
\textcolor{black}{The non-uniqueness in law for equation \eqref{eq:peanomulti_ij} will then follow from the next proposition}.

\begin{proposition}[Failure of the well posedness for the regularized Peano system]
\label{FAILURE_OF_WEAK_UNIQUENESS}
Let $\mW$ be a random process with continuous path\textcolor{black}{s} satisfying, in law, an invariance by symmetry and a self-similarity property of order $\gamma>0 $. Namely: 
$$(\mW_t,\ t \geq 0) \textcolor{black}{\overset{\rm{(law)}}{=}} (-\mW_t,\ t \geq 0), \ \textcolor{black}{\forall \rho>0, (\mW_{\rho t})_{t\ge 0}\textcolor{black}{\overset{\rm{(law)}}{=}}  (\rho^\gamma\mW_t)_{t\ge 0}}.$$
Suppose moreover that $\E\textcolor{black}{[\sup_{s\in [0,1]}|\mW_s|]} < +\infty$ and that $\mW$ and $\alpha<1$ are such that there exists a weak solution to the following SDE: 
\begin{equation}\label{eq:peanomulti}
Z_t = x + \int_0^t {\rm sign}\Bigg(\int_0^{s}  \ldots \int_0^{\textcolor{black}{s_{2}}} Z_{s_{1}} ds_{1}\ldots d s_{k}\Bigg) \Bigg|\int_0^{s}  \ldots \int_0^{\textcolor{black}{s_{2}}} Z_{s_{1}} ds_{1}\ldots d s_{k}\Bigg|^{\alpha} d s+ \mW_t,
\end{equation}
for any initial condition $x \geq 0$ where $k\in \mathbb{N}$ is given and  that it satisfies the Kolmogorov criterion.

Then, if  $\alpha <(\gamma - 1)/(k+\gamma)$, uniqueness in law fails for \eqref{eq:peanomulti}. 
\end{proposition}

Turning now to our claim, it is clear that for $\beta_i^j>0$ \eqref{eq:peanomulti_ij}, admits for all initial condition $x
\textcolor{black}{\ge} 
0$ at least one solution which satisfies the Kolmogorov Criterion. Our statement concerning the non uniqueness in law for the solution of \eqref{eq:peanomulti_ij} then readily follows from Proposition \ref{FAILURE_OF_WEAK_UNIQUENESS}. Taking  $\alpha=\beta_i^j$, $\mW=\mW^i$, which corresponds to the $(i-1)^{{\rm th}}$ iterated integrals (in time) of Brownian motion and therefore induces to take $\gamma=i-\frac 12 $, and $k=j-i$, we deduce that weak uniqueness fails as soon as 
$$\beta_i^j <  \frac{2i-3}{2j-1}.$$
It now remains to prove Proposition \ref{FAILURE_OF_WEAK_UNIQUENESS}.

\begin{proof}[Proof of Proposition \ref{FAILURE_OF_WEAK_UNIQUENESS}]
Consider the \textcolor{black}{extremal} solutions of the deterministic solutions associated with \eqref{eq:peanomulti}, that correspond to the case $ \mW_t=0$ and write $\pm c_{\alpha,k} t^{(k\alpha+1)/(\textcolor{black}{1-\alpha})} $.
The crucial point consists in comparing the fluctuations of the noise in the dynamics of \eqref{eq:peanomulti} with the \textcolor{black}{extremal} solution of the associated Peano Example. The proof follows the lines of \cite{chau:16} but we decide to reproduce it here for the sake of completeness.

For a given  parameter $\beta\in (0,1)$, we define for any continuous path $Y$ from $\R^+$ to $\R$ the variable  $\tau(Y)$ as 
$$\tau(Y) = \inf\{ t \geq 0\ : \ Y_t \leq (1-\beta)c_{\alpha,k} t^{(k\alpha + 1)/(1-\alpha)}\}.$$
\textcolor{black}{The stopping time $\tau(Y)$ then corresponds to the first passage of $Y$ below a threshold related to the (positive) extremal solution of the deterministic Peano system. Then, the key point to the proof of Proposition \ref{FAILURE_OF_WEAK_UNIQUENESS} is the following Lemma.}
\begin{lemma}\label{lemme:ce}
Let $Z$ be a weak solution of \eqref{eq:peanomulti} starting from some $x>0$ and suppose that $\alpha <(\gamma - 1)/(k+\gamma)$. Then, there exists a positive $\rho$, depending on $\alpha$, $\beta$, $\gamma$ and $\E |\mW_1|$ only\textcolor{black}{,} such that 
\begin{eqnarray}\label{eq:ce}
\mathbb{P}_x(\tau(\textcolor{black}{Z}) \geq \rho) \geq 3/4.
\end{eqnarray}
\end{lemma}
\textcolor{black}{Roughly speaking, the result tells us that, when the noise in the system is not strong enough, the solution from above the (positive) extremal solution will remain above with great probability. If weak uniqueness holds, the symmetry property implies that any solution from below the (negative) extremal solution will remain below with great probability. Letting the starting point tend to 0 (i.e. to the singularity) this leads to a contradiction. Together with the above Lemma, this last fact will allow us to conclude our counter-example.}\\

\textcolor{black}{Let $(Z,\mW)$ be a weak solution of \eqref{eq:peanomulti} with the initial condition $x=0$. Then, $(-Z,-\mW)$ is also a weak solution of \eqref{eq:peanomulti} so that, if uniqueness in law holds, $Z$ and $-Z$ have the same law.}

Let $Z^n$ be a sequence of weak solutions of \eqref{eq:peanomulti} starting from $1/n$, $n$ being a positive integer and $(\mathbb{P}_{1/n})_{n \geq 0}$ its law. Thanks to  Kolmogorov's criterion, we can extract a converging subsequence $(\mathbb{P}_{1/n_k})_{k\geq 0}$ that converges to $\mathbb{P}_0$, the law of the weak solution $Z$ of \eqref{eq:peano} starting from 0. Since the bound in \eqref{eq:ce} does not depend on the initial condition we get that 
$$\mathbb{P}_0(\tau(\textcolor{black}{Z}) \geq \rho) \geq 3/4,$$
and, thanks to uniqueness in law
$$\mathbb{P}_0(\tau(-\textcolor{black}{Z}) \geq \rho) \geq 3/4,$$
which \textcolor{black}{is obviously impossible}.\\
\end{proof}

\begin{proof}[Proof of Lemma \ref{lemme:ce}] 

 Let $Z$ be a weak solution of \eqref{eq:peanomulti} starting from $x>0$. Since it has continuous path, we have almost surely that $\tau(Z)>0$. Then, note that \textcolor{black}{for $t\in [0,\tau(Z)]$} we have:
\begin{eqnarray*}
Z_t &=& x + \int_0^t {\rm sign}\Bigg(\textcolor{black}{\int_0^{s} \ldots \int_0^{s_2}} Z_{s_1} ds_1\ldots d s_{k}\Bigg) \Bigg|\textcolor{black}{\int_0^s \ldots \int_0^{s_2}} Z_{s_1} ds_1\ldots d s_{k} \Bigg|^{\alpha} \textcolor{black}{d s}+ \mW_t\\
&\geq & (1-\beta)^{\alpha} \textcolor{black}{\bar c_{\alpha,k}} t^{(k\alpha+1)/(\textcolor{black}{1-\alpha})} + \mW_t,
\end{eqnarray*}
\textcolor{black}{for some positive constant $ \bar c_{\alpha,k}$.
Indeed, from the very definition of $\tau(Z) $, one gets $Z_{s_1}\ge (1-\beta)c_{\alpha,k}s_1^{(k\alpha+1)/(1-\alpha)} $, which once integrated in time and taking the $\alpha^{{\rm th}} $ power yields a lower bound of the form $ (1-\beta)^\alpha c_{\alpha,k}^\alpha \mathfrak c_{\alpha,k}^\alpha t^{\alpha(k+(k\alpha+1)/(1-\alpha))}= (1-\beta)^\alpha c_{\alpha,k}^\alpha \mathfrak c_{\alpha,k}^\alpha t^{\alpha(k+1)/(1-\alpha)}$ for some additional constant $\mathfrak c_{\alpha,k}>0 $ related to the iterated time integration, recalling that $\alpha<1 $ and that $t$ is small for the last inequality. One sets eventually $ \bar c_{\alpha,k}:=c_{\alpha,k}^\alpha \mathfrak c_{\alpha,k}^\alpha$}.

Hence, choosing $\eta$ such that $(1-\eta) = [(1-\beta)^{\alpha}+(1-\beta)]/2$\textcolor{black}{, we observe that $\beta-\eta + 1-\eta = (1-\beta)^\alpha$} and we get that:

\begin{eqnarray*}
Z_t  &\geq & (1-\eta)\textcolor{black}{\bar c_{\alpha,k}} t^{(k\alpha+1)/(\textcolor{black}{1-\alpha})}+ (\beta-\eta) \textcolor{black}{\bar c_{\alpha,k}} t^{(k\alpha+1)/(\textcolor{black}{1-\alpha})} + \mW_t,
\end{eqnarray*}
for all $t$ in $[0,\tau(Z)]$.

Let now $\rho$ be a positive number \textcolor{black}{to be specified later on}. Set $\textcolor{black}{\tilde c_{\alpha,k}} = (\beta-\eta) \textcolor{black}{\bar c_{\alpha,k}}$ and define
$$A = \left\{\textcolor{black}{\tilde c_{\alpha,k}} t^{(k\alpha+1)/(\textcolor{black}{1-\alpha})}+ \mW_t > 0 \text{ for all } t \text{ in } (0,\rho]\right\}.$$
\textcolor{black}{The event $A$ allows us to compare the fluctuations of the noise and those of the (positive) extremal solution. More precisely, it is the set of realizations for which the amount of noise in the system is lacking.} 
Note that on $A$ we have
\begin{eqnarray*}
Z_t  \geq (1-\eta)\textcolor{black}{\bar c_{\alpha,k}} t^{(k\alpha+1)/(\textcolor{black}{1-\alpha})} \geq  (1-\beta)\textcolor{black}{\bar c_{\alpha,k}} t^{(k\alpha+1)/(\textcolor{black}{1-\alpha})}
\end{eqnarray*}
for all $t$ in $[0,\rho]$, \textcolor{black}{recalling that $\alpha,\beta\in (0,1) $ for the last inequality}. But this is compatible only with the event $\{\tau(Z) \geq \rho\}$ so that $A \subset \{\tau(Z) \geq \rho\}$. Hence
\begin{equation}
\mathbb{P}(\tau(Z) \geq \rho) \geq \mathbb{P}(A).
\end{equation}
\textcolor{black}{It remains to choose $\rho>0$ such that $\mathbb{P}( A)\ge 3/4$. Write:
\begin{align*}
\mathbb{P}( A)=&\mathbb P[\forall t\in (0,\rho],\  \tilde c_{\alpha,k}t^{\frac{k\alpha+1}{1-\alpha}}+\mathcal{W}_t>0]=\mathbb P[\forall t\in (0,1],\  \tilde c_{\alpha,k}(\rho t)^{\frac{k\alpha+1}{1-\alpha}}+\mathcal{W}_{\rho t}>0]\\
=&\mathbb P[\forall t\in (0,1],\  \tilde c_{\alpha,k}(\rho t)^{\frac{k\alpha+1}{1-\alpha}}+\rho^\gamma \mathcal{W}_{t}>0]\\
=&\mathbb P[\forall t\in (0,1],\  \tilde c_{\alpha,k}\rho^{\frac{k\alpha+1}{1-\alpha}-\gamma}+t^{-\frac{k\alpha+1}{1-\alpha}} \mathcal{W}_{t}>0],
\end{align*}
from the self-similarity assumption on $\mathcal W$. Since by assumption $\alpha<\frac{\gamma-1}{\gamma+k}\iff\frac{k\alpha+1}{1-\alpha}-\gamma<0 $, the statement will follow taking $\rho $ small enough as soon as we prove the process $\mathcal R_t:=t^{-\frac{k\alpha+1}{1-\alpha}} \mathcal{W}_{t},\ t\in (0,1] $, which is continuous on the open set $(0,1]$, can be extended by continuity in $0$ with $\mathcal R_0=0$. Observe that $\mathbb E[|\mathcal R_t|]=t^{\gamma-\frac{k\alpha+1}{1-\alpha}} \mathbb E[|\mathcal W_1|]\underset{t\rightarrow 0}{\longrightarrow} 0$. Setting $\delta:= \gamma-\frac{k\alpha+1}{1-\alpha}>0$ and introducing $ t_n:=n^{-1/\delta(1+\eta)}$, $\eta>0$, we get that for all $\varepsilon>0 $,
$$\mathbb P[|\mathcal R_{t_n}|\ge \varepsilon]\le \varepsilon^{-1}\mathbb E[|\mathcal R_{t_n}|]=\varepsilon^{-1}t_n^\delta \mathbb E[|\mathcal W_1|]=\varepsilon^{-1}n^{-(1+\eta)}\mathbb E[|\mathcal W_1|].$$
We thus get from the Borel-Cantelli lemma that $\mathcal R_{t_n}\underset{n,\ a.s.}{\longrightarrow} 0 $. Namely, we have almost sure convergence along the subsequence $t_n$ going to zero with $n$. It now remains to prove that the process $\mathcal R_{t} $ does not fluctuate much between two successive times $t_n $ and $t_{n+1}$. Write for $t\in [t_{n+1},t_n]$:
\begin{align}
|\mathcal R_t|:=|t^{-\frac{k\alpha+1}{1-\alpha}} \mathcal{W}_{t}|\le& t_{n+1}^{-\frac{k\alpha+1}{1-\alpha}}\Big(|\mathcal{W}_{t_{n+1}}|+ \sup_{s\in [t_{n+1},t_n]}|\mathcal{W}_{s}-\mathcal{W}_{t_{n+1}}|\Big)\notag\\
\le &  t_{n+1}^{-\frac{k\alpha+1}{1-\alpha}}\Big(2|\mathcal{W}_{t_{n+1}}|+ \sup_{s\in [0,t_n]}|\mathcal{W}_{s}|\Big)\label{CTR_REMAIN_BC}.
\end{align}
The first term of the above left hand side tends almost surely to zero with $n$. Observe as well that, from the scaling properties of $\mathcal{W}$, for any $\varepsilon>0 $:
\begin{align*}
 \mathbb P[ t_{n+1}^{-\frac{k\alpha+1}{1-\alpha}} \sup_{s\in [0,t_n]}|\mathcal{W}_{s}|\ge \varepsilon]&=\mathbb P[ t_{n+1}^{-\frac{k\alpha+1}{1-\alpha}}t_n^\gamma \sup_{s\in [0,1]}|\mathcal{W}_{s}|\ge \varepsilon]\\
 \le \varepsilon^{-1}t_n^{\delta} (\frac{t_n}{t_{n+1}})^{\frac{k\alpha+1}{1-\alpha}} \mathbb E[\sup_{s\in [0,1]}|\mathcal{W}_{s}|]
 &\le C\varepsilon^{-1}n^{-(1+\eta)}\mathbb E[\sup_{s\in [0,1]}|\mathcal{W}_{s}|],
\end{align*}
which again gives from the Borel-Cantelli lemma the a.s. convergence with $n$ of the second term in the r.h.s of \eqref{CTR_REMAIN_BC}. We eventually derive
that $\mathcal R_t \underset{t\rightarrow 0,\ a.s.}{\longrightarrow} 0$.
Again, the key point is that we normalize  the process $\mathcal W$ at a rate, $t^{\frac{k\alpha+1}{1-\alpha}} $, which is lower than its own characteristic time scale, $t^{\gamma} $. This is precisely what leaves some margin to establish continuity}.


\end{proof}

\appendix

\section{\textcolor{black}{Proof of the Technical Lemmas \ref{convergence_dirac}} and \ref{convergence_LPLQ}}\label{sec:ProofLemma}

\subsection{Proof of Lemma \ref{convergence_dirac} \textcolor{black}{(Dirac convergence of the frozen density)}}

For the proof of this Lemma 
we are somehow faced with the same type of difficulties as for Lemma \ref{lemme:bilipflow}. Namely, recalling the expression of  $\tilde p^{t+\varepsilon,\y}(t,t+\varepsilon,\x,\y) $ derived from \eqref{CORRESP}, we have a dependence of the covariance  matrix $\tilde \K_{t+\varepsilon,t}^{t+\varepsilon,\y}$ and \textcolor{black}{of} the linearized flow $ \tilde \btheta_{t+\varepsilon,t}^{t+\varepsilon,\y}(\x)$ in the integration variable $\y $. 

\textcolor{black}{Let $(\btheta_{u,t}(\x)\big)_{u\in [t,t+\varepsilon]} $ be the forward flow provided by Lemma \ref{lemme:measu}. To study the sensitivity of the covariance matrix w.r.t. the flows
we now introduce, for a given point $\x\in \R^{nd}$, the linear Gaussian diffusion $(\bar \X_u)_{u\in [t,t+\varepsilon]} $ with dynamics:}
\begin{equation}
\label{BARX}
d\bar \X_u=D\gF(u,\btheta_{u,t}(\x))\bar \X_udu+B\sigma(u,\btheta_{u,t}(\x))dW_u.
\end{equation}
The associated covariance matrix between $t$ and $t+\varepsilon $ writes:
\begin{equation}
\label{COV}
\bar \K_{t+\varepsilon,t}^{t,\x}=\int_t^{t+\varepsilon} \bar \gR^{t,\x}(t+\varepsilon,u)B\sigma\sigma^*(u,\btheta_{u,t}(\x))B^*\bar \gR^{t,\x}(t+\varepsilon,u)^*du,
\end{equation}
where $(\bar\gR^{t,\x}(v,u))_{t\le u,v\le t+\varepsilon} $ stands for the resolvent associated with $$(D\gF(u,\textcolor{black}{\btheta_{u,t}(\x)}))_{u\in [t,t+\varepsilon]}.$$ \textcolor{black}{We point out \textcolor{black}{that} $(\bar\gR^{t,\x}(v,u))_{t\le u,v\le t+\varepsilon}$ and $(\tilde \gR^{t+\varepsilon,\y}(v,u))_{t\le u,v\le t+\varepsilon}$ are \textit{similar resolvents}, in the sense that they actually only differ in the flow considered in the linear dynamics. The flow is forward for $\bar\gR^{t,\x}$ and backward for $\textcolor{black}{\tilde \gR^{t+\varepsilon,\y}} $}.

Observe now that, from \A{H}, $\bar \K_{t+\varepsilon,t}^{\textcolor{black}{t,\x}} $ satisfies the \textit{good scaling property} \eqref{GSP}.
Write now:
\begin{eqnarray}
&&\int_{\R^{nd}}  \tilde p^{t+\varepsilon,\y}(t,t+\varepsilon,\x,\y) f(\y)d\y\notag\\
&=&\Bigg[\int_{\R^{nd}} \Bigg\{\frac{\exp\left(-\frac 12 \left\langle (\tilde \K_{t+\varepsilon,t}^{t+\varepsilon,\y})^{-1}(\tilde \btheta_{t+\varepsilon,t}^{t+\varepsilon,\y}(\x)-\y),\tilde \btheta_{t+\varepsilon,t}^{t+\varepsilon,\y}(\x)-\y \right\rangle \right)}{(2\pi)^{\frac{nd}2}\det(\tilde \K_{t+\varepsilon,t}^{t+\varepsilon,\y})^{\frac 12}}\notag\\
&&-\frac{\exp\left(-\frac 12 \left\langle (\tilde \K_{t+\varepsilon,t}^{t+\varepsilon,\y})^{-1}( \btheta_{t+\varepsilon,t}(\x)-\y), \btheta_{t+\varepsilon,t}(\x)-\y \right\rangle \right)}{(2\pi)^{\frac{nd}2}\det(\tilde \K_{t+\varepsilon,t}^{t+\varepsilon,\y})^\frac12}\Bigg\} f(\y) d\y\Bigg]\notag\\
&&+\Bigg[\int_{\R^{nd}} \Bigg\{\frac{\exp\left(-\frac 12 \left\langle (\tilde \K_{t+\varepsilon,t}^{t+\varepsilon,\y})^{-1}( \btheta_{t+\varepsilon,t}(\x)-\y), \btheta_{t+\varepsilon,t}(\x)-\y \right\rangle \right)}{(2\pi)^{\frac{nd}2}\det(\tilde \K_{t+\varepsilon,t}^{t+\varepsilon,\y})^{\frac 12}}\notag \\
&&-\frac{\exp\left(-\frac 12 \left\langle (\bar \K_{t+\varepsilon,t}^{t,\x})^{-1}( \btheta_{t+\varepsilon,t}(\x)-\y), \btheta_{t+\varepsilon,t}(\x)-\y \right\rangle \right)}{(2\pi)^{\frac{nd}2}\det(\bar \K_{t+\varepsilon,t}^{t,\x})^{\frac12}}\Bigg\}f(\y) d\y\Bigg]\notag\\
&&+\Bigg[\int_{\textcolor{black}{\R^{nd}}}\frac{\exp\left(-\frac 12 \left\langle (\bar \K_{t+\varepsilon,t}^{t,\x})^{-1}( \btheta_{t+\varepsilon,t}(\x)-\y), \btheta_{t+\varepsilon,t}(\x)-\y \right\rangle \right)}{(2\pi)^{\frac{nd}2}\det(\bar \K_{t+\varepsilon,t}^{t,\x})^{\frac 12}}f(\y) d\y\Bigg]\notag\\
&=:&\sum_{i=1}^3\bXi_i^\varepsilon(t,\x).\label{LAST_LABEL}
\end{eqnarray}
It is directly seen \textcolor{black}{from the dominated convergence theorem} that $\bXi_3^\varepsilon(t,\x) \underset{\varepsilon \downarrow 0}{\rightarrow} f(\x) $. It remains to prove that $\bXi_1^\varepsilon(t,\x), \bXi_2^\varepsilon(t,\x) $ can be viewed as remainders as $\varepsilon \downarrow 0 $.

Let us write
\begin{eqnarray}
|\bXi_{1}^\varepsilon(t,\x)|&\le & \|f\|_\infty \int_{\R^{nd}}^{} \frac{d\y}{(2\pi)^{\frac{nd}2}\det(\tilde \K_{t+\varepsilon,t}^{t+\varepsilon,\y})^{\frac 12}} \int_0^1 d\lambda |(\varphi_{t,\x,\y}^\varepsilon)'(\lambda)|,
\label{THE_FORM_R1}
\end{eqnarray}
\textcolor{black}{where $\lambda \in [0,1]$:}
\begin{eqnarray}
\varphi_{t,\x,\y}^\varepsilon(\lambda)&=&\exp\left(-\frac{1}{2}\left\{\langle (\tilde \K_{t+\varepsilon,t}^{t+\varepsilon,\y})^{-1} (
\textcolor{black}{\btheta_{t+\varepsilon,t}}
(\x)-\y),
\textcolor{black}{\btheta_{t+\varepsilon,t}}
(\x)-\y  \rangle +\right. \right. \nonumber\\
&& \left. \left. \lambda \left[\langle (\tilde \K_{t+\varepsilon,t}^{t+\varepsilon,\y})^{-1} ( \textcolor{black}{\tilde \btheta_{t+\varepsilon,t}^{t+\varepsilon,\y}(\x)}-\y), \textcolor{black}{\tilde \btheta_{t+\varepsilon,t}^{t+\varepsilon,\y}(\x)}-\y  \rangle  \right.\right.\right.\notag\\
&&\left.\left.\left.-\langle (\tilde \K_{t+\varepsilon,t}^{t+\varepsilon,\y})^{-1} (
\textcolor{black}{\btheta_{t+\varepsilon,t}
(\x)}-\y),
\textcolor{black}{\btheta_{t+\varepsilon,t}
(\x)}-\y  \rangle \right]\right\}\right),\nonumber\\
|(\varphi_{t,\x,\y}^\varepsilon)'(\lambda)|&\le& \left | (\tilde \K_{t+\varepsilon,t}^{t+\varepsilon,\textcolor{black}{\y}})^{-\frac12}\left\{(\btheta_{t+\varepsilon,t}(\x) -\y)+(\tilde \btheta_{t+\varepsilon,t}^{t+\varepsilon,\y}(\x) -\y) \right\}\right|\nonumber\\
&& \times \left|(\tilde \K_{\textcolor{black}{t+\varepsilon,t}}^{t+\varepsilon,\y})^{-\frac12}\left\{(\btheta_{t+\varepsilon,t}(\x) -\y)-(\tilde \btheta_{t+\varepsilon,t}^{t+\varepsilon,\y}(\x) -\y) \right\} \right|  \varphi_{t,\x,\y}^\varepsilon(\lambda),  \label{DEF_PHI_LAMBDA_THETA1} 
\end{eqnarray} 
using the Cauchy-Schwarz inequality for the last assertion. A key quantity to control for the analysis is now the linearization error  $| (\tilde \K_{t+\varepsilon,t}^{t+\varepsilon,\x})^{-\frac 12}(\btheta_{t+\varepsilon,t}(\x) -\tilde \btheta_{t+\varepsilon,t}^{t+\varepsilon,\y}(\x)) |$.
From $\eqref{GSP}$ we readily have:
\begin{eqnarray}\label{scalingsurK}\\
\left| (\tilde \K_{t+\varepsilon,t}^{t+\varepsilon,\x})^{-\frac 12}(\btheta_{t+\varepsilon,t}(\x) -\tilde \btheta_{t+\varepsilon,t}^{t+\varepsilon,\y}(\x)) \right|\le C \varepsilon^{\frac 12}|\T_{\varepsilon}^{-1}(\btheta_{t+\varepsilon,t}(\x)-\tilde \btheta_{t+\varepsilon,t}^{t+\varepsilon,\y}(\x))|.\notag
\end{eqnarray}
To bound the above r.h.s. we first introduce for $\z\in \R^{nd} $, $u\in [t,t+\varepsilon] $,
\begin{align}
\gF^{t+\varepsilon,\y}(u,\z)
:=\bigl( F_1(u,\btheta_{u,t+\varepsilon}(\y)),F_2(u,z_1,(\btheta_{u,t+\varepsilon}(\y))^{2,n}),\notag\\
F_3(u,z_2,(\btheta_{u,t+\varepsilon}(\y))^{3,n}),\cdots,
 F_n(s,z_{n-1},(\btheta_{u,t+\varepsilon}(\y))_{n})\bigr).\label{DEF_GFTY}
\end{align}
 We then write: 
\begin{align} 
\varepsilon^{\frac 12}\T_{\varepsilon}^{-1}&(\btheta_{t+\varepsilon,t}(\x)-\tilde \btheta_{t+\varepsilon,t}^{t+\varepsilon,\y}(\x))\notag\\
&:= \varepsilon^{\frac 12}\T_{\varepsilon}^{-1}\left\{ \int_{t}^{t+\varepsilon} du \biggl[\biggl( \gF(u,\btheta_{u,t}(\x))-\gF^{t+\varepsilon,\y}(u,\btheta_{u,t}(\x)) \biggr)\right.\nonumber\\ 
&\left.+\biggl(D\gF(u,\btheta_{u,t+\varepsilon}(\y))(\btheta_{u,t}(\x)-\tilde \btheta_{u,t}^{t+\varepsilon,\y}(\x) )\biggr) \right. \nonumber \\ 
&\left. + \biggl(\int_{0}^{1} d \lambda \left(  D\gF^{t+\varepsilon,\y}(u,\btheta_{u,t+\varepsilon}(\y)+\lambda (\btheta_{u,t}(\x)-\btheta_{u,t+\varepsilon}(\y) ) ) \right. \right.\nonumber\\ 
&-\left. \left. D\gF^{t+\varepsilon,\y}(u,\btheta_{u,t+\varepsilon}(\y)) \right) (\btheta_{u,t}(\x)-\btheta_{u,t+\varepsilon}(\y) )\biggr) \biggr] \right\}\nonumber \\ 
&:=(\cI_{t+\varepsilon,t}^{1}+\cI_{t+\varepsilon,t}^{2}+\cI_{t+\varepsilon,t}^{3})(\x,\y),
\label{D_EPS} 
\end{align} 
where, according to the notations of \eqref{DEF_GFTY}, for $(u,\z)\in [t,t+\varepsilon]\times \R^{nd}, D\gF^{t+\varepsilon,\y}(u,\z) $ is the $(nd)\times (nd)$ matrix with only non zero $d\times d $ matrix entries 
$(D\gF^{t+\varepsilon,\y}(u,\z))_{j,j-1}$
$:=D_{x_{j-1}}F_j(u,z_{j-1},\btheta_{u,t+\varepsilon}(\y)^{j,n})$, $ j\in \leftB 2,n\rightB $. In particular $D\gF^{t+\varepsilon,\y}(u,\btheta_{u,t+\varepsilon}(\y))=D\gF(u,\btheta_{u,t+\varepsilon}(\y)) $. 

Observe now that, from \A{S}:
\begin{eqnarray}
|\cI_{t+\varepsilon,t}^{3}(\x,\y)|&\le& C \sum_{i=2}^n \int_t^{t+\varepsilon} du \varepsilon^{-(i-1/2)} |(\btheta_{u,t}(\x)-\btheta_{u,t+\varepsilon}(\y) )_{i-1}|^{1+\eta}\notag\\
&\le& C \int_t^{t+\varepsilon} du \varepsilon^{-1+\eta/2} \Big(\varepsilon^{1/2}|\T_\varepsilon^{-1}(\btheta_{u,t}(\x)-\btheta_{u,t+\varepsilon}(\y) )|\Big)^{1+\eta}.\notag
\end{eqnarray}
From Lemma \ref{lemme:bilipflow} \textcolor{black}{(almost equivalence of the flows)} we now derive:
\begin{eqnarray}
|\cI_{t+\varepsilon,t}^{3}(\x,\y)|&\le&C \int_t^{t+\varepsilon} du \varepsilon^{-1+\eta/2} \Big(\varepsilon^{1/2}|\T_\varepsilon^{-1}(\btheta_{t+\varepsilon,t}(\x)-\y )|+1\Big)^{1+\eta}\notag\\
&\le& C\varepsilon^{\eta/2} \Big(\varepsilon^{1/2}|\T_\varepsilon^{-1}(\btheta_{t+\varepsilon,t}(\x)-\y )|+1\Big)^{1+\eta}.
\label{CTR_I3_EPS}
\end{eqnarray}
Let us now deal with $\cI_{t+\varepsilon,t}^{1}(\x,\y) $. From the previous definition of $\gF^{t+\varepsilon,\y} $ in \eqref{DEF_GFTY}, the key idea is to use the sub-linearity of $\gF $ and the appropriate H\"older exponents.
Namely, using the Young inequality we derive:
\begin{align*}
|\cI_{t+\varepsilon,t}^{1}(\x,\y)|&\le C\sum_{i=1}^n \varepsilon^{-i+1/2}\sum_{j=i}^n \int_t^{t+\varepsilon} du  |(\btheta_{u,t}(\x)-\btheta_{u,t+\varepsilon}(\y) )_{j}|^{\beta_i^j} \\
&\le  C\Bigg(\varepsilon^{-1/2}\int_t^{t+\varepsilon}du \Big(|\textcolor{black}{(\btheta_{u,t}(\x)-\btheta_{u,t+\varepsilon}(\y) )}|+1\Big) \\
&+\sum_{i=2}^n \varepsilon^{-i +1/2} \sum_{j=i}^n \int_t^{t+\varepsilon} du  \bigg\{ \bigg(\frac{|((\btheta_{u,t}(\x)-\btheta_{u,t+\varepsilon}(\y) )_{j})|}{\varepsilon^{\gamma_{i}^j}}\bigg)+\varepsilon^{\gamma_{i}^j \frac{\beta_i^j}{1-\beta_i^j}} \bigg\}\Bigg),
\end{align*}
for some parameters $\gamma_i^j >0$ to be specified below.
Hence,
\begin{align*}
|&\cI_{t+\varepsilon,t}^{1}(\x,\y)|\le  C\Bigg(\int_t^{t+\varepsilon} du \varepsilon^{1/2}|\T_\varepsilon^{-1}(\btheta_{u,t}(\x)-\btheta_{u,t+\varepsilon}(\y) )|+\varepsilon^{1/2} \\
&+\sum_{i=2}^n  \sum_{j=i}^n \int_t^{t+\varepsilon} du  \Big\{\varepsilon^{-i+j -\gamma_i^j}  \Big(\frac{|((\btheta_{u,t}(\x)-\btheta_{u,t+\varepsilon}(\y) )_{j})|}{\varepsilon^{j-1/2}}\Big)+\varepsilon^{-i+1/2+\gamma_{i}^j \frac{\beta_i^j}{1-\beta_i^j}} \Big\}\Bigg)\\
\le &C\Bigg(\int_t^{t+\varepsilon} du \varepsilon^{1/2}|\T_\varepsilon^{-1}(\btheta_{u,t}(\x)-\btheta_{u,t+\varepsilon}(\y) )|+\varepsilon^{1/2}\\
& +\sum_{i=2}^n  \sum_{j=i}^n \int_t^{t+\varepsilon} du  \Big\{\varepsilon^{-i+j -\gamma_i^j}  \varepsilon^{1/2}|\T_\varepsilon^{-1}(\btheta_{u,t}(\x)-\btheta_{u,t+\varepsilon}(\y) )|+\varepsilon^{-i+1/2+\gamma_{i}^j \frac{\beta_i^j}{1-\beta_i^j}} \Big\}\Bigg).
\end{align*}
We now use Lemma \ref{lemme:bilipflow} to derive 
$$ \varepsilon^{1/2}|\T_\varepsilon^{-1}(\btheta_{u,t}(\x)-\btheta_{u,t+\varepsilon}(\y) )|\le C(\varepsilon^{1/2}|\T_\varepsilon^{-1}(\btheta_{t+\varepsilon,t}(\x)-\y )|+1).$$ 
\begin{center}We emphasize here that in our current framework we should \textit{a priori} write\end{center}
$\btheta_{t+\varepsilon,u}(\btheta_{u,t}(\x)) $ in the above equation since we do not have a priori the flow property. Anyhow, since Lemma \ref{lemme:bilipflow} is valid for any flow starting from $\btheta_{u,t}(\x) $ at time $u$ associated with the ODE (see \textcolor{black}{equation \eqref{EQ_EQUIV_FLOW_DIFF}}
) we can proceed along the previous one, i.e. $(\btheta_{v,t}(\x))_{v\in [u,t+\varepsilon] } $.
This yields:
\begin{align}
&|\cI_{t+\varepsilon,t}^{1}(\x,\y)|\label{PREAL_CTR_I1_EPS}\\
\le & C\Big[\varepsilon^{1/2}+ (\varepsilon^{1/2}|\T_\varepsilon^{-1}(\btheta_{t+\varepsilon,t}(\x)-\y)|+1)\varepsilon\Big(1+\sum_{i=2}^n  \sum_{j=i}^n  \Big\{\varepsilon^{-i+j -\gamma_i^j}+\varepsilon^{-i+1/2+\gamma_{i}^j \frac{\beta_i^j}{1-\beta_i^j}} \Big\} \Big) \Big].\notag 
\end{align}
Choose now for $i\in \leftB 2,n\rightB $ and $j\in \leftB i,n\rightB $, 
$$-i+j -\gamma_i^j =-i+1/2+\gamma_{i}^j \frac{\beta_i^j}{1-\beta_i^j}\iff \gamma_i^j=(j-\frac 12)(1-\beta_i^j),$$ 
to balance the two previous contributions associated with the indexes $i,j$. To obtain a global smoothing effect w.r.t $\varepsilon $ in \eqref{PREAL_CTR_I1_EPS} we need to impose:
$$-i+j-\gamma_i^j>-1 \iff \beta_i^j>\frac{2i-3}{2j-1}.$$
Hence, under \eqref{COND_HOLDER}, we have that there exists $\zeta:=\zeta\big(\A{A}, (\beta_i^j)_{i\in\leftB 1,n\rightB, j\in \leftB i,n\rightB}\big) \in (0,1)$ s.t.:
\begin{eqnarray}
\label{CTR_I1_EPS}
|\cI_{t+\varepsilon,t}^{1}(\x,\y)| \le C\varepsilon^\zeta\big(1+\varepsilon^{1/2}|\T_\varepsilon^{-1}(\btheta_{t+\varepsilon,t}(\x)-\y)|\big).
\end{eqnarray}

We now get from \textcolor{black}{\eqref{scalingsurK}}, \eqref{D_EPS}, \eqref{CTR_I3_EPS}, \eqref{CTR_I1_EPS} and the Gronwall lemma that:
\begin{align}
&\left| (\tilde \K_{t+\varepsilon,t}^{t+\varepsilon,\x})^{-\frac 12}(\btheta_{t+\varepsilon,t}(\x) -\tilde \btheta_{t+\varepsilon,t}^{t+\varepsilon,\y}(\x)) \right|\notag\\
\le& C \varepsilon^{\frac 12}|\T_{\varepsilon}^{-1}(\btheta_{t+\varepsilon,t}(\x)-\tilde \btheta_{t+\varepsilon,t}^{t+\varepsilon,\y}(\x))|\notag\\
\le& C \varepsilon^{\eta/2 \wedge \zeta}\big(1+\varepsilon^{1/2}|\T_{\varepsilon}^{-1}(\btheta_{t+\varepsilon,t}(\x)-\y)| +(\varepsilon^{1/2}|\T_\varepsilon^{-1}(\btheta_{t+\varepsilon,t}(\x)-\y)|)^{1+\eta}\big).\label{LE_TERME_QUI_REGULARISE_THETA_1}
\end{align}
\textcolor{black}{Hence, recalling from \eqref{CTR_LINEA_RETRO} and Lemma \ref{lemme:bilipflow} that}
\begin{align*}
 \langle (\tilde \K_{t+\varepsilon,t}^{t+\varepsilon,\y})^{-1} &( 
\textcolor{black}{\tilde\btheta_{t+\varepsilon,t}^{t+\varepsilon,\y}
(\x)}-\y), 
\textcolor{black}{\tilde \btheta_{t+\varepsilon,t}^{t+\varepsilon,\y}
(\x)}-\y  \rangle\ge C^{-1} \varepsilon |\T_\varepsilon^{-1}(\x-\btheta_{t,t+\varepsilon}(\y))|^2\\
\ge& C^{-1} (\varepsilon |\T_\varepsilon^{-1}(\btheta_{t+\varepsilon,t}(\x)-\y)|^2-1) ,
\end{align*}
we get from the definition in \eqref{DEF_PHI_LAMBDA_THETA1} that for all $\lambda\in [0,1]$:
\begin{eqnarray*}
\varphi_{t,\x,\y}^\varepsilon(\lambda) \le C\exp(\textcolor{black}{-}C^{-1} \varepsilon |\T_\varepsilon^{-1}(\btheta_{t+\varepsilon,t}(\x)-\y)|^2)   .
\end{eqnarray*}
We finally  obtain from \eqref{DEF_PHI_LAMBDA_THETA1} and \eqref{LE_TERME_QUI_REGULARISE_THETA_1} that there exists $C_2:=C_2(T,\A{A})\ge 1$ s.t.:
\begin{eqnarray}
\label{CTR_DER_PHI}
|(\varphi_{t,\x,\y}^\varepsilon)'(\lambda)| &\le& 
C_2 \varepsilon^{\eta/2\wedge \zeta} \exp\left(-C_2^{-1}\varepsilon |\T_\varepsilon^{-1}(\btheta_{t+\varepsilon,t}(\x)-\y)|^2 \right).
\end{eqnarray}
Plugging Equation \eqref{CTR_DER_PHI} into \eqref{THE_FORM_R1}, we derive that, since $\tilde \K_{t+\varepsilon,t}^{t+\varepsilon,\y} $ satisfies \eqref{GSP},  $|\bXi_{1}^\varepsilon(t,\x)|\underset{\varepsilon \downarrow 0}{\rightarrow} 0 $.

Let us \textcolor{black}{now consider the term} $\bXi_2^\varepsilon(t,\x) $ \textcolor{black}{in \eqref{LAST_LABEL}}. Write first $\bXi_2^\varepsilon(t,\x):=(\bXi_{21}^\varepsilon+\bXi_{22}^\varepsilon)(t,\x) $ where:
\begin{align}
\bXi_{21}^\varepsilon(t,\x)&:=  \int_{\R^{nd}}^{}\frac{d\y}{(2\pi)^{\frac{nd}2}} \frac{f(\y)}{\det(\bar \K_{t+\varepsilon,t}^{t,\x})^{\frac 12}}\int_0^1 d\lambda (\psi_{t,\x,\y}^\varepsilon)'(\lambda),
\notag\\
\forall \lambda \in [0,1], \ \psi_{t,\x,\y}^\varepsilon(\lambda)&:=\exp\left(-\frac 12 \biggl\{ \langle (\bar \K_{t+\varepsilon,t}^{t,\x})^{-1} (\btheta_{t+\varepsilon,t}(\x)-\y),  \btheta_{t+\varepsilon,t}(\x)-\y\rangle    \right.\nonumber\\
&+\left.\lambda \biggl[ \langle (\tilde \K_{t+\varepsilon,t}^{t+\varepsilon,\y})^{-1} (\btheta_{t+\varepsilon,t}(\x)-\y), \btheta_{t+\varepsilon,t}(\x)-\y\rangle\right.\notag\\
& \left. -\langle (\bar \K_{t+\varepsilon,t}^{t,\x})^{-1} (\btheta_{t+\varepsilon,t}(\x)-\y),  \btheta_{t+\varepsilon,t}(\x)-\y\rangle\biggl] \biggr\}   \right),\nonumber \\
\bXi_{22}^\varepsilon(t,\x)&:=  \int_{\R^{nd}}^{}\frac{d\y}{(2\pi)^{\frac{nd}2}}f(\y)\Bigg[ \frac{1}{\det(\tilde \K_{t+\varepsilon,t}^{t+\varepsilon,\y})^{\frac 12}}-\frac{1}{\det(\bar \K_{t+\varepsilon,t}^{t,\x})^{\frac 12}}\Bigg] (\psi_{t,\x,\y}^\varepsilon)(1).\label{DEF_BXI_22}
\end{align}
Observe that for all $\lambda \in [0,1] $,
\begin{eqnarray*}
|(\psi_{t,\x,\y}^\varepsilon)'(\lambda)|&\le&  \left |  \langle   ( (\tilde \K_{t+\varepsilon,t}^{t+\varepsilon,\y})^{-1}-(\bar \K_{t+\varepsilon,t}^{t,\x})^{-1})( \btheta_{t+\varepsilon,t}(\x)-\y), \btheta_{t+\varepsilon,t}(\x)-\y \rangle \right | 
\psi_{t,\x,\y}^\varepsilon(\lambda)\nonumber.
\end{eqnarray*}
Equation \eqref{GSP}, which holds for $\bar \K_{t+\varepsilon,t}^{t,\x} $ as well,  yields:
\begin{eqnarray}
|(\psi_{\textcolor{black}{t},\x,\y}^\varepsilon)'(\lambda)|&\le& C \left |  \langle   ( (\tilde \K_{t+\varepsilon,t}^{t+\varepsilon,\y})^{-1}-(\bar \K_{t+\varepsilon,t}^{t,\x})^{-1})( \btheta_{t+\varepsilon,t}(\x)-\y), \btheta_{t+\varepsilon,t}(\x)-\y \rangle \right |\notag\\
&&\times \exp(-C\varepsilon |\T_{\varepsilon}^{-1}(\btheta_{t+\varepsilon,t}(\x)-\y)|^2)\nonumber \\
&\textcolor{black}{=:}&C |Q_\varepsilon| \exp(-C\varepsilon |\T_{\varepsilon}^{-1}( \btheta_{t+\varepsilon,t}(\x)-\y)|^2),\label{CTR_PSI_DER}
\end{eqnarray}
for $C:=C(\A{A},T)$. 

Now,  the covariance matrices explicitly write
\begin{eqnarray*}
\tilde \K_{t+\varepsilon,t}^{t+\varepsilon,\y}&=&\int_{t}^{t+\varepsilon} du \tilde \gR^{t+\varepsilon,\y}(t+\varepsilon,u)B a(u,\btheta_{u,t+\varepsilon}(\y))B^* \tilde \gR^{t+\varepsilon,\y}(t+\varepsilon,u)^* ,\\
 \bar \K_{t+\varepsilon,t}^{t,\x}&=&\int_t^{t+\varepsilon} du \bar \gR^{t,\x}(t+\varepsilon,u)B a(u,\textcolor{black}{\btheta_{u,t}(\x)})B^* \bar \gR^{t,\x}(t+\varepsilon,u)^* ,
 \end{eqnarray*}
 where $\tilde \gR^{t+\varepsilon,\y},\ \bar \gR^{t,\x} $ respectively denote the resolvents associated with the linear parts of equations \eqref{LIN_AROUND_BK_FLOW} and \eqref{BARX}.
Thus, setting:
\begin{equation}
\label{DEF_K_MACRO}
\tilde \K_{t+\varepsilon,t}^{t+\varepsilon,\y}=\varepsilon^{-1}\T_\varepsilon \widehat {\tilde \K}_1^{t+\varepsilon,t,t+\varepsilon,\y} \T_\varepsilon ,\ \bar \K_{t+\varepsilon,t}^{t,\x}=\varepsilon^{-1} \T_\varepsilon \widehat{\bar  \K}_1^{t+\varepsilon,t,t,\x} \T_\varepsilon , \end{equation}
we write:
\begin{eqnarray*}
&&\left |\langle (\tilde \K_{t+\varepsilon,t}^{t+\varepsilon,\y}-\bar \K_{t+\varepsilon}^{t,\x}) (\btheta_{t+\varepsilon,t}(\x)-\y), \btheta_{t+\varepsilon,t}(\x)-\y \rangle\right|\\
&&=\left| \langle (\widehat{\tilde \K}_1^{t+\varepsilon,t,t+\varepsilon,\y}-\widehat {\bar \K}_1^{t+\varepsilon,t,t,\x})(\varepsilon^{-1/2}\T_\varepsilon (\btheta_{t+\varepsilon,t}(\x)-\y)),\varepsilon^{-1/2}\T_\varepsilon (\btheta_{t+\varepsilon,t}(\x)-\y) \rangle \right| \\
&&\le C |\widehat{\tilde \K}_1^{t+\varepsilon,t,t+\varepsilon,\y}-\widehat {\bar  \K}_1^{t+\varepsilon,t,t,\x}| \textcolor{black}{\times} \varepsilon^{-1}|\T_\varepsilon ( \btheta_{t+\varepsilon,t}(\x)-\y)|^2.
\end{eqnarray*}
It remains to control the term
\begin{align}
&\widehat{\tilde \K}_1^{t+\varepsilon,t,t+\varepsilon,\y}-\widehat {\bar \K}_1^{t+\varepsilon,t,t,\x}=(\Delta_1^{t+\varepsilon,t}-\Delta_2^{t+\varepsilon,t})(\x,\y),\notag\\
&\Delta_1^{t+\varepsilon,t}(\x,\y):=\varepsilon \int_{t}^{t+\varepsilon} \!\!\!\!du \T_\varepsilon^{-1} \tilde \gR^{t+\varepsilon,\y}(t+\varepsilon,u)B \Delta a(u,t+\varepsilon)(\x,\y)
B^* \tilde \gR^{t+\varepsilon,\y}(t+\varepsilon,u)^*\T_\varepsilon^{-1},\notag\\
&\Delta a(u,t+\varepsilon)(\x,\y)=\Big(a(u,\btheta_{u,t+\varepsilon}(\y))-a(u,\btheta_{u,t}(\x))\Big),\notag\\
&\Delta_2^{t+\varepsilon,t}(\x,\y):=\varepsilon \int_{t}^{t+\varepsilon} \!\!\!\! du \T_\varepsilon^{-1} 
\Delta {\widetilde {\bar \gR}}^{t+\varepsilon,u,t,\x,\y}(t+\varepsilon,u)
B a(u,\btheta_{u,t}(\x))B^* 
\bar \gR^{t,\x}(t+\varepsilon,u)^*\T_\varepsilon^{-1}\notag\\
& -\varepsilon \int_{t}^{t+\varepsilon} du\bigg\{ \T_\varepsilon^{-1} 
\tilde  \gR^{t+\varepsilon,\y}(t+\varepsilon,u)
B a(u,\btheta_{u,t}(\x))B^*
\Big(\Delta {\widetilde {\bar \gR}}^{t+\varepsilon,u,t,\x,\y}(t+\varepsilon,u)\Big)^*\T_\varepsilon^{-1}\bigg\},\notag\\
&\Delta {\widetilde {\bar \gR}}^{t+\varepsilon,u,t,\x,\y}(t+\varepsilon,u)=\Big(\bar \gR^{t,\x}(t+\varepsilon,u) -\tilde \gR^{t+\varepsilon,\y}(t+\varepsilon,u)\Big).
\label{BIG_DECOUPAGE}
\end{align}
\textcolor{black}{From the scaling properties of the resolvent, see e.g. Lemma 6.2 in \cite{meno:17} for details, we have that:
\begin{equation}
\label{scaled_resolvent}
 \bar \gR^{t,\x}(t+\varepsilon,u)=\T_\varepsilon \widehat{\bar  \gR}_1^{t+\varepsilon,u,t,\x} \T_{\varepsilon}^{-1} ,\ \tilde \gR^{t+\varepsilon,\y}(t+\varepsilon,u)=\T_\varepsilon \textcolor{black}{\widehat{\tilde  \gR}_1^{t+\varepsilon,u,t+\varepsilon,\y}} \T_{\varepsilon}^{-1} ,
 \end{equation}
where $\widehat{\bar  \gR}_1^{t+\varepsilon,u,t,\x},\ \textcolor{black}{\widehat{\tilde  \gR}_1^{t+\varepsilon,u,t+\varepsilon,\y}}$ are non-degenerate bounded matrices, uniformly in $u\in [t,t+\varepsilon] $. Hence, from \eqref{scaled_resolvent} and the definitions in \eqref{BIG_DECOUPAGE}:}
\begin{align}
&|\Delta_1^{t+\varepsilon,t}(\x,\y)|\notag\\
\le &C\varepsilon^{-1} \int_{t}^{t+\varepsilon} |\Delta a(u,t+\varepsilon)(\x,\y)| du\le C\varepsilon^{-1} \int_{t}^{t+\varepsilon} |\btheta_{u,t}(\x)-\btheta_{u,t+\varepsilon}(\y)|^\eta du\label{CTR_DELTA_1}\\
\le  &C\varepsilon^{\eta/2-1}\int_t^{t+\varepsilon} du|\varepsilon^{1/2}\T_\varepsilon^{-1}(\btheta_{u,t}(\x)-\btheta_{u,t+\varepsilon}(\y))|^\eta \notag\\
\le& C\varepsilon^{\eta/2}(|\varepsilon^{1/2}\T_\varepsilon^{-1}(\btheta_{t+\varepsilon,t}(\x)-\y)|^\eta+1),\notag
\end{align}
using again Lemma \ref{lemme:bilipflow} for the last inequality.
\textcolor{black}{Still from \eqref{scaled_resolvent}, the definitions in \eqref{BIG_DECOUPAGE} \textcolor{black}{and recalling as well that $\T_{\varepsilon}^{-1}B a(u,\btheta_{u,t}(\x)) B^* \T_{\varepsilon}^{-1}=\varepsilon^{-2}B a(u,\btheta_{u,t}(\x)) B^*$}}, we now write  that:
\begin{eqnarray}
\label{CTR_DELTA_2}
|\Delta_2^{t+\varepsilon,t}(\x,\y)|\le C\varepsilon^{-1}\int_t^{t+\varepsilon} |\T_\varepsilon^{-1} 
\Delta {\widetilde {\bar \gR}}^{t+\varepsilon,u,t,\x,\y}(t+\varepsilon,u)\T_\varepsilon| du.
\end{eqnarray}
Note then:
\begin{eqnarray*}
&&|\T_\varepsilon^{-1} 
\Delta {\widetilde {\bar \gR}}^{t+\varepsilon,u,t,\x,\y}(t+\varepsilon,u)\T_\varepsilon|\\
&&=\Big|\T_{\varepsilon}^{-1}\int_u^{t+\varepsilon} \Big(D\gF(v,\btheta_{v,t}(\x))\bar \gR^{t,\x}(v,u)-D\gF(v,\btheta_{v,t+\varepsilon}(\y))\tilde  \gR^{t+\varepsilon,\y}(v,u)\Big) dv\T_\varepsilon\Big|\\
&&\le \int_{u}^{t+\varepsilon} |\T_{\varepsilon}^{-1}D\gF(v,\btheta_{v,t}(\x))\T_\varepsilon  | |\T_\varepsilon^{-1}\Delta {\widetilde {\bar \gR}}^{t+\varepsilon,u,t,\x,\y}(v,u)\T_\varepsilon| dv\\
&& \quad +\int_{u}^{t+\varepsilon} \Big| \T_{\varepsilon}^{-1}\Big(D\gF(v,\btheta_{v,t}(\x))-D\gF(v,\btheta_{v,t+\varepsilon}(\y))\Big)\T_\varepsilon \Big| |\T_\varepsilon^{-1}\tilde \gR^{t+\varepsilon,\y}(v,t)\T_\varepsilon|dv\\
&&\le C\int_{u}^{t+\varepsilon}  {\varepsilon}^{-1}|D\gF(v,\btheta_{v,t}(\x))-D\gF(v,\btheta_{v,t+\varepsilon}(\y))| dv,
\end{eqnarray*}
using the Gronwall lemma and the structure of the resolvent for the last inequality. 

Pay attention that we only know from \A{S} that for all $i\in \leftB 2,n\rightB, \forall \z^{i:n}=(z_i,\cdots,z_n)\in \R^{(n-i+1)d} $, $z_{i-1}\mapsto D_{x_{i-1}}F_i(z_{i-1},\z^{i:n}) $ is $C^{\eta}(\R^d,\R^d\otimes \R^d) $-H\"older continuous for $\eta>0$. We thus have to handle the above term with some care. Write with the notations of \eqref{DEF_GFTY}:
\begin{eqnarray*}
&&|\T_\varepsilon^{-1} \Delta {\widetilde {\bar \gR}}^{t+\varepsilon,u,t,\x,\y}(t+\varepsilon,u)\T_\varepsilon|\\
&\le& 
C\int_{u}^{t+\varepsilon}  {\varepsilon}^{-1}\sum_{i=2}^n\Big(  |D_{x_{i-1}}F_i(v,\btheta_{v,t}(\x))-D_{x_{i-1}}F_i^{t+\varepsilon, \y}(v,\btheta_{v,t}(\x))| \\
&&+|D_{x_{i-1}}F_i^{t+\varepsilon, \y}(v,\btheta_{v,t}(\x))-  D_{x_{i-1}}F_i(v,\btheta_{v,t+\varepsilon}(\y))| \Big)dv\\
&\le & C\int_{u}^{t+\varepsilon}  {\varepsilon}^{-1}\sum_{i=2}^n\Big( 
|D_{x_{i-1}}F_i(v,\btheta_{v,t}(\x))-D_{x_{i-1}}F_i(v,\btheta_{v,t}(\x)_{i-1},(\btheta_{v,t+\varepsilon}(\y))^{i:n})|\\
&&+ |(\btheta_{v,t}(\x)-\btheta_{v,t+\varepsilon}(\y))_{i-1}|^\eta\Big) dv=:(R_1+R_2)(t+\varepsilon,u,t,\x,\y).
\end{eqnarray*} 
We get $|R_2(t+\varepsilon,u,t,\x,\y)|\le   C\int_{u}^{t+\varepsilon}  {\varepsilon}^{-1}|\btheta_{v,t}(\x)-\btheta_{v,t+\varepsilon}(\y)|^\eta dv$ which can be handled similarly to \eqref{CTR_DELTA_1}. This yields: 
$$|R_2(t+\varepsilon,u,t,\x,\y)|\le C\varepsilon^{\eta/2}\big(\varepsilon^{1/2}|\T_{\varepsilon}^{-1}(\btheta_{t+\varepsilon,t}(\x)-\y)| +1)^\eta .$$
 On the other hand, using a reverse Taylor expansion, for positive parameters $(\delta_i)_{i\in \leftB 2,n\rightB} $ to be specified:
\begin{align*}
&|R_{\textcolor{black}{1}}(t+\varepsilon,u,t,\x,\y)|\\
\le&  C\int_{u}^{t+\varepsilon}  {\varepsilon}^{-1}\sum_{i=2}^n\Big( 
\Big|\big\{F_{i}\big(v,\btheta_{v,t}(\x)_{i-1}+\delta_i,(\btheta_{v,t}(\x))^{i:n}\big)-F_{i}(v,\btheta_{v,t}(\x))  \big\} \\
&-\big\{ F_{i}\big(v,\btheta_{v,t}(\x)_{i-1}+\delta_i,(\btheta_{v,t+\varepsilon}(\y))^{i:n}\big)-F_{i}\big(v,\btheta_{v,t}(\x)_{i-1},(\btheta_{v,t+\varepsilon}(\y))^{i:n}\big)  \big\}\Big|\delta_i^{-1} + \delta_i^\eta\Big) dv\\
\le& C\int_{u}^{t+\varepsilon}  {\varepsilon}^{-1}\sum_{i=2}^n\Big( \sum_{j=i}^n |(\btheta_{v,t}(\x)-\btheta_{v,t+\varepsilon}(\y))_j|^{\beta_{i}^j} \delta_i^{-1}+\delta_i^\eta\textcolor{black}{\Big)} \\
\le & C \bigg\{  \sum_{i=2}^n\sum_{j=i}^n  \textcolor{black}{\Big(|\varepsilon^{1/2}(\T_\varepsilon^{-1}(\btheta_{t+\varepsilon,t}(\x)-\y))|^{\beta_i^j} +1\Big)}\varepsilon^{(j-1/2)\beta_i^j}\delta_i^{-1}  +\max_{i\in \leftB 2,n\rightB}\delta_i^\eta\bigg\},
\end{align*} 
using again Lemma \ref{lemme:bilipflow} for the last inequality. For this contribution to be a remainder it therefore suffices to choose $\delta_i= \max_{j\in \leftB i,n\rightB} \varepsilon^{(j-1/2)\beta_i^j-\gamma}$, for  $\gamma>0  $ small enough. From the above computations we eventually derive that there exists $\zeta':=\zeta'(\A{A}, (\beta_i^j)_{i\in\leftB 1,n\rightB, j\in \leftB i,n\rightB}) \in (0,1) $ s.t.
\begin{eqnarray*}
|\T_\varepsilon^{-1} 
\Delta {\widetilde {\bar \gR}}^{t+\varepsilon,u,t,\x,\y}(t+\varepsilon,u)\T_\varepsilon|
&\le& C\varepsilon^{\frac \eta2 \wedge \zeta'}\big(\varepsilon^{\frac 12}|\T_{\varepsilon}^{-1}(\btheta_{t+\varepsilon,t}(\x)-\y)| +1).
\end{eqnarray*}
Plugging this bound into \eqref{CTR_DELTA_2}, we then derive from \eqref{CTR_DELTA_1} and \eqref{BIG_DECOUPAGE} that:
\begin{equation}
\label{DIFF_K}
|\widehat {\tilde \K}_1^{t+\varepsilon,t,t+\varepsilon,\y}-\widehat {\bar \K}_1^{t+\varepsilon,t,t,\x}|\le C\varepsilon^{\frac \eta2 \wedge \zeta'}\big(\varepsilon^{\frac 12}|\T_{\varepsilon}^{-1}(\btheta_{t+\varepsilon,t}(\x)-\y)| +1).
\end{equation}

\textcolor{black}{Recalling now that $\widehat {\tilde \K}_1^{t+\varepsilon,t,t+\varepsilon,\y} ,\ \widehat {\bar \K}_1^{t+\varepsilon,t,t,\x}$ are because of \eqref{GSP} and \eqref{DEF_K_MACRO} non-degenerate uniformly w.r.t. the parameter $ \varepsilon$, we deduce that the inverse matrices $\Big(\widehat {\tilde \K}_1^{t+\varepsilon,t,t+\varepsilon,\y}\Big)^{-1}, \Big(\widehat {\bar \K}_1^{t+\varepsilon,t,t,\x}\Big)^{-1} $ have the same H\"older regularity. Indeed, 
\begin{align*}
&\Big(\widehat {\tilde \K}_1^{t+\varepsilon,t,t+\varepsilon,\y}\Big)^{-1}- \Big(\widehat {\bar \K}_1^{t+\varepsilon,t,t,\x}\Big)^{-1}\\
=&\Big(\widehat {\bar \K}_1^{t+\varepsilon,t,t,\x}\Big)^{-1}\Bigg(\widehat {\bar \K}_1^{t+\varepsilon,t,t,\x}-\widehat {\tilde \K}_1^{t+\varepsilon,t,t+\varepsilon,\y}\Bigg)\Big(\widehat {\tilde \K}_1^{t+\varepsilon,t,t+\varepsilon,\y}\Big)^{-1},
\end{align*}
and we eventually conclude from \eqref{DIFF_K}. Hence}, \textcolor{black}{from the definition in \eqref{CTR_PSI_DER}}
\begin{align*}
|Q_\varepsilon|:=&\left |\langle \bigl ( (\tilde \K_{t+\varepsilon,t}^{t+\varepsilon,\y})^{-1}-(\bar \K_{t+\varepsilon,}^{t,\x})^{-1}\bigr) ( \btheta_{t+\varepsilon,t}(\x)-\y), \btheta_{t+\varepsilon,t}(\x)-\y \rangle\right|\\
=&\left| \langle \bigl(   \left(\widehat{\tilde \K}_1^{t+\varepsilon,t,t+\varepsilon,\y}\right)^{-1}- \bigl(\widehat {\bar \K}_1^{t+\varepsilon,t,t,\x}\bigl)^{-1} \bigl)(\varepsilon^{1/2}\T_\varepsilon^{-1} (\btheta_{t+\varepsilon,t}(\x)-\y))\right.\\
&,\left.\varepsilon^{1/2}\T_\varepsilon^{-1} (\btheta_{t+\varepsilon,t}(\x)-\y) \rangle \right| \\
\le & C\varepsilon^{ \eta/2\wedge \zeta'}\big(\varepsilon^{1/2}|\T_{\varepsilon}^{-1}(\btheta_{t+\varepsilon,t}(\x)-\y)| +1)
 \varepsilon |\T_\varepsilon^{-1} ( \btheta_{t+\varepsilon,t}(\x)-\y)|^2,\\
  C:=&C(\A{A},T).
\end{align*} 
We eventually get:
 \begin{align}
&|\bXi_{21}^\varepsilon(t,\x)|\notag\\
\le& C  \varepsilon^{\eta/2 \wedge \zeta'}\int_{\R^{nd}}^{}\frac{d\y}{\varepsilon^{\frac{n^2d}2}}(\varepsilon^{1/2}|\T_\varepsilon^{-1} (\btheta_{t+\varepsilon,t}(\x)-\y) |+1)
\exp(-C\varepsilon|\T_\varepsilon^{-1}(\btheta_{t+\varepsilon,t}(\x)-\y)|^2) \notag\\
\le& C\varepsilon^{ \eta/2\wedge \zeta'},\ C:=C(\A{A},T).
\label{CTR_BTH_2}
\end{align}
This yields $|\bXi_{21}^\varepsilon(t,\x)|\underset{\varepsilon \rightarrow 0}{\rightarrow} 0 $. Arguments similar to those employed for  $\bXi_{21}^\varepsilon(t,\x) $ can be used to prove that for the term $\bXi_{22}^\varepsilon(t,\x) $ defined in \eqref{DEF_BXI_22}. Namely,
\textcolor{black}{
\begin{eqnarray}
|\bXi_{22}^\varepsilon(t,\x)|&\le & \|f\|_\infty\int_{\R^{nd}}^{}\frac{d\y}{(2\pi)^{\frac{nd}2}}\Bigg| \frac{1}{\det(\tilde \K_{t+\varepsilon,t}^{t+\varepsilon,\y})^{\frac 12}}-\frac{1}{\det(\bar \K_{t+\varepsilon,t}^{t,\x})^{\frac 12}}\Bigg| (\psi_{t,\x,\y}^\varepsilon)(1)\notag\\
&\le & C\|f\|_\infty \int_{\R^{nd}}  d\y \Bigg|1-\frac{\det(\tilde \K_{t+\varepsilon,t}^{t+\varepsilon,\y})^{\frac 12}}{\det(\bar \K_{t+\varepsilon,t}^{t,\x})^{\frac 12}} \Bigg|\bar p_{C^{-1}}(t,\textcolor{black}{t+\varepsilon},\x,\y)\notag\\
&\le & C\|f\|_\infty \int_{\R^{nd}}  d\y \Bigg|1-\det\Big(\tilde \K_{t+\varepsilon,t}^{t+\varepsilon,\y}\big(\bar \K_{t+\varepsilon,t}^{t,\x}\big)^{-1}\Big)^{\frac 12} \Bigg|\bar p_{C^{-1}}(t,\textcolor{black}{t+\varepsilon},\x,\y). \label{THE_CTR_bXI_22}
\end{eqnarray}
Using again \eqref{DEF_K_MACRO}, we now write:
\begin{eqnarray*} 
 \det\Big(\tilde \K_{t+\varepsilon,t}^{t+\varepsilon,\y}\big(\bar \K_{t+\varepsilon,t}^{t,\x}\big)^{-1}\Big)&=&\det\Big(\widehat {\tilde \K}_1^{t+\varepsilon,t,t+\varepsilon,\y}\big(\widehat {\bar \K}_1^{t+\varepsilon,t,t,\x}\big)^{-1}\Big) \\
 &=&\det\Big(I+(\widehat {\tilde \K}_1^{t+\varepsilon,t,t+\varepsilon,\y}- \widehat {\bar \K}_1^{t+\varepsilon,t,t,\x})\big(\widehat {\bar \K}_1^{t+\varepsilon,t,t,\x}\big)^{-1}\Big).
 \end{eqnarray*}
}
\textcolor{black}{Plugging this identity into \eqref{THE_CTR_bXI_22}, we thus derive from \eqref{DIFF_K} that} $\bXi_{22}^\varepsilon(t,\x)\underset{\varepsilon \rightarrow 0}{\rightarrow} 0  $. The proof is complete.\qed
\subsection{Proof of Lemma \ref{convergence_LPLQ}}
{\color{black} Up to approximation argument, we can suppose without loss of generality that $f\in C_0^{1,2}([0,T)\times \R^{nd},\R)$.} Write for $p',q'>1 $,
\begin{eqnarray*}
\|f_\varepsilon-f\|_{L^{q'}([0,T],L^{p'}(\R^{nd}))}^{q'}=\int_0^T \|f_\varepsilon(t,\cdot)-f(t,\cdot)\|_{L^{p'}(\R^{nd})}^{q'}dt.
\end{eqnarray*}
\textcolor{black}{Note that, up to a middle point type argument, the indicator part in the very definition of $f_\varepsilon$ can be easily dealt. With a slight abuse of notation, we thus start from the following expression for $f_\varepsilon$:}
$$
\textcolor{black}{
\forall (t,\x)\in [0,T]\times \R^{nd},\quad \int_{\R^{nd}}f(t+\varepsilon,\y) \tilde p^{t+\varepsilon,\y}(t,t+\varepsilon,\x,\y) d\y.
}
$$

Now,
\begin{align*}
&\|f_\varepsilon(t,\cdot)-f(t,\cdot)\|_{L^{p'}(\R^{nd})}^{p'}\\
=&\int_{\R^{nd}} \Big|\int_{\R^{nd}}f(t+\varepsilon,\y) \tilde p^{t+\varepsilon,\y}(t,t+\varepsilon,\x,\y) d\y-f(t,\x)\Big|^{p'} d\x \\
\le& 2^{p'-1}\Big(\int_{\R^{nd}} \Big|\int_{\R^{nd}}f(t+\varepsilon,\y) \tilde p^{t+\varepsilon,\y}(t,t+\varepsilon,\x,\y) d\y-f(t+\varepsilon,\btheta_{t+\varepsilon,t}(\x)\textcolor{black}{)}\Big|^{p'} d\x\\
&+ \int_{\R^{nd}} |f(t+\varepsilon,\btheta_{t+\varepsilon,t}(\x))-f(t,\x)|^{\textcolor{black}{p'}}d\x\Big)=:2^{p'-1}(I_1^\varepsilon\textcolor{black}{(t)}+I_2^\varepsilon\textcolor{black}{(t)}).
\end{align*}
Recalling that $f$ is smooth and has compact support in time and space, we readily get $I_2^\varepsilon\textcolor{black}{(t)}\underset{\varepsilon\rightarrow 0}{\longrightarrow}0$. Let us now turn to $I_1^\varepsilon\textcolor{black}{(t)} $. Write:
\begin{align*}
&I_1^\varepsilon\textcolor{black}{(t)}\\
\le&  2^{p'-1}\Big(\int_{\R^{nd}}\Big|\int_{\R^{nd}}\Big(f(t+\varepsilon,\y)-f(t+\varepsilon,\btheta_{t+\varepsilon,t}(\x)\Big) \tilde p^{t+\varepsilon,\y}(t,t+\varepsilon,\x,\y) d\y\Big|^{p'}d\x\\
&+\int_{\R^{nd}}\Big|f(t+\varepsilon,\btheta_{t+\varepsilon,t}(\x)) \int_{\R^{nd}} \big(\tilde p^{t+\varepsilon,\y}(t,t+\varepsilon,\x,\y)- \bar p(t,t+\varepsilon,\x,\y )\big)d\y\Big|^{p'} d\x\Big)\\
=:&I_{11}^\varepsilon\textcolor{black}{(t)}+I_{12}^\varepsilon\textcolor{black}{(t)},
\end{align*}
where $\bar p(t,t+\varepsilon,\x,\y ):= \frac{\exp\left(-\frac 12 \left\langle (\bar \K_{t+\varepsilon,t}^{t,\x})^{-1}( \btheta_{t+\varepsilon,t}(\x)-\y), \btheta_{t+\varepsilon,t}(\x)-\y \right\rangle \right)}{(2\pi)^{\frac{nd}2}\det(\bar \K_{t+\varepsilon,t}^{t,\x})^{\frac 12}}$ is a \textit{true} density in $\y$ and already appears in the term $\bXi_3^\varepsilon $ \textcolor{black}{defined in \eqref{LAST_LABEL}} in the proof of Lemma \ref{convergence_dirac}. The previous analysis of this term readily gives $I_{12}^\varepsilon\textcolor{black}{(t)}\underset{\varepsilon\rightarrow 0}{\longrightarrow}0 $. On the other hand, from the good scaling property of equation \eqref{GSP}, Lemma \ref{lemme:bilipflow} and the H\"older inequality,  there exists $C:=C(p',\A{A})$ s.t.:
\begin{align*}
I_{11}^\varepsilon \textcolor{black}{(t)}\le& C \int_{\R^{nd}} \int_{\R^{nd}}\Big|f(t+\varepsilon,\y)-f(t+\varepsilon,\btheta_{t+\varepsilon,t}(\x))\Big|^{p'}\bar p_{C^{-1}}(t,t+\varepsilon,\x,\y) d\y d\x\\
\le & C \int_{\R^{nd}} (\I_{d(\btheta_{t+\varepsilon,t}(\x),{\rm supp}( f(t+\varepsilon,\cdot)))\le \delta}+\I_{d(\btheta_{t+\varepsilon,t}(\x),{\rm supp}( f(t+\varepsilon,\cdot)))> \delta})\\
&\times \int_{\R^{nd}}\Big|f(t+\varepsilon,\y)-f(t+\varepsilon,\btheta_{t+\varepsilon,t}(\x))\Big|^{p'}\bar p_{C^{-1}}(t,t+\varepsilon,\x,\y) d\y d\x\\
\le & C\Big(\varepsilon^{p'/2}\|Df\|_\infty+\int_{\R^{nd}} \I_{d(\btheta_{t+\varepsilon,t}(\x),{\rm supp}( f(t+\varepsilon,\cdot)))> \delta}\\
&\times \int_{\R^{nd}} 
|f(t+\varepsilon,\y)|^{p'}
\bar p_{C^{-1}}(t,t+\varepsilon,\x,\y) d\y d\x\Big)\\
\le &C\Big(\varepsilon^{p'/2}\|Df\|_\infty
+
\|f\|_\infty \exp\Big(-C^{-1}\frac{\delta^2}\varepsilon\Big)\int_{\R^{nd}}  \I_{\y\in {\rm supp}(f(t+\varepsilon,\cdot))} \\
&
\times\int_{\R^{nd}}\bar p_{C^{-1}}(t,t+\varepsilon,\x,\y)  d\x d\y\Big)\\
\le &C\Big(\varepsilon^{p'/2}\|Df\|_\infty
+
\|f\|_\infty \exp\Big(-C^{-1}\frac{\delta^2}\varepsilon\Big) \Big), 
\end{align*}
where in the above computations we recall that $\bar p_{C^{-1}} $ introduced in \eqref{CTR_GRAD} is actually a density w.r.t. $\x$. \textcolor{black}{We conclude the proof from the above convergence of $I_1^\varepsilon(t)$, $I_2^\varepsilon(t)$ thanks to the dominated convergence theorem.}

\section*{Acknowledgments.}
For the second author, the article was prepared within the framework of a subsidy granted to the HSE by the Government of the Russian Federation for the implementation of the Global Competitiveness Program.
\bibliographystyle{alpha}
\bibliography{bibli}

\end{document}